\documentclass[11pt]{elsarticle}
\usepackage{graphicx}
\usepackage{amssymb}    
\usepackage{epstopdf}

\usepackage[toc,page,header]{appendix}
\usepackage{lscape}
\usepackage{epsfig}
\usepackage{fullpage}
\usepackage{graphicx}
\usepackage{fancyhdr}
\usepackage{amsmath,amssymb,xspace}
\usepackage{tabularx}

\newfont{\bbb}{msbm10 scaled\magstep1}

\newtheorem{prop}{Proposition}[section]

\newtheorem{rem}{Remark}[section]

\newtheorem{definition}{Definition}[section]

\let \leq \leqslant
\let \geq \geqslant

\let \epsilon \varepsilon

{ \par \medskip \par
  \noindent \textit{\textbf{Demonstration\/}} : }{\null \hfill $\Box$ \par }
 
\newcommand{\R} {\ensuremath{\mathbb{R}}}
\newcommand{\Z} {\ensuremath{\mathbb{Z}}}
\newcommand{\N} {\ensuremath{\mathbb{N}}}
\newcommand{\C} {\ensuremath{\mathbb{C}}}

\renewcommand{\thefootnote}{\arabic{footnote}}

\journal{}

\begin{document}

\begin{frontmatter}

\title{Domain Decomposition Method for the N-body Time-Independent and Time-Dependent Schr\"odinger Equation}

\author[carl,crm]{E. Lorin}
\ead{elorin@math.carleton.ca}

\address[carl]{School of Mathematics and Statistics, Carleton University, Ottawa, Canada, K1S 5B6}
\address[crm]{Centre de Recherches Math\'{e}matiques, Universit\'{e} de Montr\'{e}al, Montr\'{e}al, Canada, H3T~1J4}

\renewcommand{\thefootnote}{\arabic{footnote}}

\begin{abstract}
This paper is devoted to the derivation of a pleasingly parallel Galerkin method for the time-independent $N$-body Schr\"odinger equation, and its time-dependent version modeling molecules subject to an external electric field \cite{BAN,gauge,CCT}.  In this goal, we develop a Schwarz Waveform Relaxation (SWR) Domain Decomposition Method (DDM) for the $N$-body Schr\"odinger equation. In order to optimize the efficiency and accuracy of the overall algorithm, i) we use mollifiers to regularize the singular potentials and to approximate the Schr\"odinger Hamiltonian, ii) we select appropriate orbitals, and iii) we carefully derive and approximate the SWR transmission conditions. Some low dimensional numerical experiments are presented to illustrate the methodology.
\end{abstract}

\begin{keyword} $N$-body Schr\"odinger equation, domain decomposition method, mollifiers, parallel computing.
\end{keyword}

\end{frontmatter}
\tableofcontents
\section{Introduction}\label{NBODY}
This paper is devoted to the derivation of a pleasingly parallel real-space algorithm for solving the $N$-body Schr\"odinger equation. It is well-known that the numerical computation to the solution to this equation faces the curse of dimensionality, as it requires computations in a $3N$-dimensional space. Even for 2 electrons ($N=2$), smart parallel numerical algorithms must then be developed in order to tackle this problem. At the discrete level, the dimension of the Hamiltonian is basically dependent on the chosen basis functions: the less information is contained in the basis functions, the larger the dimension of the approximate Hamiltonian. The ``worst-case scenario'', would be a finite difference approximation. However, elaborated alternatives exist such as the Full Configuration Interaction (FCI) \cite{ostlund}, which requires highly non-trivial basis functions, and which then allows for the construction of relatively compact discrete Hamiltonians. The time-independent and dependent Schr\"odinger equations are computed using algorithms (linear system and eigenvalue solvers) basically requiring many high-dimensional matrix-vector products. The parallelization of such algorithms is a major research field, and numerous parallel libraries exist among which, we can for instance cite {\tt lapack, arpack, sparselib, petsc, iml,...}. However, the performance of the parallel implementation is often far from ideal, and other approaches should be explored. In this goal, we propose a Schwarz domain decomposition method for solving the $N$-body Schr\"odinger equation. The principle consists of solving a large number of Schr\"odinger equations on ``small'' spatial domains. The interest is double. First, the parallel implementation is expected to be highly efficient, as local Schr\"odinger equations are solved independently and their corresponding solution is connected to the others only through the so-called transmission conditions \cite{DoleanBook}. Secondly, we can benefit from a {\it scaling effect}, as the computational complexity of real-space numerical (time-independent or dependent) Schr\"odinger equation solvers, is usually polynomial in time. We then solve in parallel, several small linear systems associated to local Schr\"odinger equations, rather than a large one. The price to pay though, is the need for computing several times (Schwarz iterations) the same equations each Schwarz iteration, but with different boundary conditions. More specifically, the chosen DDM is a Schwarz Waveform Relaxation method \cite{halpern3,GanderHalpernNataf,halpern2,jsc}, which is an elaborated fixed point algorithm. Notice that SWR could as well be used as a preconditioning technique, but we would a priori then not benefit from the {\it scaling effect}, as using the latter approach we will still have to solve a huge linear system, although of course very preconditioned. This method is characterized by a choice of transmission conditions (or boundary conditions) on each subdomain, derived from the solution to local wave equations.  This popular method is for instance analyzed for a 2-subdomain Schr\"odinger equation in  \cite{halpern2,jsc,lorin-TBS,lorin-TBS2,AML}.\\
A SWR-Galerkin method for solving the $N$-body Schr\"odinger equation is developed in this paper, using 2 different types of basis functions. The first basis is composed by Gaussian functions. The second basis which is used, is constituted by local Slater's determinants in the FCI formalism. The corresponding SWR method requires i) the computation of 1-electron orbitals, that is eigenfunctions of local $1$-body Schr\"odinger Hamiltonians, ii) from which we construct local Slater's determinants. These Slater determinants can be used as local basis functions. The SWR methods developed in this paper are first applied to solve the stationary Schr\"odinger equation using the Normalized Gradient Flow (NFG) method \cite{bao}. The NGF is a minimization method, which consists in solving a normalized Schr\"odinger equation in imaginary time (this is why, it is also referred in the literature as the {\it imaginary time method}), or equivalently a normalized heat equation, with variable integration times. The SWR method is next applied to the time-dependent Schr\"odinger equation, modeling a molecule subject to an external electric field. For sufficiently intense fields, the $N$-body wavefunction is expected to be delocalized, requiring then large $3N$-dimensional computational domains and then justifying the use of DDM. The purpose of this paper is not to show some high-dimensional simulations but rather to precisely describe aefficient general strategy for addressing the N-body problem. Some numerical results are however proposed in one dimension ($d=1$) for 2 electrons ($N=2$) to illustrate the proposed approach.
\subsection{$N$-body Time-Independent Schr\"odinger Equation}\label{TISE}
The stationary $N$-particle Schr\"odinger equation, under the Born-Oppenheimer approximation, reads in $d$ dimensions \cite{CAM15-10,ostlund}, as follows
\begin{eqnarray}\label{etise}
H_0\psi(\widetilde{{\bf x}}_1,\cdots,\widetilde{{\bf x}}_N) = \lambda \psi(\widetilde{{\bf x}}_1,\cdots,\widetilde{{\bf x}}_N)
\end{eqnarray}
with $\widetilde{{\bf x}}_i=({\bf x}_i,\omega_i)$,  where ${\bf x}_i \in \R^d$ is the spatial coordinates of the $i$th electron, and $\omega_i = \{-1/2,1/2\}$ its spin. In \eqref{etise}, the wavefunction $\psi$ is an eigenstate associated to the eigenvalue $\lambda.$ The Schr\"odinger Hamiltonian $H_0$, for $N$ electrons and $P$ fixed nuclei (Born-Oppenheimer approximation) is given by
\begin{eqnarray*}
H_0 = -\cfrac{1}{2}\sum_{i=1}^N\triangle_i -\sum_{i=1}^N\sum_{A=1}^P\cfrac{Z_A}{|{\bf x}_i-{\bf x}_A|} + \sum_{i=1}^N\sum_{j>i}^N\cfrac{1}{|{\bf x}_i-{\bf x}_j|}.
\end{eqnarray*}
where ${\bf x}_A \in \R^d$ denotes the position of the $A$th nucleus, and $Z_A$ its charge. In order to ensure the antisymmetry of the wavefunction $\psi$, due to Pauli's exclusion principle
\begin{eqnarray*}
\psi(\widetilde{{\bf x}}_1,\cdots,\widetilde{{\bf x}}_p,\cdots,\widetilde{{\bf x}}_q,\cdots,\widetilde{{\bf x}}_N) = -\psi(\widetilde{{\bf x}}_1,\cdots,\widetilde{{\bf x}}_q,\cdots,\widetilde{{\bf x}}_p,\cdots,\widetilde{{\bf x}}_N).
\end{eqnarray*}
for any $1\leq p < q \leq N$, one can consider the traditional FCI approach \cite{ostlund} based on Slater's determinants. Assume that $\{\phi_j({\bf x})\}_{j=1}^M$ is a set of $M$ orthonormal spatial orbitals in $\R^d$, and that $\alpha(\omega)$, $\beta(\omega)$ denote the spin coordinates, we can then define $2M$ orbitals as follows: $\chi_{2j-1}(\widetilde{{\bf x}}) = \phi_j({\bf x})\alpha(\omega)$, $\chi_{2j}(\widetilde{{\bf x}}) = \phi_j({\bf x})\beta(\omega)$, for $j=1,\cdots,M$. Next, $N$-spin orthogonal orbitals are defined among ${2M \choose N}$ Slater's (antisymmetric) determinants.
\begin{eqnarray*}
w(\widetilde{{\bf x}}_1,\cdots,\widetilde{{\bf x}}_N) = 
\cfrac{1}{\sqrt{N!}}\left|
\begin{array}{cccc}
\chi_1(\widetilde{{\bf x}}_1) & \chi_2(\widetilde{{\bf x}}_1) & \cdots & \chi_N(\widetilde{{\bf x}}_1) \\
\chi_1(\widetilde{{\bf x}}_2) & \chi_2(\widetilde{{\bf x}}_2) & \cdots & \chi_N(\widetilde{{\bf x}}_2) \\
\cdot & \cdot & \cdot & \cdot \\
\cdot & \cdot & \cdot & \cdot \\
\cdot & \cdot & \cdot & \cdot \\
\chi_1(\widetilde{{\bf x}}_N) & \chi_2(\widetilde{{\bf x}}_N) & \cdots & \chi_N(\widetilde{{\bf x}}_N)
\end{array}
\right|.
\end{eqnarray*}
 Slater's determinants $\{w_j\}_j$, can then be used as basis functions in order to compute the eigenfunctions of $H_0$. However, it is possible to construct more simple local basis functions $\{w_j\}_j$. In particular, we will also use in this paper, Gaussian basis functions, in Section \ref{GBF}. Once, the basis functions are selected, any wavefunction $\psi$ is then expanded as $\sum_{j} c_jw_j$. {\it From now on} and for the sake of the presentation, we will omit the spin variable and only consider the spatial coordinates.
\subsection{$N$-body Time-Dependent Schr\"odinger Equation}\label{TDSE}
The $N$-body time-dependent Schr\"odinger equation (TDSE), for $t \in (0,T)$ in {\it length gauge} (LG) \cite{gauge} reads
\begin{eqnarray}\label{LG}
{\tt i}\partial_t \psi({\bf x},t) = \big(H_0 + \sum_{i=1}^N{\bf x}_i\cdot {\bf E}(t) \big)\psi({\bf x},t)
\end{eqnarray}
where ${\bf E}(t)$ denotes a given external electric field, under the dipole approximation (wavelength of the electric field much larger than the spatial scale of the $N$-body system).
That is, for $1 \leq l\leq K$
\begin{eqnarray*}
\left.
\begin{array}{lcl}
{\tt i}\sum_{j=1}^K\langle w_j,w_l\rangle\dot c_j(t) & = & \sum_{j=1}^K\langle H_0 w_j,w_l\rangle c_j(t) + \sum_{j=1}^K\sum_{i=1}^N\langle {\bf x}_i w_j,w_l\rangle \cdot {\bf E}(t) c_j(t)
\end{array}
\right.
\end{eqnarray*}
which can be rewritten
\begin{eqnarray*}
{\tt i}{\bf A}\dot{\bf c}(t) = \big({\bf H}_0 + {\bf T}(t)\big){\bf c}(t)
\end{eqnarray*}
with ${\bf c}=(c_1,\cdots,c_N)^T$, and where 
\begin{eqnarray*}
T_{jl}(t) = \sum_{i=1}^N\langle  {\bf x}_i w_j,w_l\rangle \cdot {\bf E}(t)
\end{eqnarray*}
and ${\bf A}=\big(\langle w_j,w_l\rangle\big)_{1\leq j,l \leq K}$. If the basis functions are Slater's determinants, their orthogonality implies that ${\bf A}$ is the identity matrix (when considering Neumann boundary conditions). When we are interested in the interaction of $N$-electrons with an intense external field ${\bf E}(t)$, it is necessary to include a very large number of 1-electron orbitals ($M \gg 1$, $K \gg 1$), in particular ``non-local'' ones, as ionization is also expected \cite{PBC,AB12}. As a consequence, actual computations {\it in all} $\R^{dN}$ is, in principle, necessary. 
\subsection{Normalized Gradient Flow (NGF) method}\label{ITM}
In order to solve the time-independent Schr\"odinger equation, the method which is proposed is the imaginary time method, also referred in the Mathematics literature as a Normalized Gradient Flow (NGF) method. Let us rewrite the time-dependent Schr\"odinger in a compact form in $\R^{3N}$:
\begin{eqnarray}\label{LG}
{\tt i}\partial_t \psi = \big(-\cfrac{1}{2}\sum_{i=1}^N\triangle_i+V\big)\psi.
\end{eqnarray}
The NGF method for computing the ground state $\phi$ of the Schr\"odinger Hamiltonian, consists of solving the time-dependent Schr\"odinger equation in imaginary time, and normalizing the solution at each time iteration. The converged state minimizes the energy functional
$$E(\phi):=\min_{\|\chi \|_{L^{2}(\mathbb{R}^{dN})}=1}E(\chi)$$ 
defined by 
\begin{equation}\label{dennrj}
E(\chi):=\int_{ \R^{dN}}|\nabla \chi({\bf x}_1,\cdots,{\bf x}_N)|^2+V({\bf x}_1,\cdots,{\bf x}_N)|\chi({\bf x}_1,\cdots,{\bf x}_N)|^2d{\bf x}_1,\cdots,d{\bf x}_N.
\end{equation}
More specifically the ground state is constructed by solving for $({\bf x}_1,\cdots,{\bf x}_N) \in  \R^{dN}$ and $t_{n}< t <t_{n+1}$,
\begin{eqnarray}\label{e1bis}
\left\{
\begin{array}{l}
\partial_t \phi({\bf x}_1,\cdots,{\bf x}_N,t) = -\nabla_{\phi^{*}}E(\phi) \\ \hspace{0.9cm} =\triangle \phi({\bf x}_1,\cdots,{\bf x}_N,t) - V({\bf x}_1,\cdots,{\bf x}_N)\phi({\bf x}_1,\cdots,{\bf x}_N,t), \\
\displaystyle \phi({\bf x}_1,\cdots,{\bf x}_N,t_{n+1}):=\phi({\bf x}_1,\cdots,{\bf x}_N,t^{+}_{n+1})=\frac{\phi({\bf x}_1,\cdots,{\bf x}_N,t^{-}_{n+1})}{\| \phi(\cdot ,t^{-}_{n+1})\|_{L^{2}(\mathbb{R}^{dN})}},\\
\phi({\bf x}_1,\cdots,{\bf x}_N,t)=\phi_0({\bf x}), \, ({\bf x}_1,\cdots,{\bf x}_N) \in  \R^{dN}, \textrm{with $\|\phi_{0} \|_{L^{2}(\mathbb{R}^{dN})}=1$.}
\end{array}
\right.
\end{eqnarray} 
In the above system of equation, $t_{0}:=0<t_{1}<...<t_{n+1}<...$ are discrete times, $\phi_{0}$ is an initial data for the time marching algorithm discretizing the projected
gradient method
 and pointwise $\lim_{t\rightarrow t_{n}^{\pm}}\phi({\bf x}_1,\cdots,{\bf x}_N,t)=\phi({\bf x}_1,\cdots,{\bf x}_N,t_{n}^{\pm})$, see \cite{bao} for $N=1$ and $d=3$. \\
However by construction, the computed state is not antisymmetric. As a consequence a constraint must be added in order to ensure that Pauli's exclusion principle is well satisfied. We denote by $\mathcal{A}$ an antisymmetrization operator. For instance, if $N=2$, and $d=1$ we have
\begin{eqnarray*}
\mathcal{A}\phi(x_1,x_2)=
\left\{
\begin{array}{cl}
\phi(x_1,x_2) & \hbox { if }  x_1 > x_2, \\
-\phi(x_2,x_1) & \hbox { if }  x_1 \leq x_2.
\end{array}
\right.
\end{eqnarray*}
We assume that $V$ is symmetric, that is for any $i,j$ in $\{1,\cdots,N\}^2$
\begin{eqnarray*}\label{sym}
V({\bf x}_1,\cdots,{\bf x}_i,\cdots,{\bf x}_j,\cdots,{\bf x}_N) = V({\bf x}_1,\cdots,{\bf x}_j,\cdots,{\bf x}_i,\cdots,{\bf x}_N).
\end{eqnarray*}
In the general situation, $\phi$ is antisymmetrized using odd permutations. We notice that $\mathcal{A}\phi$ satisfies the heat equation
\begin{eqnarray}\label{HA}
\partial_t \mathcal{A}\phi =\triangle \mathcal{A}\phi - V({\bf x}_1,\cdots,{\bf x}_N)\mathcal{A}\phi.
\end{eqnarray}
as $\partial_t$, $\triangle$ are linear operators, and as $\mathcal{A}V=V$. We have now to show that the following algorithm is energy decreasing for any $({\bf x}_1,\cdots,{\bf x}_N) \in  \R^{dN}$ and $t_{n}< t <t_{n+1}$:
\begin{eqnarray}\label{e1ter}
\left\{
\begin{array}{l}
\partial_t \phi  =\triangle \phi({\bf x}_1,\cdots,{\bf x}_N,t) - V({\bf x}_1,\cdots,{\bf x}_N)\phi({\bf x}_1,\cdots,{\bf x}_N,t), \\
\displaystyle \phi({\bf x}_1,\cdots,{\bf x}_N,t_{n+1}):=\phi({\bf x}_1,\cdots,{\bf x}_N,t^{+}_{n+1})=\frac{\mathcal{A}\phi({\bf x}_1,\cdots,{\bf x}_N,t^{-}_{n+1})}{\|\mathcal{A}\phi(\cdot ,t^{-}_{n+1})\|_{L^{2}(\mathbb{R}^{dN})}},\\
\phi({\bf x}_1,\cdots,{\bf x}_N,0)=\phi_0({\bf x}_1,\cdots,{\bf x}_N), \textrm{with $\|\phi_{0} \|_{L^{2}(\mathbb{R}^{dN})}=1$.}
\end{array}
\right.
\end{eqnarray} 
We notice first that:
\begin{eqnarray*}
\left.
\begin{array}{lcl}
\cfrac{d}{dt}\|\mathcal{A}\phi\|_{L^2(\R^{dN})}^2 & = & 2\int_{\R^{dN}}\mathcal{A}\phi\partial_t(\mathcal{A}\phi) = 2\int_{\R^{dN}}\mathcal{A}\phi \big(\cfrac{1}{2}\triangle -V \big)\mathcal{A}\phi. \\
& =&  -2\int_{\R^{dN}}\cfrac{1}{2}|\nabla \mathcal{A}\phi|^2+V\mathcal{A}\phi^2 \leq 0.
\end{array}
\right.
\end{eqnarray*}
Then following Theorem 2.1 in \cite{bao}, the energy defined in \eqref{dennrj} satisfies
\begin{eqnarray*}
\left.
\begin{array}{lcl}
\cfrac{d}{dt}E\Big(\cfrac{\mathcal{A}\phi}{\|\mathcal{A}\phi\|}\Big) & = & \int_{\R^{dN}}\cfrac{|\nabla \phi|^2}{\|\mathcal{A}\phi\|^2_{L^2(\R^{dN})}} + \cfrac{V \mathcal{A}\phi^2}{\|\mathcal{A}\phi\|^2_{L^2(\R^{dN})}} \\
& = & 2 \int_{\R^{dN}}\cfrac{\nabla \mathcal{A}\phi \cdot \partial_t(\nabla \mathcal{A}\phi)}{2\|\mathcal{A}\phi\|^2_{L^2(\R^{dN})}} +\cfrac{V\mathcal{A}\phi\partial_t(\mathcal{A}\phi)}{\|\mathcal{A}\phi\|^2_{L^2(\R^{dN})}} \\
& & -\Big(\cfrac{d}{dt}\|\mathcal{A}\phi\|^2_{L^2(\R^{dN})}\Big)\int_{\R^{dN}}\Big(\cfrac{|\nabla \mathcal{A}\phi|^2}{2\|\mathcal{A}\phi\|^4_{L^2(\R^{dN})}} + \cfrac{V\mathcal{A}\phi^2}{\|\mathcal{A}\phi\|^4_{L^2(\R^{dN})}}\Big)\\
& = & -2\cfrac{\|\mathcal{A}\phi_t\|^2_{L^2(\R^{dN})}}{\|\mathcal{A}\phi\|^2_{L^2(\R^{dN})}}\Big(\langle \mathcal{A},\phi \mathcal{A}\phi_t\rangle\Big)^2-\|\mathcal{A}\phi\|^2_{L^2(\R^{dN})}\|\mathcal{A}\phi_t\|^2_{L^2(\R^{dN})}\Big)\\
& \leq & 0.
\end{array}
\right.
\end{eqnarray*}
We then conclude that
\begin{prop}
Assuming that $V$ is a symmetric potential, the algorithm \eqref{e1ter} is convergent to an antisymmetry state of minimal energy.
\end{prop}
 the NGF algorithm will converge to the minimum energy antisymmetric state. 
\begin{rem}
A more straightforward approach is simply to notice that if the initial data is antisymmetric ($\mathcal{A}\phi_0=\phi_0$), then the solution to the heat equation will be antisymmetric as long as the potential is symmetric. This is a simple consequence of the uniqueness of the Cauchy problem associated to \eqref{HA}. Then, as mentioned in Theorem 2.2 from \cite{bao},
\begin{eqnarray}\label{e1ter}
\left\{
\begin{array}{l}
\partial_t \phi =\triangle \phi - V({\bf x}_1,\cdots,{\bf x}_N)\phi + \mu_{\phi}\phi, \, ({\bf x}_1,\cdots,{\bf x}_N) \in  \R^{dN},\, t \geq 0, \\
\phi({\bf x}_1,\cdots,{\bf x}_N,0)=\cfrac{\phi_0}{\|\phi_0\|_{L^{2}(\mathbb{R}^{dN})}}, \, ({\bf x}_1,\cdots,{\bf x}_N) \in  \R^{dN}, \textrm{with $\|\phi_{0} \|_{L^{2}(\mathbb{R}^{dN})}=1$.}
\end{array}
\right.
\end{eqnarray} 
where $\mu_{\phi}$ is defined as 
\begin{eqnarray*}
\mu_{\phi}(t) = \cfrac{1}{\|\phi(\cdot,t)\|^2_{L^2(\R^{dN})}}\int_{\R^{dN}}\cfrac{1}{2}|\nabla\phi|^2+V({\bf x}_1,\cdots,{\bf x}_N)\phi^2.
\end{eqnarray*}
\end{rem}
\subsection{Organization of the paper}
This paper is organized as follows. In Section \ref{GBF}, we present the construction of Gaussian local basis functions. We then propose in Section \ref{1D-2E}, a methodology to construct local Slater's determinants which can be used as local basis functions. Some properties of local Slater's determinants, as well as the efficient construction of local Hamiltonians is discussed in this section as well as in \ref{APXA}. Section \ref{SWR} is devoted to the derivation and implementation of the Schwarz Waveform Relaxation algorithm for solving the $N$-body Schr\"odinger equation. Some mathematical properties of the SWR will be recalled in this section, and their computational complexity will be discussed in \ref{APXC}. Sections \ref{NumGauss} and \ref{NumSlater} are devoted to some numerical experiments for solving the time-independent and time-dependent $2$-body Schr\"odinger equations, in one dimension. More specifically, the experiments are performed using local Gaussian basis functions in Section \ref{NumGauss}, and local Slater basis functions in Section \ref{NumSlater}. We finally conclude in Section \ref{conclusion}.
\section{Local Gaussian basis functions}\label{GBF}
\begin{figure}[!ht]
\begin{center}
\hspace*{1mm}\includegraphics[height=6cm, keepaspectratio]{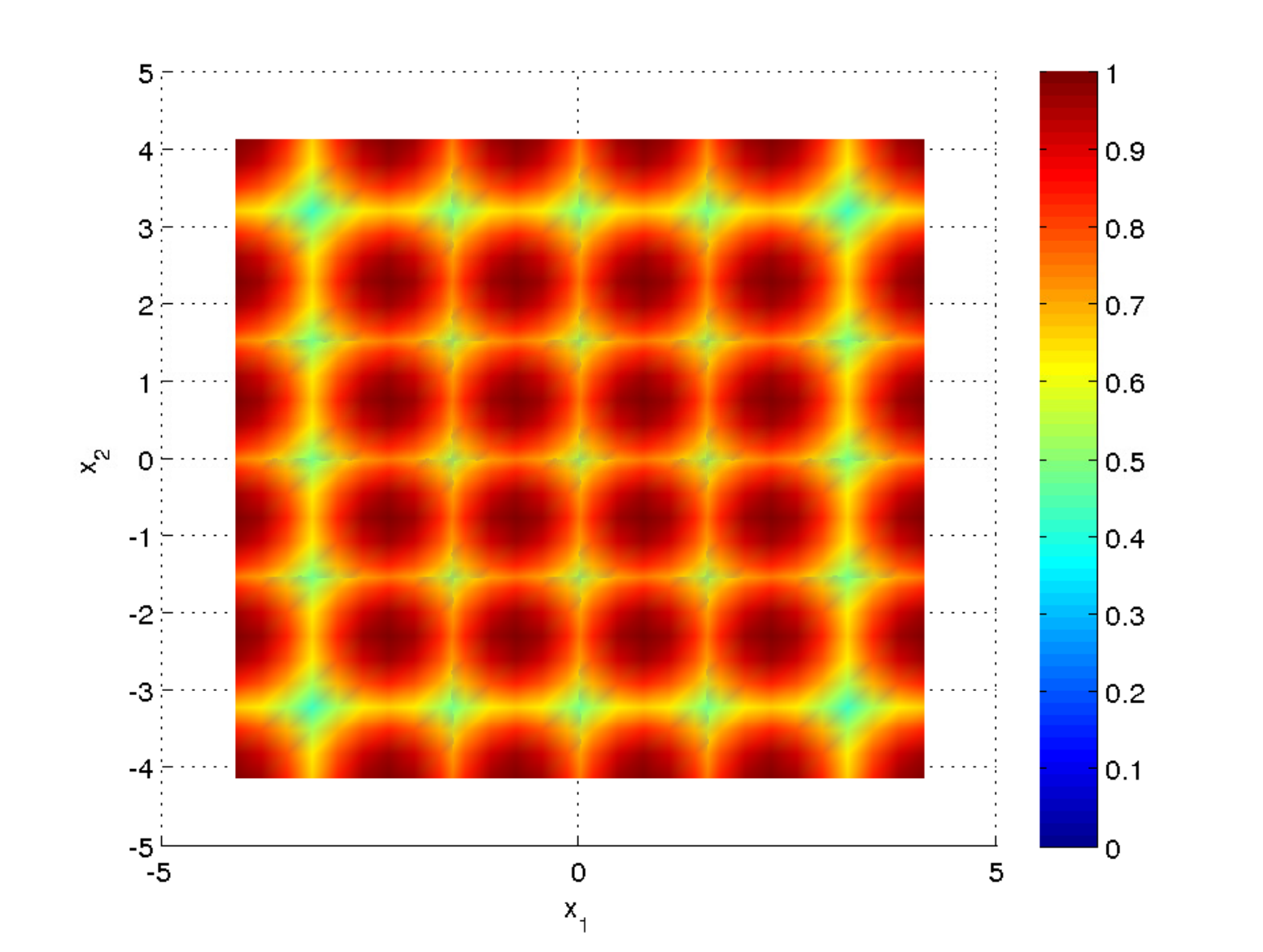}
\caption{$36$ Gaussian basis functions in one subdomain for $N=2$ and $d=1$.}
\label{GBF0}
\end{center}
\end{figure}

The domain decomposition method for solving the $N$-body Schr\"odinger equation which is proposed in this paper is based on a Galerkin method. The choice of the local basis functions is of crucial matter in order to make the computation as efficient as possible. Before considering complex basis functions in Section \ref{1D-2E}, we study the methodology with simple basis functions. A natural choice is to use Gaussian functions. \\
In order to simplify the notations, we will consider here, the case $N=2$, $d=1$. The extension of the following ideas is straightforward for arbitrary $N$ and $d$ and is shortly discussed at the end of this subsection. We denote by $\big\{D_j\big\}_{j\in\Z}$ an infinite sequence of open intervals, such that: $\R=\cup_{j\in \Z}\overline{D}_j$ and $D_i\cap D_j=\emptyset$, for $i\neq j$, and $\Lambda_{i,j}=D_i\times D_j \subsetneq \R^2$, for any $i$ and $j$ in $\Z$. Naturally we have $\cup_{(i,j)\in \Z^2}\overline{\Lambda}_{i,j} = \R^2$. We denote by $\big\{\phi^{i}_j\big\}_{(i,j)\in \N\times\Z}$ the set of one-dimensional Gaussian functions, defined by
\begin{eqnarray}\label{gauss1d}
\phi_j^{i}(x_k) = \exp\big(-\delta_k^{(i)}(x_k-\alpha^{(i)}_j)^2\big).
\end{eqnarray}
where $\delta_k^{(i)}$ is a subdomain-dependent ($i$-index) positive number for Electron $k$ ($k=1,2$), and $\alpha_{j}^{(i)} \in D_i$ is a sequence of Gaussian centers.  When the $\delta_j^{(i)}$'s are subdomain and particle independent, we will use the notation $\delta$. Now, we can construct local basis functions for any $\Lambda_{i,j}$. From any localized orbitals $\phi^i_l$, $\phi_p^j$, with $p,l$ in $\N$ (basis function indices) and $i,j$ in $\Z$ (subdomain indices), we define $v^{i,j}_{l,p}$ by:
\begin{eqnarray*}
v^{i,j}_{l,p}(x_1,x_2) = \phi^i_l(x_1)\phi^j_p(x_2).
\end{eqnarray*}
In term of support, we have
\begin{eqnarray*}
\left.
\begin{array}{lcl} 
\mbox{Supp}_{(x_1,x_2)} v^{i,j}_{l,p} & = &\mbox{Supp}\big(\phi^i_l(x_1)\phi_p^j(x_2)\big).\\
&=& \mbox{Supp}\phi_l^i\times\mbox{Supp}\phi_p^j.\\
&\subsetneq &\cup_{k=-1}^1\big(\Lambda_{i+k,j}\cup\Lambda_{i,j+k}\big).
\end{array}
\right.
\end{eqnarray*}
If $\delta^{(i)}=\delta$ and $\alpha_j^{(i)}=\alpha^{(i)}$ is taken subdomain independent, the local basis functions are actually identical in all the subdomains, which is quite convenient from a computational point of view, as we only need to construct once for all, a unique free-particle Hamiltonian. The weakness of this approach is that naturally, as the local basis functions do not contain any particular information, a large number should be used. In Fig. \ref{GBF0}, we present in a given subdomain, the local Gaussian basis functions. The construction to Gaussian basis functions for $N$ particles in $d$ dimensions is naturally straightforward by considering the tensor products of $N$ local Gaussian functions: $\Pi_{k=1}^N\phi({\bf x}_k)$. The analysis of convergence of the Galerkin method applied to the Schr\"odinger equation, and using Gaussian basis functions was presented in \cite{GalerSchro}.\\
In order to directly construct antisymmetric basis functions (at least locally) it is possible to construct (spinless) Slater-like Gaussian basis functions \cite{ostlund}, from any localized Gaussian functions $\phi^i_l$, $\phi_p^j$, with $p,l$ in $\N$ and $i,j$ in $\Z$,
\begin{eqnarray}\label{LSD}
v_{l,p}^{i,j}(x_1,x_2) = 
\cfrac{1}{\sqrt{2}}\left|
\begin{array}{cc}
\phi^i_l(x_1) & \phi^j_p(x_1) \\
\phi^i_l(x_2) & \phi^j_p(x_2) 
\end{array}
\right| = \cfrac{1}{\sqrt{2}}\Big(\phi_l^i(x_1)\phi_p^j(x_2) -\phi_l^i(x_2)\phi_p^j(x_1)\Big)
\end{eqnarray}
In the next section, we consider more elaborated antisymmetric basis functions using the traditional Slater's determinants computed from 1-electron orbitals.
\section{Local orbitals and Local Slater's Determinants as basis functions}\label{1D-2E}
This section is devoted to the construction of Local Orbitals (LO's) and Local Slater's Determinants (LSD's).\\
As we have done in the previous section, we will detail the case $d=1$ and $N=2$, that is a two-body problem in one-dimension. This is motivated by the fact that the extension to the general case (arbitrarily $N$ and $d$ case) can be deduced from \cite{CAM15-09,CAM15-10} and does not present any fundamental difficulty, but would complexify the notations.  The material presented here will be used for the Schwarz Waveform Relaxation (SWR) Domain Decomposition Method (SWR-DDM) presented in Section \ref{SWR}. The local orbitals and Slater's determinants will allow for the construction of local Hamiltonians and local wavefunctions, from which we will reconstruct a global wavefunction. The basic idea is to construct local, in the sense subdomain dependent, Slater's determinants from local 1-electron orbitals. This procedure can be applied to any subdomain, or only in some of the subdomains, typically those containing the nuclei. \\
\\
We denote by $x_i$ ($i=1,2$) the coordinate of the $i$th particle. The Schr\"odinger Hamiltonian reads, for $2$ fixed nuclei
\begin{eqnarray*}
H_0 = -\cfrac{1}{2}\sum_{i=1}^2\triangle_i -\sum_{i=1}^2\sum_{A=1}^2\cfrac{Z_A}{|x_i-x_A|} + \cfrac{1}{|x_1-x_2|}
\end{eqnarray*}
where $x_A \in \R$ denotes the position of the $A$th nucleus and $Z_A$ its charge. Antisymmetry of the wavefunction reads
\begin{eqnarray*}
\psi(x_1,x_2) = -\psi(x_1,x_2), \qquad \forall (x_1,x_2) \in \R^2.
\end{eqnarray*}
\subsection{Local FCI procedure}
We denote by $\big\{\phi_j\big\}_{j\in \Z}$ the set of 1-electron orbitals, which will allow for the construction of the compact support {\it localized orbitals}, LO's, denoted by $\big\{\phi_j^{i}\big\}_{(i,j)\in \N\times \Z}$\footnote{Top index refers to subdomain $D_i$, and bottom index to full orbital $\phi_j$ index}. Typically $\phi_j^i$ should satisfy, for any $i \in \N$ and $j \in \Z$
\begin{eqnarray}\label{phi_ij}
\left.
\begin{array}{l}
\phi^i_{j}(x) = \phi_j(x), \qquad \mbox{ if } x \in D_i\\
\mbox{Supp}\phi_j^i \subsetneq D_{i-1}\cup D_i\cup D_{i+1}\\
\phi_j^i \in C^2\big(D_{i-1}\cup D_i\cup D_{i+1},\R\big).
\end{array}
\right.
\end{eqnarray}
By construction, we will assume that $x_A,x_B \in D_0$. In order to solve the stationary Schr\"odinger equation, we choose the FCI model for a $2$-electron problem. The latter is based on (spinless) Slater's Determinant basis functions (SD's) which are defined as follows. 
From any localized orbitals $\phi^i_l$, $\phi_p^j$, with $p,l$ in $\N$ and $i,j$ in $\Z$, Slater's determinants as follows:
\begin{eqnarray}\label{LSD}
v_{l,p}^{i,j}(x_1,x_2) = 
\cfrac{1}{\sqrt{2}}\left|
\begin{array}{cc}
\phi^i_l(x_1) & \phi^j_p(x_1) \\
\phi^i_l(x_2) & \phi^j_p(x_2) 
\end{array}
\right| = \cfrac{1}{\sqrt{2}}\Big(\phi_l^i(x_1)\phi_p^j(x_2) -\phi_l^i(x_2)\phi_p^j(x_1)\Big).
\end{eqnarray}
Notice that in practice the number of determinants to compute can be reduced. For instance, for $l=p$ only indices $j\geq i+1$ should be considered. In the following, we will denote by $\mbox{Supp}_{(x_1,x_2)}$, the support of any function with respect to its $(x_1,x_2)$-variables. As a consequence:
\begin{eqnarray}\label{supp}
\left.
\begin{array}{lcl} 
\mbox{Supp}_{(x_1,x_2)} v^{i,j}_{p,l} & = &\mbox{Supp}\Big(\phi_l^i(x_1)\phi_p^j(x_2) -\phi_l^i(x_2)\phi_p^j(x_1)\Big)\\
&=& \mbox{Supp}\phi_p^i\times\mbox{Supp}\phi_p^j\\
&\subsetneq &\cup_{k=-1}^1\Big(\Lambda_{i+k,j}\cup\Lambda_{i,j+k}\Big).
\end{array}
\right.
\end{eqnarray}
In other words, the support of any $v^{i,j}_{p,l}$ is compact and is strictly included in the union of $\overline{\Lambda}_{i,j}$ with the subdomains having an edge in common with $\overline{\Lambda}_{i,j}$. By construct, $v^{i,j}_{p,l}$ is naturally antisymmetric.
\subsection{Local orbital construction}\label{subsec:SLO}
\begin{figure}[!ht]
\begin{center}
\hspace*{1mm}\includegraphics[height=6cm, keepaspectratio]{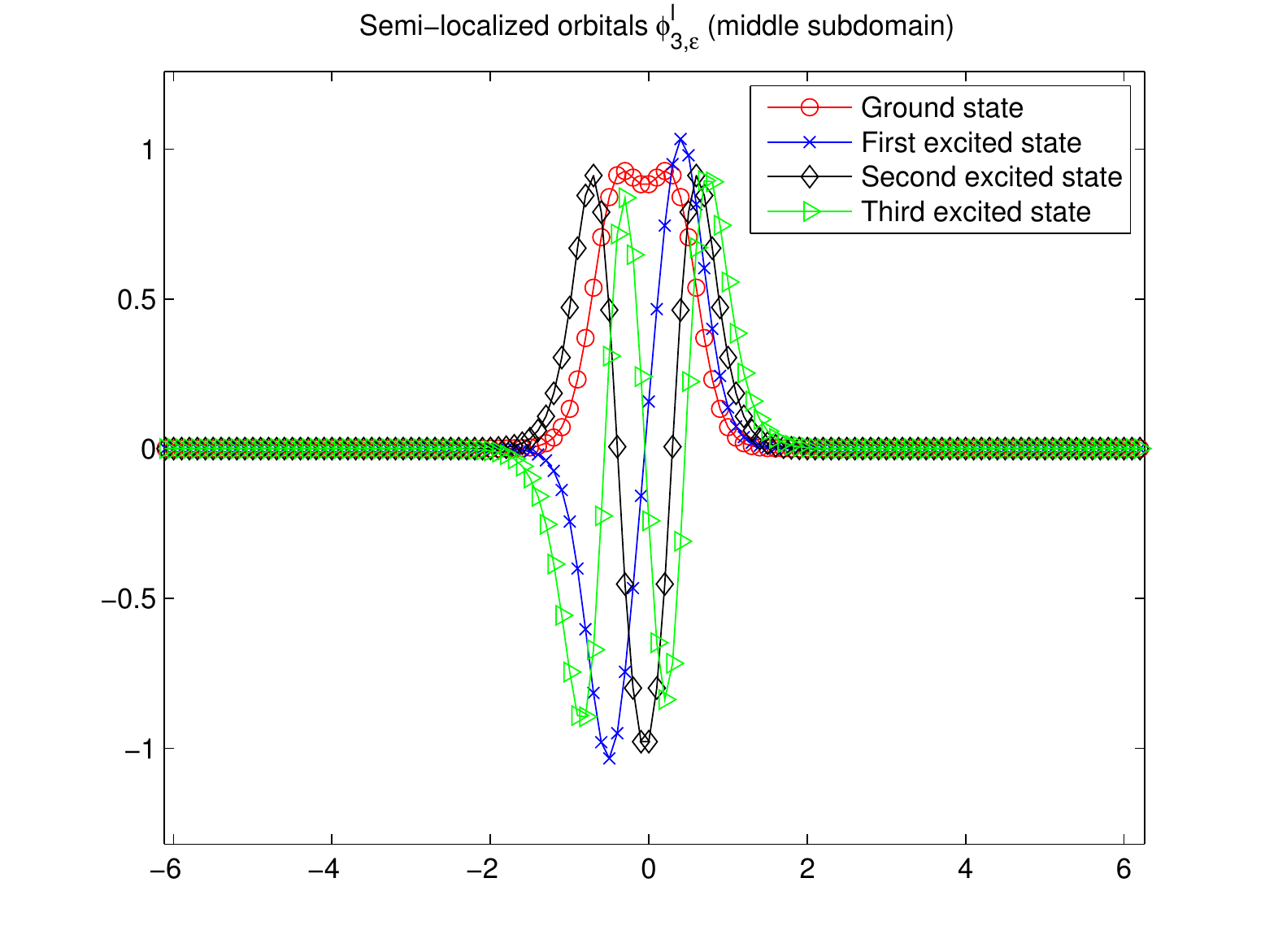}
\hspace*{1mm}\includegraphics[height=6cm, keepaspectratio]{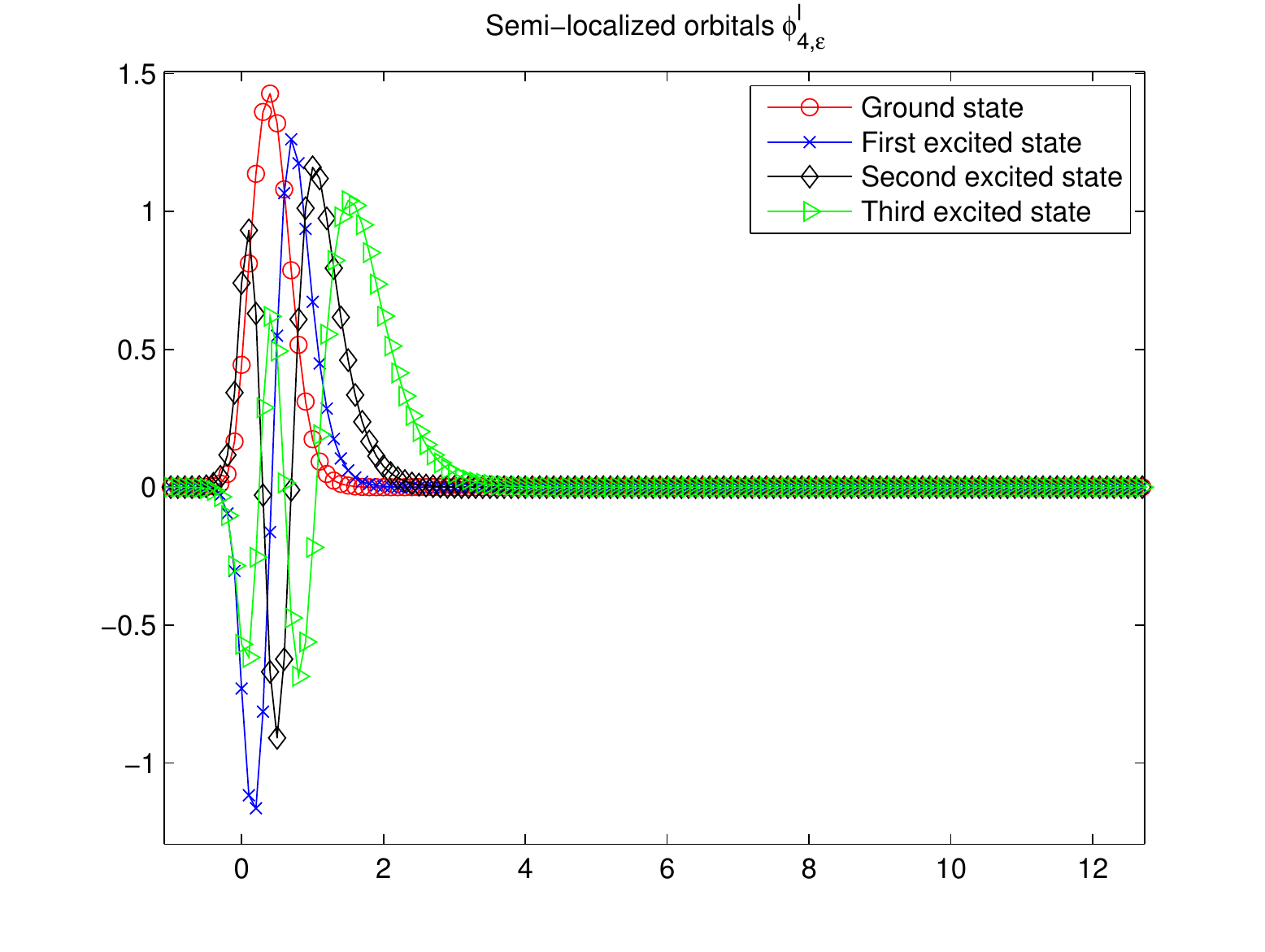}
\caption{First 4 eigenstates LO's $\phi_{l,\epsilon}^3$ (left), $\phi_{l,\epsilon}^4$ (right) for $l=1,\cdots,4$.}
\label{SLO_app2}
\end{center}
\end{figure}

The domain decomposition introduced above, allows for an adaptive selection of 1-electron orbitals per-subdomain. The key points are i) the number $P$ of nuclei, ii) their location, and iii) in the time-dependent case, the intensity of the external electric field. Notice that for any subdomain $D_i$, we select $M_i$ 1-electron localized orbitals, $\big\{\phi_l^i\big\}_{l=1}^{M_i}$. Then, from two sets of localized orbitals, $\big\{\phi_l^i\big\}_{l=1}^{M_i}$, $\big\{\phi_p^j\big\}_{p=1}^{M_j}$, we can construct ${M_i+M_j \choose 2}$ LSD's \eqref{LSD}. From a practical point of view, we consider a finite number $L$, of one-dimensional subdomains partially covering $\R$: $\cup_{i=1}^LD_{i} \subsetneq \R$. Notice that this will force us to impose absorbing conditions at the global computational domain boundary \cite{ABC,MOLPHYS}. Then, for each subdomain $\Lambda_{i,j}$, we will select $K_{i,j}$ LSD's $\big\{v_{k}^{i,j}\big\}_{k=1}^{K_{i,j}}$, among ${M_i+M_j \choose 2}$ determinants. Notice however that the procedure which is presented below, may only be relevant for subdomains containing at least one nucleus. In the other subdomains, local Gaussian basis functions could be considered. The stationary wavefunction $\psi$, solution to the Schr\"odinger equation, will then be searched in each $\Omega_i$, in the form 
\begin{eqnarray*}
\psi_i(x_1,x_2) = \sum_{k=1}^{K_{i,j}}c_{k}^{i,j}v_{k}^{i,j}(x_1,x_2).
\end{eqnarray*}
where $c_{k}^{i,j}$ are the unknown coefficients. 
We now detail the procedure to construct the localized orbitals $\big\{\phi_j^{i}\big\}_{j \in \N}$ under the condition \eqref{phi_ij},  for $i \in \{1,\cdots,L\}$. We consider as a generic example the case of the $H_2$-molecule, corresponding to $Z_A=Z_B=1$.\\
\\
\\
The approach which is proposed is based on ideas presented in \cite{CAM15-09}. Rather than post-processing the full domain 1-electron orbitals, we directly construct the smooth localized orbitals with compact support, and with orthogonality properties. This is possible thanks to the use of i) infinite potentials at the subdomain boundary, and ii) of mollifiers when a subdomain contains a nucleus singularity. We proceed as follows.  We consider the two following situations, for a given subdomain $D_i$, with $2\leq i \leq  L-1$.
\begin{itemize}
\item {\it $D_i$ contains a nucleus singularity}. Only a few subdomains belong to this first category, in particular when we are interested in the time-dependent Schr\"odinger equation for intense field-particle interaction. In that case, mollifiers will allow for an arbitrarily accurate smoothing of the nucleus singularities. Notice that in 1-d, the singularity treatment is different than in 3-d. Indeed in the latter case, we benefit from the fact that a Coulomb potential, up to a multiplicative constant is a fundamental solution to Poisson's equation. This property allows for an accurate and efficient treatment of the localized orbitals. The Coulomb potential is then approximated by a smooth function $G_{\epsilon}$, thanks to mollifiers $B_{\epsilon}$ as defined in \cite{CAM14-43}, and such that:
\begin{eqnarray}\label{GEPS}
G_{\epsilon}({\bf x}) = \cfrac{1}{4\pi}\big(V*B_{\epsilon}\big)({\bf x})
\end{eqnarray}
where $V({\bf x}) = -1/|{\bf x}|$, which also satisfies
\begin{eqnarray}\label{PE1D}
4\pi \triangle G_{\epsilon}({\bf x}) = B_{\epsilon}({\bf x}).
\end{eqnarray}
As a consequence, a smooth approximation of the Coulomb potential $V$ using \eqref{PE1D}, can be constructed with $G_{\epsilon}\rightarrow_{\epsilon \rightarrow 0}V$ in $\mathcal{D}'(\R^3)$. Notice that this property is also fundamental for efficiently computing the $6$-dimensional integrals in order to construct the global discrete Hamiltonian \cite{CAM14-43}. In 1-d, the fundamental solution of the Poisson equation is $|x|$ and the latter property does not occur anymore. Instead, we directly computed $G_{\epsilon}$ using \eqref{GEPS} with $B_{\epsilon}$ defined by:
\begin{eqnarray}
\label{Beps}
B_{\epsilon}(x)=\left\{
\begin{array}{cc}
\cfrac{1}{\epsilon}\sigma_{M(1)}\Big(1-\big(\cfrac{x}{\epsilon}\big)^2\Big)^M, & |x|\leq \epsilon,\\
0, & |x|>\epsilon
\end{array}
\right.
\end{eqnarray}
where $M$ refers to the order of the mollifier and the scaling factors $\sigma_M(1)$ are explicitly defined in \cite{CAM14-43}. For instance, for $\sigma(1)=3/4$, $\sigma(2)=15/16$ we represent in Fig. \ref{beps05} for a unique domain $B_{\epsilon=0.5}$ (Left) and $G_{\epsilon=0.5}$ (Right). In particular it is proven in \cite{CAM14-43}, that for any smooth function $f$
\begin{eqnarray*}
\|f-f*B_{\epsilon}\|_2 = \sum_{k \geq 1}c_{k}\epsilon^{2k}
\end{eqnarray*}
for some positive sequence $\{c_{k}\}_k$.
\begin{figure}[!ht]
\begin{center}
\hspace*{1mm}\includegraphics[height=6cm, keepaspectratio]{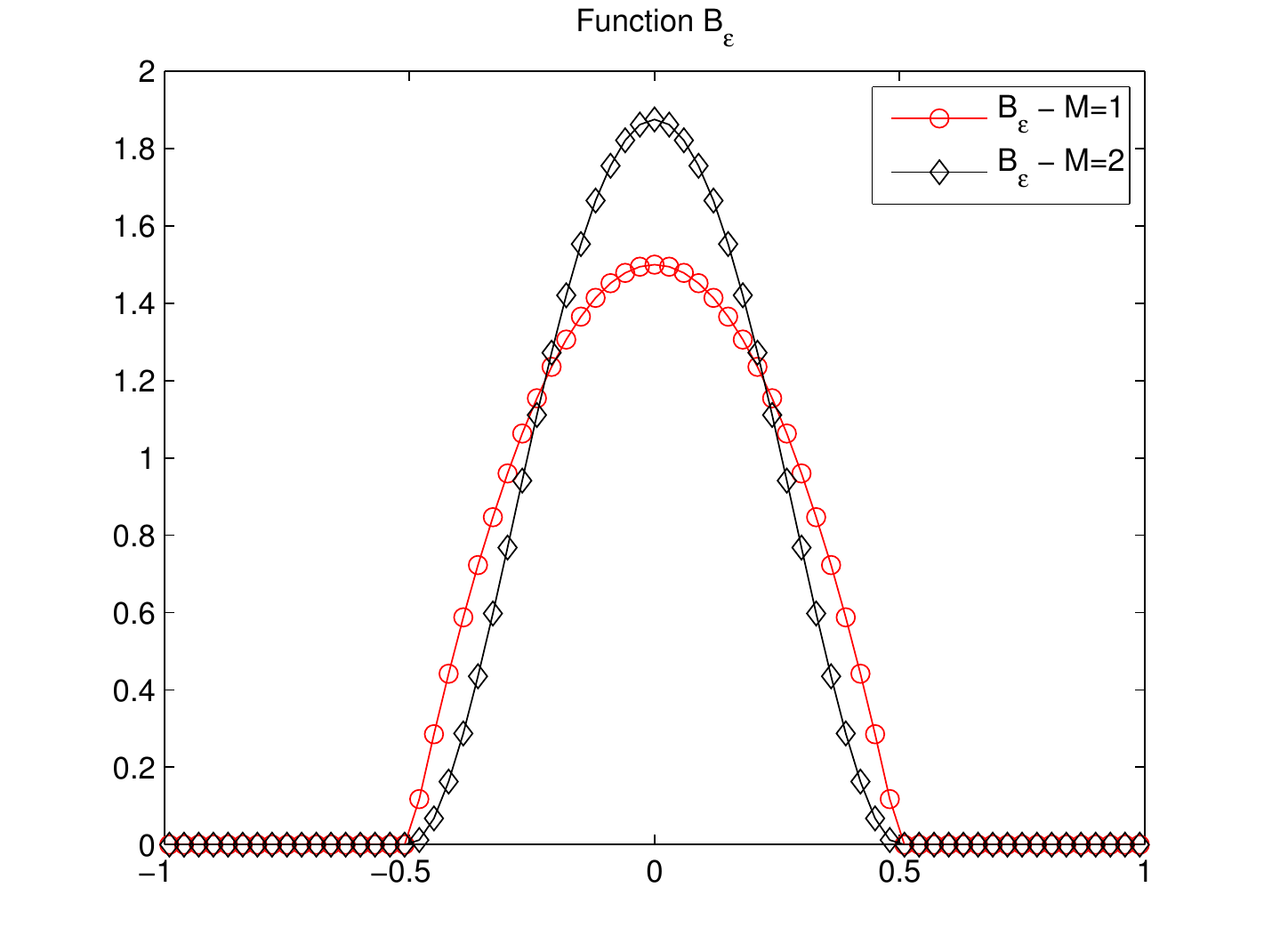}
\hspace*{1mm}\includegraphics[height=6cm, keepaspectratio]{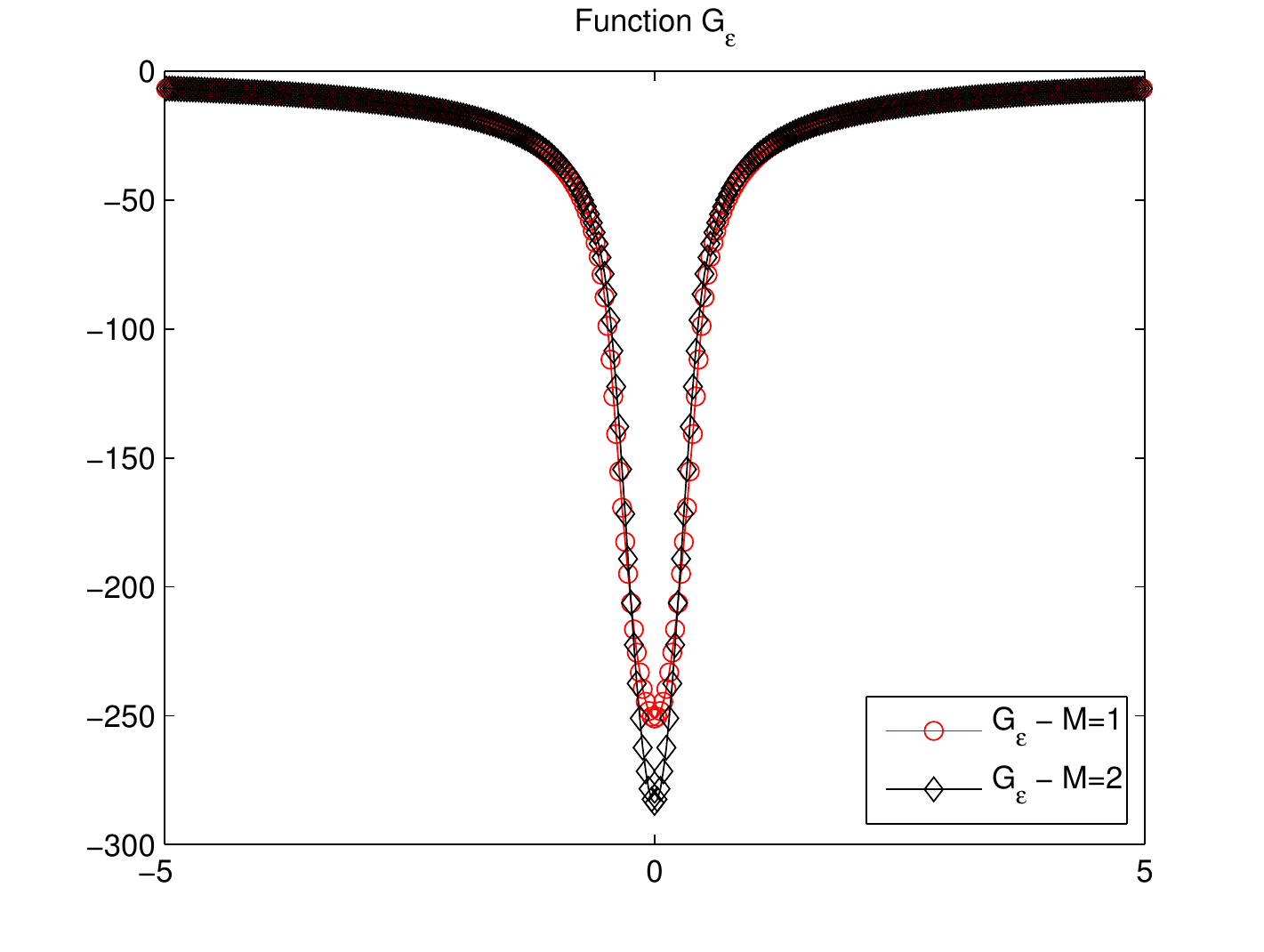}
\caption{(Left) Mollifiers $B_{0.5}$ ($M=1,2$), (Right) and $G_{0.5}$ for $d=1$.}
\label{beps05}
\end{center}
\end{figure}
Once $G_{\epsilon}$ is computed, we introduce a barrier potential as in \cite{CAM15-09}
\begin{eqnarray}\label{BARPOT}
V_{\textrm{b}}(x) = s_{\epsilon_b}(x-x_b)V_{\infty}
\end{eqnarray}
where i) the smooth function $s_{\epsilon_s}$ is equal to $0$ for $x<x_b-\epsilon_b/2$ and $1$ for $x>x_b-\epsilon_b/2$, ii) $\epsilon_b>0$, and iii) $V_{\infty}$ and $x_b$ are imposed. The support  of the localized orbitals is then $(x_{c_i}-x_b,x_{c_i}+x_b)$, where $x_{c_i}$ denotes the coordinates of the center of the subdomain $D_i$.  We typically choose $x_b > |D_i|/2$ to ensure that the localized orbitals are not null at $D_i$'s boundary. A contrario, taking $x_b$ too large will lead to a loss of computational efficiency due to a large localized orbital support. In $D_i$, we then solve the following one-dimensional one-electron eigenvalue problem 
\begin{eqnarray*}
\big(-\cfrac{1}{2}\partial_{x}\big(a_{\epsilon}(x)\partial_x\big)  + G_{\epsilon}(x-x_A) + G_{\epsilon}(x-x_B) + V_{\textrm{b}}(x-x_{c_i})\big)\varphi_{l}^{i}(x) = \lambda^{i}_{l,\epsilon}\varphi_{l}^{i}(x)
\end{eqnarray*}
where $1\leq l\leq M_i$ (resp. $2 \leq i\leq L-1$) is the orbital (resp. subdomain) index and where $a_{\epsilon}(x):=1-s_{\epsilon_b}(x-x_b)$.
Notice that the choice of the localized orbitals is motivated by physical considerations.  When we are interested in field-particle interaction, for subdomains containing the nuclei, we will select the localized orbitals corresponding to the lower energy states, as they will be predominant in the overall wavefunction in the vicinity of the nucleus singularities.
\item {\it $D_i$ does not contain any nucleus singularity}. In that case, the regularization  of the Coulomb potential through mollifiers is naturally useless. The localized orbitals are then directly obtained by solving
\begin{eqnarray*}
\big(-\cfrac{1}{2}\partial_{x}\big(a_{\epsilon}(x)\partial_x\big)  - \cfrac{1}{|x-x_A|} - \cfrac{1}{|x-x_B|}+ V_{\textrm{b}}(x)\big)\varphi_{l}^{i}(x) = \lambda^{i}_{l,\epsilon}\varphi_{l}^{i}(x)
\end{eqnarray*}
Similarly to the previous case (subdomain containing the nuclei), the selected localized orbitals will strongly depend on the relative position of the nuclei $/$ $D_i$. Alternatively, for those subdomains, we can use local Gaussian basis functions as described in Section \ref{GBF}.
\end{itemize}
Once the localized orbitals are computed, we can construct the discrete Schr\"odinger Hamiltonian. The approach which is proposed below will benefit from i) the fact that the localized orbitals are selected accordingly to the position of the nuclei, ii) the compact support of the LO's and iii) their orthogonality property (more specifically their extension by $0$ to all $\Omega$). This last point necessitates some precisions. First, we notice that by construction for any $i \in \{2,\cdots,L-1\}$, the supports of $\{\varphi_l^{i}\big\}_{1\leq l\leq M_i}$ and of $\{\varphi_m^{j}\big\}_{1\leq m\leq M_{j}}$ with $j \neq i-1,i,i+1$ are disjoint, so that these LO's are trivially orthogonal. By construction, the LO's $\{\varphi_l^{i}\big\}_l$ of any $D_i$ are also orthogonal to each other. For $j=i-1$ or  $j=i+1$, the orthogonality of the LO's $\{\varphi_l^{i}\big\}_{1\leq l\leq M_i}$ and of $\{\varphi_m^{j}\big\}_{1\leq m\leq M_j}$ is not, a priori, ensured. However, by construction for any $1\leq l \leq M_i$ and $1\leq m \leq M_{i\pm 1}$
\begin{eqnarray*}
\Big|\hbox{Supp}\big(\varphi_l^i\cap\varphi^{i\pm 1}_m\big)\Big| \leq 2x_b-\big|D_i \cup D_{i\pm 1}\big| \, .
\end{eqnarray*}
Then, as the LO's (smoothly) vanish at the boundary of their support, for $x_b-|D_{i,i\pm 1}|/2$ small enough, we expect that $\int_{\R}\varphi_l^{i}(x)\varphi^{i\pm1}_m(x)dx$ to be small. For $L=2$ and $x_b=8$ (which is relatively very large) and ($L=$) 5 subdomains, we represent for subdomain $D_{i=0}=(-10,10)$ (resp. $D_{i=1}=(-2,18)$) $\varphi_l^0$ (resp. $\varphi^1_l$), for $l=1,\cdots,4$ in Figs. \ref{SLO_app2}.
\subsection{Augmented local bases}\label{subsec:CB}
\begin{figure}[!ht]
\begin{center}
\hspace*{1mm}\includegraphics[height=4cm, keepaspectratio]{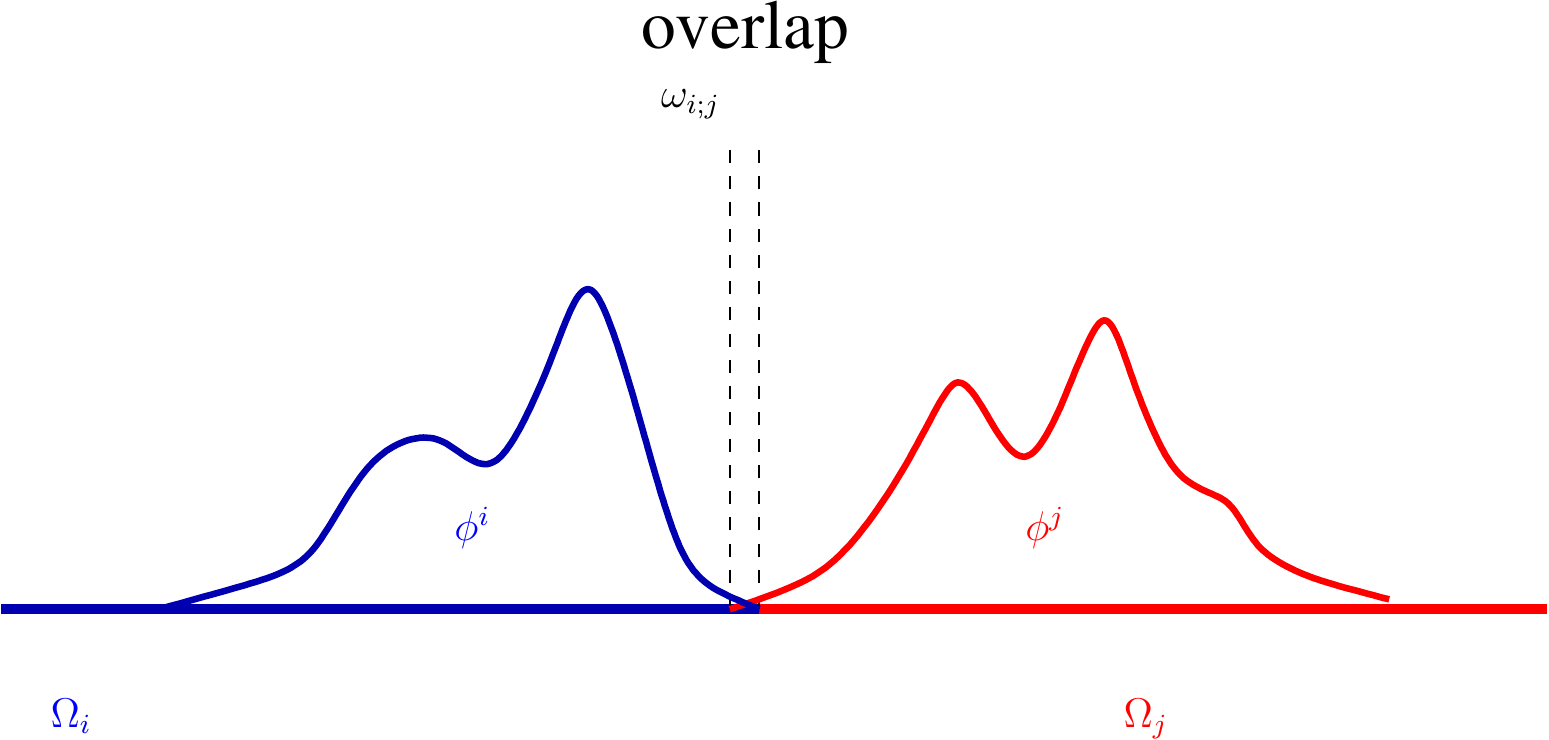}
\hspace*{1mm}\includegraphics[height=4cm, keepaspectratio]{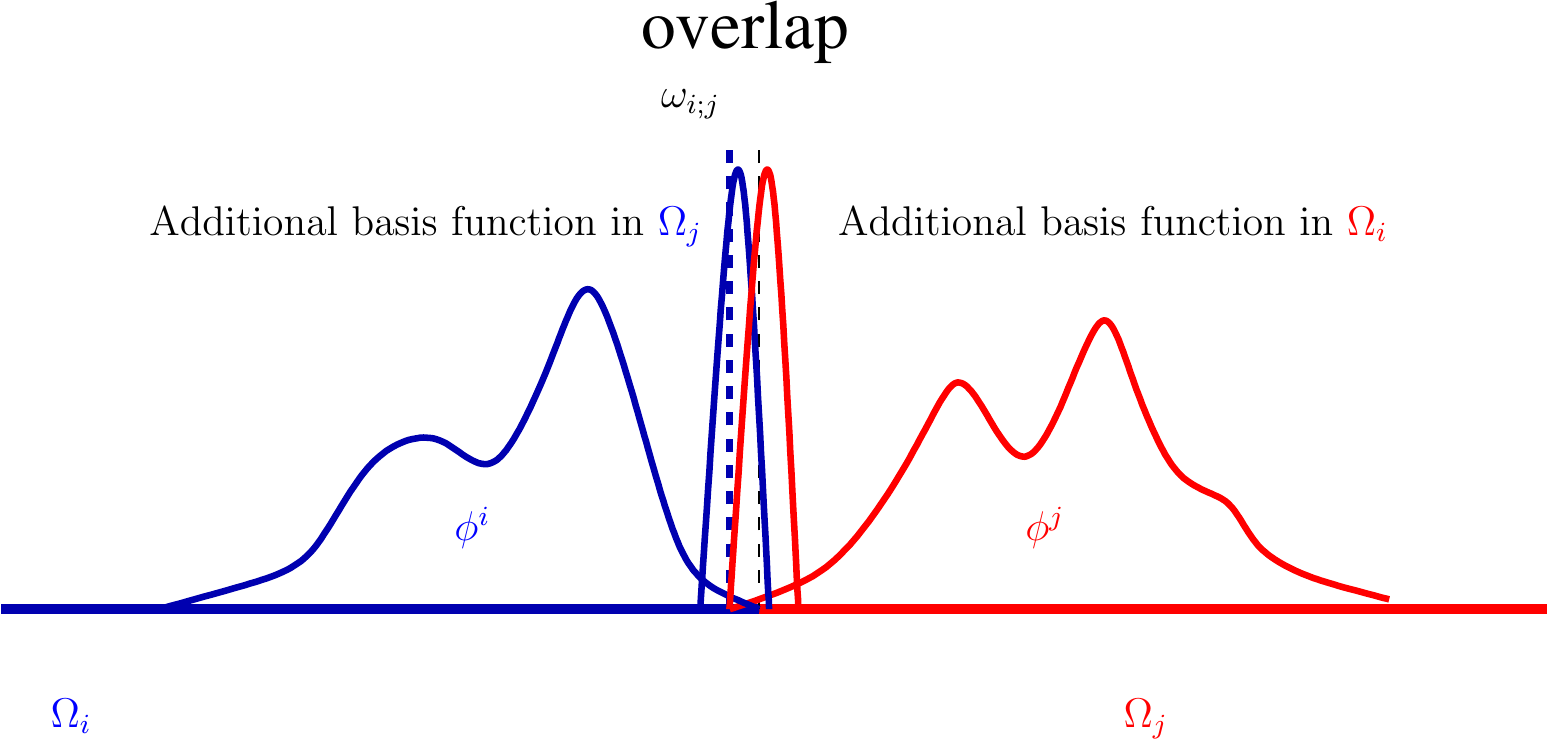}
\caption{(Left) Local basis functions overlapping issue. (Right) Additional Gaussian basis functions ensuring a proper transmission.}
\label{overPB}
\end{center}
\end{figure}

By construction the local Slater's determinant basis functions are null at the boundary of the subdomains. This can constitute an issue if the overlapping zone between two subdomains is too narrow, as in such zones the basis functions are basically null or very small, see Fig. \ref{overPB} (Left). In order to fix this issue, a simple solution consists of adding Gaussian basis functions all around the subdomains Fig. \ref{overPB} (Right). It will then ensure that in any overlapping zone the local wavefunctions could be properly transmitted from one subdomain to another thanks to the transmission conditions.
Notice that in \ref{APXA}, we present a general strategy to efficiently compute the integrals involved in the construction of the local Hamiltonians, using the formalism proposed in \cite{CAM15-10}.
\subsection{Important remarks about subdomain and local basis function indices}\label{notations}
In order to lighten the presentation, some compact notations will be used along the paper. 
\begin{itemize}
\item Functions $\{w_l\}_l$ will systematically refer to basis functions for a unique domain problem, that is without DDM.
\item  In a given subdomain $\Omega_i$ ($i$-index), the local basis functions ($l$-index) can be also denoted by $v_l^i$. For a given two-dimendional subdomain $\Lambda_{i,j}$, the local basis functions could be also denoted by $v^{i,j}_{l}$, where $l$ denotes the basis function index, and $i,j$ the one-dimensional subdomain indices. This notation was already used in Subsection \ref{subsec:SLO} to define local Slater's determinants.
\item  For and $N=2$, $d=1$, $1\leq i\leq L$ and $1\leq j \leq L$, $\Lambda_{i,j}$ will also be denoted by $\Omega_{i+jL}$. In this case, the $L^2$ subdomains will be reindexed as $\{\Omega_i\}_{1\leq i \leq L^2}$.
\end{itemize}
In general, for the local basis functions the top index always refers to the basis function index, and the bottom one to the subdomain index.
\section{Schwarz waveform relaxation domain decomposition method for the Schr\"odinger equation}\label{SWR}
\begin{figure}[!ht]
\begin{center}
\hspace*{1mm}\includegraphics[height=6cm, keepaspectratio]{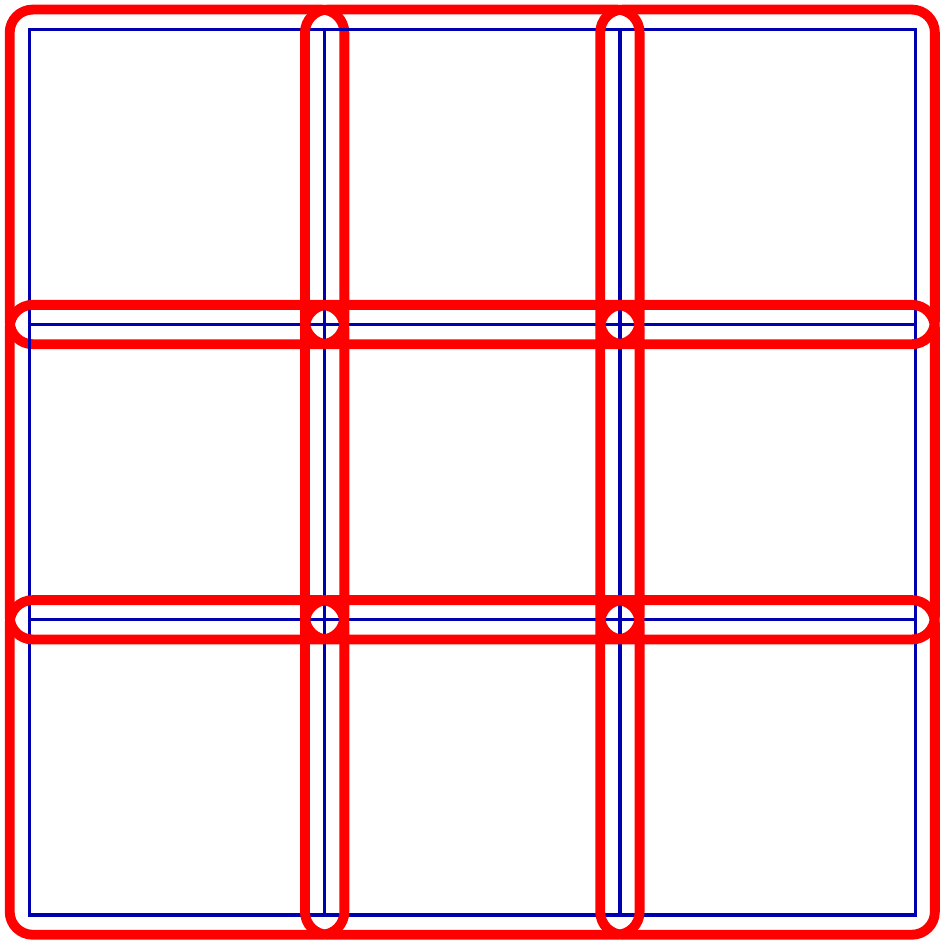}
\hspace*{1mm}\includegraphics[height=6cm, keepaspectratio]{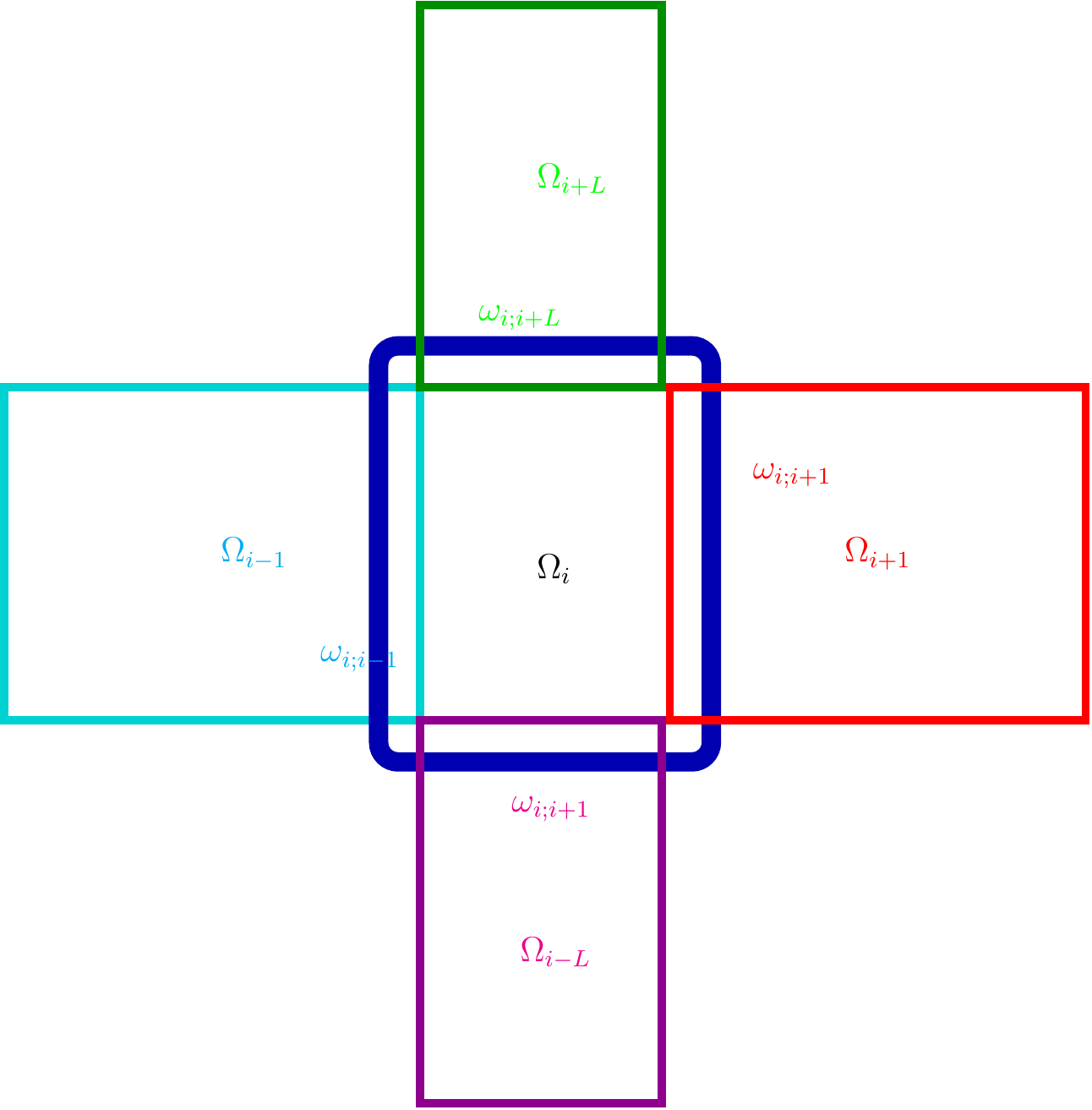}
\caption{(Left) Domain decomposition: Overlapping subdomains are represented in red. Blue subdomains do not overlap. (Right) Domain decomposition with overlapping region on $\Omega_i$ with $\Omega_{j,k,l,m}$ in $\R^2$}
\label{generalDDM}
\end{center}
\end{figure}

We first decompose $\R^{dN}\ni ({\bf x}_1,\cdots,{\bf x}_N)$, in $L^{dN}$ overlapping hypercubes $\Omega_i$ where $L$ is an integer parameter, $\cup_{i=1}^{L^{dN}}\Omega_i \subseteq \R^{dN}$ and apply a Schwarz waveform relaxation algorithm \cite{lorin-TBS,lorin-TBS2}. We present two different approaches. The first one leads to an a posteriori antisymmetric wavefunction, and second one ensures a priori Pauli's exclusion principle (see \ref{APXA}). In the following, we denote i) the artificial interfaces by $\Gamma_{i;j}=\partial \Omega_i \cap \Omega_j \subset \R^{dN-1}$, for any $i\neq j$, and ii) by $\omega_{i;j}$ the overlapping regions $\omega_{i;j}=\Omega_i\cap \Omega_j \subset \R^{dN}$, for any $(i,j) \in \{1,\cdots,L^{dN}\}^2$. For convenience, we also denote the Cartesian product $\R^{dN} = \R_1^d \times \R_{2}^d\cdots \times \R_N^d$, such that ${\bf x}_j \in \R^d_j$ for all $1\leq j\leq N$. We now denote by $\psi^{(k)}_i({\bf x}_1,\cdots,{\bf x}_N,t)$ the solution to the $N$-body TDSE in $\Omega_i$, at time $t$ and Schwarz iteration $k \geq 0$. For any $\Omega_i$, $1\leq i\leq L^{dN}$, we construct a basis of $K_i$ local basis functions (Gaussian functions or Slater's determinants) in $\Omega_i$, denoted by $\big\{v^i_j\big\}_{1\leq j\leq K_i}$ in order to compute $\psi_i^{(k)}$.   Basically, we will solve local time-dependent or time-independent local Schr\"odinger equations and reconstruct a global solution to the global Schr\"odinger equation. We then never compute the global solution from a global discrete Hamiltonian, but rather by computing $L^{dN}$ local wavefunctions (one per subdomain) using discrete local Hamiltonians, Figs. \ref{generalDDM}.  SWR algorithms are in particular, studied in \cite{halpern3,GanderHalpernNataf} and allow for a consistent decoupling on smaller subproblems of high dimensional (non-local) classical, quantum and relativistic wave equations.
\subsection{Schwarz Waveform Relaxation algorithm for the TDSE}\label{SWR1}
\begin{figure}[!ht]
\begin{center}
\hspace*{1mm}\includegraphics[height=6cm, keepaspectratio]{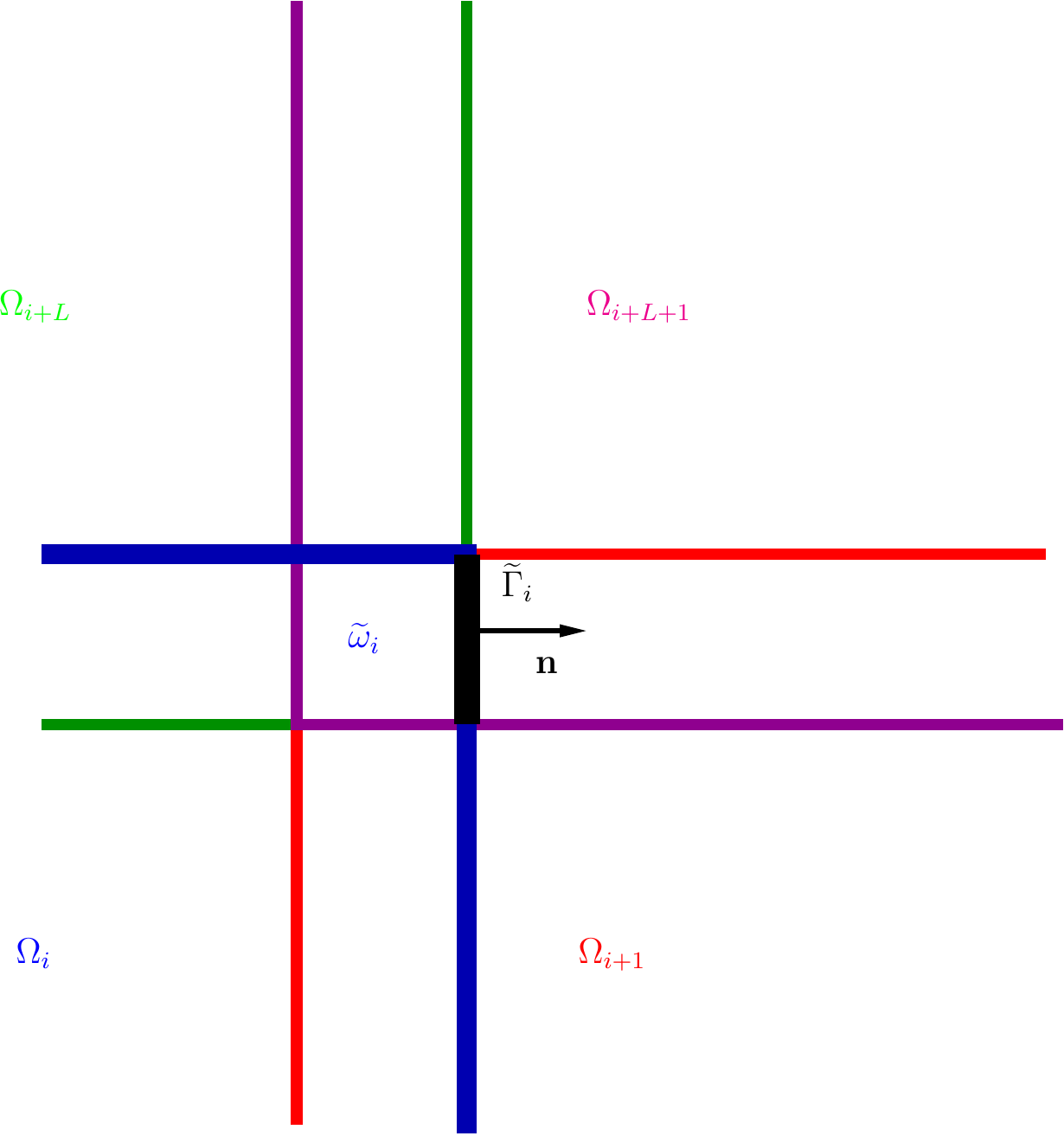}
\hspace*{1mm}\includegraphics[height=6cm, keepaspectratio]{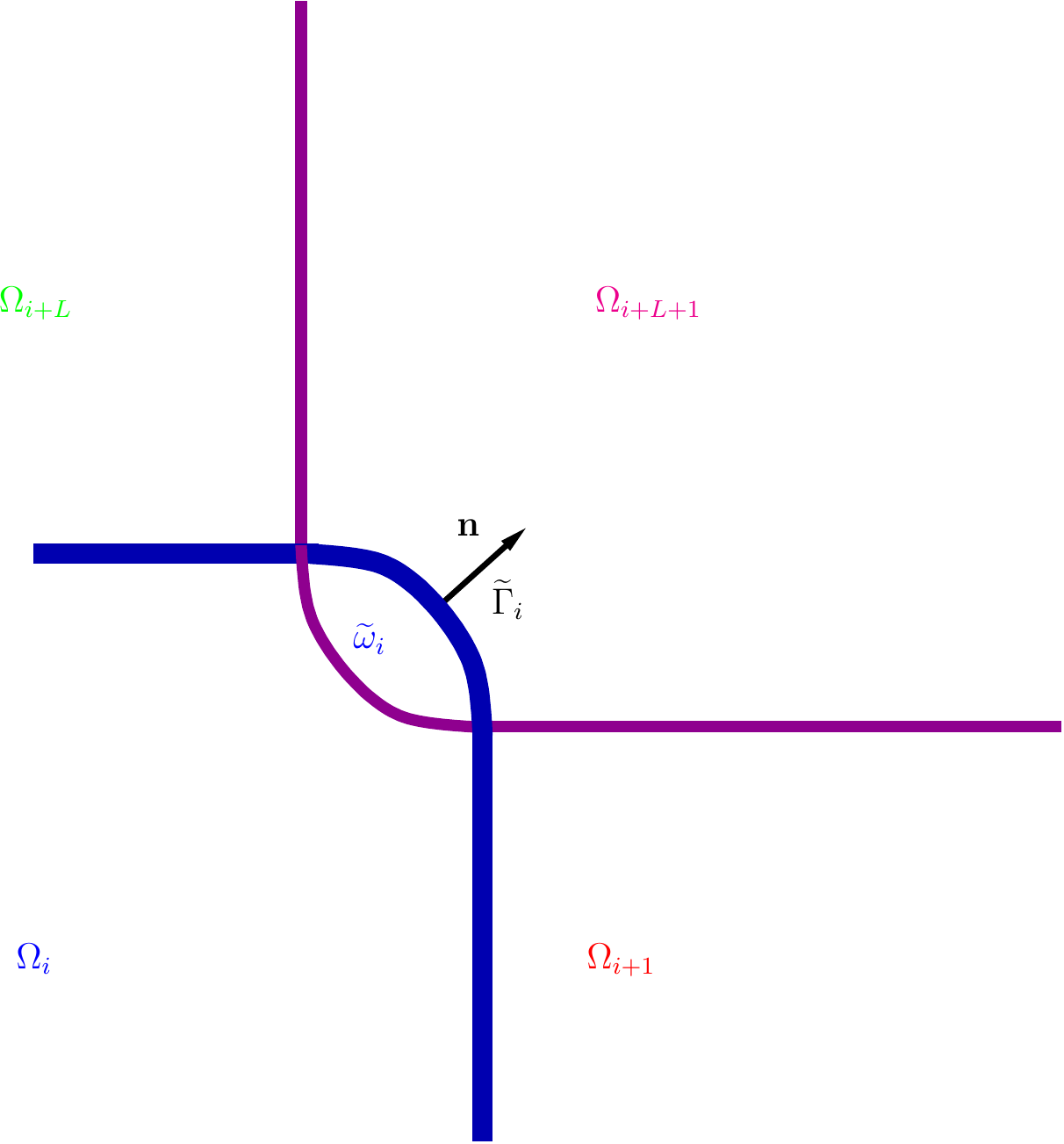}
\caption{(Left) Domain decomposition with $d=1$ and $N=2$: treatment of transmission conditions on interfaces belonging to more than $2$ subdomains. The black segment belonging to $\Omega_{i},\Omega_{i+1},\Omega_{i \pm L}$. (Right) Domain decomposition with $d=1$ and $N=2$: smooth subdomain boundary}.
\label{ddm3}
\end{center}
\end{figure}

We detail the DDM algorithm first for 2 subdomains $\Omega_{i}$, $\Omega_j$ with $i\neq j$, then for zones where more than $2$ subdomains overlap.  \\
\\
\noindent{\bf Two-subdomain overlapping zones.} Assume first that $\psi_i^{(k=0)}$ is a given function. The Schwarz Waveform Relaxation algorithm (SWR) $\mathcal{P}_i^{(k)}$ with $1\leq i\leq L^{dN}$ and $k \geq 1$, reads in LG and for $2$ subdomains, as 
\begin{eqnarray}\label{S1}
\hspace*{1cm}\mathcal{P}_i^{(k)} \,: \qquad \left\{
\begin{array}{lcll}
{\tt i}\partial_t \psi_i^{(k)} &  = & \Big(H_0 + \sum_{i=1}^N{\bf x}_i\cdot {\bf E}(t)\Big)\psi_i^{(k)} & \hbox{ on } \Omega_i \times (0,T),\\
\\
\psi_i^{(k)}(\cdot,0) &  = & \psi_{0|\Omega_i} & \hbox{ on } \Omega_i,\\
\\
\mathcal{B}_{i;j}\psi_i^{(k)} & = & \mathcal{B}_{i;j}\psi_j^{(k-1)} &  \hbox{ on } \Gamma_{i;j}\times (0,T)
\end{array}
\right.
\end{eqnarray}
where $\mathcal{B}_{i;j}$ is a transmission operator defined at $({\bf x}_1\cdots,{\bf x}_N) \in \Gamma_{i;j}=\partial \Omega_i \cap \Omega_j$.  \\
\\
\noindent{\bf Multi-subdomain overlapping zones.} The proposed decomposition requires a special treatment in zones, generically denoted $\widetilde{\omega}_i \subseteq \Omega_i$ see Fig. \ref{ddm3} (Left), where more than 1 subdomain overlap with $\Omega_i$. We denote by $\mathcal{O}(\widetilde{\omega}_i)$ the set of indices of the subdomains, distinct from $\Omega_i$, sharing the zone $\widetilde{\omega}_{i}$ with $\Omega_i$.  Notice that for interior subdomains (that is excluding the subdomains of the external layer) $\textrm{Card}\mathcal{O}(\widetilde{\omega}_i)=2^{dN}-1$.  The approach which is proposed is actually an averaging process. Let us generically denote by $\widetilde{\Gamma}_{i}$ the interface of $\widetilde{\omega}_i$ involved in the transmission conditions. The condition we impose at $\widetilde{\Gamma}_i$, thanks to the operator $\widetilde{\mathcal{B}}_i$, is defined by:
\begin{eqnarray*}
\widetilde{\mathcal{B}}_{i}\psi^{(k)}_i=\cfrac{1}{\textrm{Card}\mathcal{O}(\widetilde{\omega}_i)}\sum_{j \in \mathcal{O}(\widetilde{\omega}_i)}\mathcal{B}_{i;j}\psi^{(k-1)}_j.
\end{eqnarray*}
 In order to clarify the process, let us detail the case $d=1$ and $N=2$, with a total of $L^2$ subdomains. At a given interface $\widetilde{\Gamma}_i$ of $\widetilde{\omega}_i$, located at the right$/$top of a given subdomain $\Omega_i$ with $i\leq L(L-1)-1$, we assume that there are $2^{2\times 1}-1=3$ subdomains involved in the transmission condition, namely $\Omega_{i+1},\Omega_{i+L},\Omega_{i+L+1}$, see Fig. \ref{ddm3}. Then we impose at $\widetilde{\Gamma}_i$:
\begin{eqnarray*}
\widetilde{\mathcal{B}}_{i}\psi^{(k)}_i=\cfrac{1}{3}\big(\mathcal{B}_{i;i+1}\psi^{(k-1)}_{i+1}+\mathcal{B}_{i;i-L+1}\psi^{(k-1)}_{i-L+1}+\mathcal{B}_{i;i+L+1}\psi^{(k-1)}_{i+L+1}\big) \, .
\end{eqnarray*}
Notice also that a special treatment of the TC at the cross-points can improve the convergence deterioration at these locations. We do not explore this issue, but we refer to the related literature \cite{CP1,CP2}. A simple way to circumvent this difficulty consists in regularizing the corners of the subdomains, as shown in Fig. \ref{ddm3} (Right). In a Galerkin-method framework, managing such smooth regions is then straightforward. Another simple approach is presented in \cite{stcyr}, allows to avoid the discretization of the right-hand-side of the transmission conditions.
\\
\\
\noindent{\bf Selection of the transmission conditions.} From the convergence and computational complexity points of view, the selection of the transmission operator $\mathcal{B}_{i;j}$, is of crucial matter \cite{lorin-TBS}, \cite{halpern2}. The usual types of transmission conditions (TC) are now shortly recalled. The most simple approach is a Dirichlet-based TC's, where $\mathcal{B}_{i;j}$ is an identity operator. That is, we impose:
\begin{eqnarray*}
\psi^{(k)}_i = \psi^{(k-1)}_j \qquad \hbox{ on } \Gamma_{i;j}\times(0,T).
\end{eqnarray*}
In the literature, this is referred as the Classical Schwarz Waveform Relaxation (CSWR) algorithm. CSWR is in general convergent, and trivial to implement, but unfortunately exhibits usually very slow convergence rate and possibly stability issues at the discrete level \cite{halpern2,halpern3,jsc}. The CSWR method also necessitates an overlapping between the subdomains. As a consequence, more appropriate TC's should be derived, such as Robin-based TC's. Denoting ${\bf n}_{ij} \in \R^{dN}$, the outward normal vector to $\Gamma_{i;j}$, and for $\mu_{ij} \in \R_+^*$ in imaginary time (heat equation \cite{halpern3}), and $\mu_{ij} \in {\tt i}\R_{-}^*$ (real time, Schr\"odinger equation \cite{halpern2}), the Robin TC's read:
\begin{eqnarray}\label{RTC}
\big(\partial_{{\bf n}_{ij}}+\mu_{ij}\big)\psi^{(k)}_i = \big(\partial_{{\bf n}_{ij}}+\mu_{ij}\big)\psi^{(k-1)}_j \qquad \hbox{ on } \Gamma_{i;j}\times(0,T).
\end{eqnarray}
In the numerical simulations, $\mu_{ij}$ will be taken interface-independent, and will simply be denoted by $\mu$. This method will be referred as the Robin-SWR method. Along this paper, the Robin-SWR we will be used, as it is known for allowing for a good compromise between convergence rate and computational complexity, see \cite{halpern2}. Robin-SWR can be seen as an approximation of Optimal SWR (OSWR) which are based on transparent or high order absorbing TC's and reads, at $\Gamma_{i;j}\times(0,T)$
\begin{eqnarray*}
\big(\partial_{{\bf n}_{ij}}+{\tt i}\Lambda^{\pm}_{ij}\big)\psi^{(k)}_i = \big(\partial_{{\bf n}_{ij}}+{\tt i}\Lambda^{\pm}_{ij}\big)\psi^{(k-1)}_j \qquad \hbox{ on } \Gamma_{i;j}\times(0,T)
\end{eqnarray*}
where $\Lambda^{\pm}_{ij}$ is a pseudodifferential Dirichlet-to-Neumann (DtN) operator, \cite{nir,hor1,hor2} and derived from the Nirenberg factorization: 
\begin{eqnarray*}
P = \big(\partial_{{\bf n}_{ij}}+{\tt i}\Lambda^{+}_{ij}\big)\big(\partial_{{\bf n}_{ij}}+{\tt i}\Lambda^{-}_{ij}\big) + R
\end{eqnarray*}
where $P$ is the time-dependent $N$-particle Schr\"odinger operator, $R\in$ OPS$^{-\infty}$, and $\Lambda^{\pm}_{ij}$ are operators associated to outgoing$/$incoming waves.  Robin-SWR then consists in approximating the pseudodifferential operators $\Omega_{ij}$ by a an algebraic operator $\mu_{ij}$. OSWR and quasi-OSWR are applied and analyzed to linear Schr\"odinger equations in \cite{halpern2,jsc}, and show much faster convergence than CSWR, but are also much more computationally complex to implement.
\\
Finally, a commun convergence criterion for the Schwarz DDM is set for all $i\neq j$ in $\{1,\cdots,L^{dN}\}$ by
\begin{eqnarray}\label{CVTOT}
\big\| \hspace{0.2cm} \sum_{i=1}^{L^{dN}}\|\psi^{(k)}_{i|\Gamma_{i;j}}-\psi^{(k)}_{j|\Gamma_{i;j}}\|_{L^2(\Gamma_{i;j})}\big\|_{L^{2}(0,T)} \leq  \delta^{\textrm{Sc}}.
\end{eqnarray}
 When the convergence of the whole iterative SWR algorithm is obtained at Schwarz iteration $k^{(\textrm{cvg})}$, then one gets the converged global solution 
$\psi^{\textrm{cvg}}:=\psi^{(k^{\textrm{cvg}})}$, typically with $\delta^{\textrm{Sc}}=10^{-14}$ (''Sc'' for Schwarz). 
\\
\\
\noindent{\bf Construction of the local approximate solutions.} The construction of approximate solutions to \eqref{S1} is now performed. The local wavefunction $\psi_i^{(k)}$ is expanded as follows:
\begin{eqnarray}\label{exp2}
\psi^{(k)}_i(\cdot,t) = \sum_{j=1}^{K_i}c_j^{i,(k)}(t)v_j^{i}
\end{eqnarray}
where $\big\{v_j^{i}\big\}_{1\leq j \leq K_i}$ are the local basis functions\footnote{constructed from $2M_i$ orbitals $\big\{\phi^i_j({\bf x})\big\}_{1\leq j\leq 2M_i}$, or as Gaussian functions, following the same strategy as presented in Subsections \ref{subsec:SLO}. When the local basis functions are chosen as local Slater's determinants, we should have $M_i \ll M$, where $M$ is the number of basis functions, for the 1-domain FCI.} associated to $\Omega_i$. In $c_j^{i,(k)}$ and $\psi^{(k)}_i$ the index $k$ refers to the Schwarz iteration, $i$ refers to the subdomain index, and $j$ to the basis function index. The expansions \eqref{exp2} lead to $L^{dN}$-independent linear systems of ODE's, $1\leq i\leq L^{dN}$
\begin{eqnarray*}
{\tt i}{\bf A}_i\dot{\bf c}^{i,(k)}(t) = \big(\widetilde{{\bf H}}_{i} + {\bf T}_i(t)\big){\bf c}^{i,(k)}(t)
\end{eqnarray*}
where ${\bf A}_{i} \in M_{K_i}(\R)$ stands for $A_{i;(j,l)} = \langle v^j_i , v^l_i\rangle$ and $\widetilde{{\bf H}}_{i} \in M_{K_i}(\R)$ stands for the discrete local Hamiltonian $\widetilde{H}_{i;(j,l)} = \langle v^j_i , (H_{0}+R_i) v^l_i\rangle$ including the contribution to the Robin transmission operator $R_i$, where $1\leq j,l \leq K_i$. For $d=3$, with ${\bf E}(t)=\big(E_x(t),E_y(t),E_z(t)\big)^T$ with $t\in (0,T)$, we have:
\begin{eqnarray*}
{\bf T}_i(t) = \widetilde{{\bf H}}_{i} + E_x(t){\bf Q}^{x}_i + E_y(t){\bf Q}^{y}_i + E_z(t){\bf Q}^{z}_i
\end{eqnarray*}
where 
\begin{eqnarray}\label{Q}
\left.
\begin{array}{lcl}
\big\{{\bf Q}^{x}_{i;(j,l)}\big\}_{1\leq j,l\leq K_i} & = & \Big\{\sum_{m=1}^N\langle x_m v_j^{i}, v_l^i\rangle\Big\}_{1\leq j,l\leq K_i}, \\
 \big\{{\bf Q}^{y}_{i;(j,l)}\big\}_{1\leq j,l\leq K_i} & = & \Big\{\langle y_mv_j^{i}, v_l^i\rangle\Big\}_{1\leq j,l\leq K_i}, \\ 
\big\{{\bf Q}^{z}_{i;(j,l)}\big\}_{1\leq j,l\leq K_i} & = & \Big\{\langle z_mv^i_{j}, v^i_l\rangle \Big\}_{1\leq j,l\leq K_i}.
\end{array}
\right.
\end{eqnarray}
Matrices ${\bf Q}_i^{x,y,z}$, $\widetilde{{\bf H}}_{i}$ are computed once for all, in each subdomain. Efficient computation of $\widetilde{{\bf H}}_{i}$ is presented in \cite{CAM15-09}, \cite{CAM15-10}.
\subsection{Schwarz Waveform Relaxation algorithm for the TISE}\label{SWR2}
The DDM which is proposed above for the time-dependent Schr\"odinger equation can be implemented for computing Schr\"odinger Hamiltonian's bound states, using the {\it imaginary time} method also referred in the literature as the {\it Normalized Gradient Flow} (NGF) method \cite{bao,lorin-TBS,lorin-TBS2}.  It basically consists of replacing {\it $t$ by ${\tt i}t$} in TDSE, {\it normalizing the solution} at each time iteration, which is finally convergent to an eigenfunction to $H_0$, by default the ground state. In the imaginary time framework, the SWR domain decomposition is similar as above, and the same notation as Sections \ref{SWR}, \ref{NAWF} and \ref{APXA}. In any subdomain $\Omega_i$, we define for $K_i \in \N^*$ 
\begin{eqnarray}\label{exp3}
\phi^{(k)}_i(\cdot,t) = \sum_{j=1}^{K_i}c_j^{i,(k)}(t)v_j^{i}
\end{eqnarray}
where $\big\{v^i_{j}\big\}_{1\leq j \leq K_i}$ are the local basis functions. Denoting the discrete times $t_{0}:=0<t_{1}<...<t_{n+1}<...$ with $t_{n+1}=t_n+\Delta t$ for some $\Delta t>0$, and the initial guess $\phi_{0}$, the SWR method $\mathcal{I}_i^{(k)}$ for $1\leq i\leq L^{dN}$, $1\leq j \leq L^{dN}$ with nonempty $\omega_{i;j}:=\Omega_i\cap \Omega_j$, and $k \geq 1$,  when only 2 subdomains overlap, then reads:
\begin{eqnarray}\label{SIT}
\hspace*{1cm}\mathcal{I}_i^{(k)} \,: \qquad \left\{
\begin{array}{lcl}
\partial_t \phi_i^{(k)} & = & -H_0 \phi_i^{(k)}, \, \hbox{ on }\Omega_i\times (t_{n},t_{n+1}), \\
\\
\mathcal{B}_{i;j}\phi_i^{(k)} & =& \mathcal{B}_{i;j}\phi_j^{(k-1)}, \, \hbox{ on } \Gamma_{i;j}\times (t_{n},t_{n+1}),\\
\\
\phi_i^{(k)}(\cdot,0) & = & \phi_0, \hbox{ on }\Omega_{i},\\
\\
\displaystyle \phi_i^{(k)}(\cdot,t_{n+1}) & =& \phi_i^{(k)}(\cdot,t^{+}_{n+1})=\frac{\phi_i^{(k)}(\cdot,t^{-}_{n+1})}{|| \sum_{j=1}^{L^{dN}}\tilde{\phi}_j^{(k)}(\cdot ,t^{-}_{n+1})||_{2}}, \, \hbox{ on }\Omega_i
\end{array}
\right.
\end{eqnarray}
where $\mathcal{B}_{i;j}$ is a transmission operator written in imaginary time, and defined at $({\bf x}_1\cdots,{\bf x}_N) \in \Gamma_{i;j}$, and  $\tilde{\phi}^{(k)}_i$ denotes the extension by $0$ to $\R^{dN}$ of $\phi^{(k)}_i$. As discussed in Subsection \ref{SWR1}, in the regions covered by more than $2$ subdomains, a special treatment of the transmission condition is necessary, but is strictly identical to the time-dependent case, see Section \ref{SWR1}.\\
For a given Schwarz iteration $k$, we stop the NGF computations when the {\it reconstructed} approximate solution $\phi^{n+1,(k)}$ satisfies at time $t_{n+1}$
\[
|| \phi^{n+1,(k)} - \phi^{n,(k)}||\leq \delta \, .
\]
with $\delta$ a small parameter.  When the convergence is reached, then the stopping time is such that:
$T^{(k)}:=T^{\textrm{cvg},(k)}=n^{\textrm{cvg},(k)}\Delta t$  for a converged solution $\phi^{\textrm{cvg},(k)}$ reconstructed from the $L^{dN}$ subdomain solutions $\phi_i^{\textrm{cvg},(k)}$. A convergence criterion for the Schwarz DDM is set, for all $i\neq j$ in $\{1,\cdots,L^{dN}\}$ by
\begin{eqnarray}\label{CVTOT2}
\big\| \hspace{0.2cm} \sum_{i=1}^{L^{dN}} \|\phi^{\textrm{cvg},(k)}_{i|\Gamma_{i;j}}-\phi^{\textrm{cvg},(k)}_{j|\Gamma_{i;j}}\|_{L^2(\Gamma_{i;j})}\big\|_{L^{2}(0,T^{(k^{\textrm{cvg}})})} \leq  \delta^{\textrm{Sc}} \,.
\end{eqnarray}
In the numerical experiments, we will use a bit different criterion. When the convergence of the whole iterative SWR$/$NGF algorithm is obtained at Schwarz iteration $k^{\textrm{cvg}}$ one then gets the converged global eigenstate $\phi^{\textrm{cvg}}:=\phi^{\textrm{cvg},(k^{\textrm{cvg}})}$ typically with $\delta^{\textrm{Sc}}=10^{-14}$.
\subsection{Wavefunction reconstruction}\label{NAWF}
The main weakness of the basic decomposition above is that by default, it does not ensure the antisymmetry of the overall wavefunction in $\R^{dN}$. Indeed, in each subdomain $\Omega_i$ and any time $t\in (0,T$)
\begin{eqnarray*}
\phi_i^{(k)}({\bf x}_1,\cdots,{\bf x}_N,t) = \sum_{j=1}^{K_i}v_j^i({\bf x}_1,\cdots,{\bf x}_N)c_j(t)
\end{eqnarray*} 
is a priori not antisymmetric, and a fortiori the reconstructed solution $\phi^{(k)}$ (in $\R^{dN}$). Indeed
\begin{enumerate} 
\item if $({\bf x}_1,\cdots,{\bf x}_N) \in \omega_{i;j}$ where $j$ is unique (that is zones where only two subdomains overlap), then 
\begin{eqnarray*}
\phi^{(k)}({\bf x}_1,\cdots,{\bf x}_N,t) = 
\left\{
\begin{array}{l}
\phi_i^{(k)}({\bf x}_1,\cdots,{\bf x}_N,t), \, ({\bf x}_1,\cdots,{\bf x}_N) \in \Omega_i-\omega_{i;j}, \, \forall (i,j) \in \{1,\cdots,L^{dN}\}^2\\
\\
\cfrac{\phi_i^{(k)}+\phi_j^{(k)}}{2}({\bf x}_1,\cdots,{\bf x}_N,t), \, ({\bf x}_1,\cdots,{\bf x}_N) \in \omega_{i;j}, \, \forall (i,j) \in \{1,\cdots,L^{dN}\}^2.
\end{array}
\right.
\end{eqnarray*}
\item if $({\bf x}_1,\cdots,{\bf x}_N) \in \omega_{i;j}$ where $j$ is not unique, that is there exists a zone denoted by $\widetilde{\omega}_i$, where $\textrm{Card}\mathcal{O}(\widetilde{\omega}_i) \leq 2^{dN}-1$ subdomains, $\{\Omega_{i_j}\}_{j \in\mathcal{O}(\widetilde{\omega}_i)}$, have a common overlap with $\Omega_i$. Then for $({\bf x}_1,\cdots,{\bf x}_N) \in \widetilde{\omega}_{i}$
\begin{eqnarray*}
\phi^{(k)}({\bf x}_1,\cdots,{\bf x}_N,t) = \cfrac{1}{\textrm{Card}\mathcal{O}(\widetilde{\omega}_i)+1}\big(\phi_i^{(k)}+\sum_{j=1}^{\mathcal{O}(\widetilde{\omega}_i)}\phi_{i_j}^{(k)}\big)({\bf x}_1,\cdots,{\bf x}_N,t).
\end{eqnarray*}
In fact, if the local basis functions are Slater's determinants basis functions, then:
\begin{eqnarray*}
\phi^{(k)}({\bf x}_1,\cdots,{\bf x}_p,\cdots,{\bf x}_q,\cdots,{\bf x}_N,t) = -\phi^{(k)}({\bf x}_1,\cdots,{\bf x}_q,\cdots,{\bf x}_p,\cdots,{\bf x}_N,t)
\end{eqnarray*}
occurs only if $({\bf x}_1,\cdots,{\bf x}_p,\cdots,{\bf x}_q,\cdots,{\bf x}_N)$ {\it and} $({\bf x}_1,\cdots,{\bf x}_q,\cdots,{\bf x}_p,\cdots,{\bf x}_N)$ belong to $\Omega_i$. Then, even when local Slater's determinants are constructed, a more careful decomposition is then necessary to ensure a global antisymmetry.
\end{enumerate}
One possible approach is to antisymmetrize at each time step the reconstructed wavefunction, thanks to the operator $\mathcal{A}$ defined in Subsection \ref{ITM}. We also propose in \ref{APXA}, an algorithm to ensure the antisymmetry of the reconstructed wavefunction within the SWR formalism.
\subsection{Numerical algorithm}\label{expcons}
We give details about the explicit construction of the numerical solver. Let us denote by ${\boldsymbol \psi}_i^{n,(k)}(x_1,x_2) = \sum_{j=1}^{K_i}c_{j}^{i,n,(k)}v_j^i(x_1,x_2)$ the approximate solution in $\Omega_i$, at Schwarz iteration $k$, and time $t_n$. We denote by $\widetilde {\bf H}_i$ the discrete Hamiltonian in $\Omega_i$ including the transmission condition contribution. The discrete parallel algorithm in real (resp. imaginary) time reads as follows.\\
\\
\noindent Schwarz iteration loop, from $k=0$ to convergence $k=k^{\textrm{(cvg)}}$:
\begin{enumerate}
\item At initial real (resp. imaginary) time $t=0$, we restrict $\phi_0$ to $\Omega_i$, then project $\phi_{0|\Omega_i}$ onto the local basis functions $\big\{v^i_j\}$, where $i \in \{1,\cdots,L^{dN}\}$ is the subdomain index, and where $j \in \{1,\cdots,K_i\}$ is the local basis function index, in order to construct the local Cauchy data $\phi_{0;i}^{(k)} = \sum_{j=1}^{K_i}c_{j}^{i,0,(k)}v_j^i$. Additional details will also be presented in Subsection \ref{testA}.
\item Real (resp. imaginary) time iterations, from $n=0$ to $n=n_T$ (resp. $n=0$ to $n=n^{\textrm{cvg},(k)}$), that is from time $t_0=0$ to time $t_{n_{T}}=T$ (resp. $t_{n^{\textrm{cvg},(k)}}=T^{\textrm{cvg},(k)}$) to update the basis coefficients ${\bf c}^{i,n+1,(k)} = \{c_{j}^{i,n+1,(k)}\}_{1 \leq j \leq K_i}$ (resp. $\widetilde{{\bf c}}^{i,n+1,(k)}$) from ${\bf c}^{i,n,(k)} = \{c_{j}^{i,n,(k)}\}_{1 \leq j \leq K_i}$ (resp. ${\bf c}^{i,n,(k)}$), by solving, in real time:
\begin{eqnarray*}
\Big({\bf A}_i+{\tt i}\cfrac{\Delta t}{2}\widetilde{{\bf H}}_i+{\tt i}\cfrac{\Delta t}{2}{\bf T}_i^{n+1}\Big){\bf c}^{i,n+1,(k)} = \Big({\bf A}_i-{\tt i}\cfrac{\Delta t}{2}\widetilde{{\bf H}}_i-{\tt i}\cfrac{\Delta t}{2}{\bf T}_i^{n}\Big){\bf c}^{i,n,(k)}
\end{eqnarray*}
(resp. imaginary time: $\Big({\bf A}_i+\Delta t\widetilde{{\bf H}}_i\Big)\widetilde{{\bf c}}^{i,n+1,(k)} = {\bf A}_i{\bf c}^{i,n,(k)}$)
\item Reconstruction of the global TDSE wavefunction ${\boldsymbol \psi}^{i,n+1,(k)}$ (resp. $\widetilde{{\boldsymbol \phi}}^{n+1,(k)}$).
\item In imaginary time, {\it only}:  $L^2$-normalization of the local wavefunctions in imaginary time, that is
\begin{eqnarray*}
{\bf c}^{i,n+1,(k)} = \cfrac{\widetilde{{\bf c}}^{i,n+1,(k)}}{\|\sum_{j=1}^{L^{dN}}\widetilde{{\boldsymbol \phi}}_j^{i,n+1,(k)}\|_{2}}
\end{eqnarray*}
and antisymmetrization, thanks to the operator $\mathcal{A}$, see Subsection \ref{ITM}.
\item At final real (resp. imaginary) time $T$ (resp. $T^{\textrm{(cvg)},(k)}$) and Schwarz iteration $k$, we have determined $\psi^{n_T,(k)}$ (resp. $\phi^{n^{\textrm{(cvg)}},(k)}$).
\end{enumerate}
\noindent At convergence of the Schwarz algorithm, we get $\psi^{n_T,(k^{\textrm{(cvg)}})}$ (resp. $\phi^{n^{\textrm(cvg)},(k^{\textrm{(cvg)}})}$) which is then an approximation of $\psi(\cdot,T)$ (resp. ground state of $H_0$).\\
\\
Notice that i) the implicit Euler scheme guarantees the local unconditional stability of the imaginary time solver \cite{bao}, and ii) a Crank-Nicolson scheme guarantees also the unconditional stability, as well as the convergence at order $2$ in space and time, for the real time solver, see Step 2.
\subsection{Convergence of the SWR algorithms}
Although a rigorous analysis of the presented SWR-DDM solver applied to the $N$-body Schr\"odinger equation is out of reach, we can provide some useful references and some mathematical properties of the presented algorithms. We here summarize some of the known results about the convergence of SWR algorithm, as well as their rate of convergence for the Schr\"odinger equation in real and imaginary time. Notice that these results are usually established for two subdomains in 1-d. Basically, it consists of i) reformulating the SWR method as a fixed point algorithm, and when necessary ii)  using pseudodifferential calculus in order to derive the corresponding contraction factor, as a function of the overlap size and of the frequencies of the wavefunction. The principle of proof is similar in real and imaginary time, although a finer analysis is necessary in real time, as it requires a closer study of the contraction factor in three different zones (hyperbolic, elliptic glancing zones), see \cite{nir}. 
In 1-d, it was proven in \cite{halpern2}, that the convergence of the Classical, Robin and q-Optimal Schwarz (CSWR, Robin-SWR, q-OSWR) methods for the real time Schr\"odinger equation, with differentiable and bounded potential with a bounded derivative. Notice, that unlike the CSWR method, the Robin and q-Optimal SWR do not require an overlap between the subdomains. These methods were used in \cite{jsc} in a laser-particle setting involving in particular recollision and ionization, where transmission conditions were derived from high order absorbing boundary conditions. We also refer to \cite{BesseXing} for some numerical implementation and performance of CSWR and Robin-SWR methods for the time-dependent Schr\"odinger equation in 2-d. In imaginary time with space-dependent potential, the convergence of the CSWR and q-OSWR methods have been established in 1-d for high frequency problems, as well as their rate of convergence as a function of the size  of the overlapping zone $\epsilon$, see \cite{lorin-TBS}. It was proven that in the case of the CSWR algorithm, the rate of convergence is (at first order) exponential in $-\epsilon \sqrt{|\tau|}$, where $\tau$ denotes the co-variable (frequency) associated to $t$ (time), and that for positive potentials accelerate the convergence of the algorithm. Quasi-OSWR methods are shown to accelerate the CSWR by a factor $|\tau|^{-p}$, for some $p \in \N^*$, dependent on the order of approximation $p$, of the q-OSWR method. In fine, we get
\begin{eqnarray*}
\hspace*{1cm}\lim_{k \rightarrow +\infty}\|\psi_{|\Omega_i} - \psi_i^{(k)}\|_{L^2(\Omega_i\times(0,T))} = 0 \, .
\end{eqnarray*}
The convergence of the CSWR in 2-d for two subdomains with smooth convex$/$concave boundary is established in imaginary time in \cite{lorin-TBS2}. As in the one-dimensional setting, the rate of convergence is exponential in $-\epsilon \sqrt{|\tau|}$, where $\epsilon$ characterizes the thickness of the overlapping region, but it is also established a deceleration effect of the interface curvature, suggesting that flat interfaces are preferable than curved ones. Similarly the rate of convergence for the CSWR, and q-OSWR methods can be established for the time dependent Schr\"odinger with space-dependent potentials \cite{lorin-TBS3}. In \ref{APXC}, we analyze the computational complexity of the SWR method applied to the TISE and TDSE. \\
Notice that in order to accelerate the rate of convergence of the domain decomposition method, a multilevel approach should be coupled to the proposed method. Indeed, although to our knowlegde there is no rigorous proof, we expect that the larger the number of subdomains, the larger the number Schwarz iterations to converge, in particular in the time-dependent case. As it is was proposed in \cite{multilevel}, a multilevel strategy helps to accelerate the convergence of the SWR method for the Schr\"odinger equation. It was shown that in the case of the NGF method, a multilevel approach allows for an acceleration of the convergence of the NGF algorithm at each Schwarz iteration. Regarding the time-dependent case, a substantial acceleration of the Schwarz algorithm is observed.

\section{Numerical experiments: Gaussian local basis functions}\label{NumGauss}
This section is devoted to some numerical experiments in imaginary and real time using local Gaussian basis functions in order to validate the methodology developed in this paper for a one-dimensional 2-electron problem, that is with $d=1$ and $N=2$. Numerically this corresponds to two-dimensional time-dependent problems. Naturally, DDM is not necessary for this low-dimensional problem, but we intend here to show that the presented method is indeed convergent, and is a good candidat in higher dimension where DDM becomes relevant. Realistic simulations in higher dimensions will be presented in a forthcoming paper.\\
We assume that the global domain $\Omega=[a-\epsilon^{(x_1)}/2,b+\epsilon^{(x_1)}/2]\times[c-\epsilon^{(x_2)}/2,d+\epsilon^{(x_2)}/2]$ is uniformly decomposed in $L^2$ subdomains, $\Omega_{i+Lj}=[a_i-\epsilon^{(x_1)}/2,b_i+\epsilon^{(x_1)}/2]\times[c_j-\epsilon^{(x_2)}/2,d_j+\epsilon^{(x_2)}/2]$, for all $i,j=1,\cdots,L$ and $\epsilon^{(x_{1,2})}>0$. That is $a_i=a+(i-1)(b-a)/L$ (resp. $c_i=c+i(d-c)/L$) and $b_i=a+i(b-a)/L$ (resp. $d_i=c+i(d-c)/L$), for $i=1,\cdots,L$. The overlapping zone in each direction (North, West, South, East) is a band of size $\epsilon^{(x_2)}\times (b-a)/L$ (South, North) or $\epsilon^{(x_1)}\times (d-c)/L$ (East, West).  We denote by $N^{(x_1)}$ and $N^{(x_2)}$ the total number of grid points, in each coordinate. The total number of local basis functions will be assumed, for the sake of simplicity, to be subdomain-independent, and is denoted by $K:=K_i$, for $i=1,\cdots,L^2$.
\subsection{Test 1.a: Gaussian local basis function construction}\label{testA}
In this first test, we simply represent the local basis functions, and the reconstructed Gaussian Cauchy data. More specifically, we assume that $L^2=25$ subdomains cover a two-dimensional global domain $(-15,15)$, that is $a=-15$ and $b=15$. We construct  $K=N_{\phi}^2=6^2$ Gaussian local basis functions per subdomain. Say for a subdomain $\Omega_{i+jL}=[a_i-\epsilon^{(x_1)}/2,b_i+\epsilon^{(x_1)}/2]\times[c_j-\epsilon^{(x_2)}/2,d_j+\epsilon^{(x_2)}/2]$, the basis functions are constructed as:
\begin{eqnarray*}
v_{l,p}^{i,j}(x_1,x_2) = \exp\big(-0.4(x_1-\alpha^{(l)}_i)^2-0.4(x_2-\beta^{(p)}_i)^2\big)
\end{eqnarray*}
where $\alpha^{(l)}_i$ (resp. $\beta^{(p)}_j$) for $l=1,\cdots,N_{\phi}$ (resp. $p=1,\cdots,N_{\phi}$), are uniformly distributed numbers in $[a_i-\epsilon^{(x_1)}/2,b_i+\epsilon^{(x_1)}/2]$ (resp. $[c_j-\epsilon^{(x_2)}/2,d_j+\epsilon^{(x_2)}/2]$).  From now on, in order to simplify the notations, we will denote the basis functions $v^{i}_l$ (see Subsection \ref{notations}), for $l=1,\cdots,N_{\phi}^2$ and for $i=1,\cdots,L^2$.  The overlapping zone on each subdomain, represents $\approx 10\%$ of the overall subdomain. The total number of grid points is fixed at $N^{(x_1)}\times N^{(x_2)}=201^2$. We represent in Fig. \ref{CauchyRec} (Left), the coverage of a given subdomain by the $N_{\phi}^2$ basis functions, from above.  The reconstructed function $\phi^{(0)}$, defined on $\Omega$ from $\phi_0$
\begin{eqnarray*}
\phi_0(x_1,x_2) = \exp\big(-0.2(x_1^2+x_2^2)\big)
\end{eqnarray*}
is computed as follows:
\begin{itemize}
\item $i \in \{1,\cdots,L^2\}$, we construct the sparse matrices ${\bf A}_{i}=\{A_{i;(l,p)}\}_{1\leq l \leq N_{\phi}^2,1 \leq p\leq N_{\phi}^2}$,
\begin{eqnarray*}
A_{i;(l,p)} = \langle v_{l}^{i},v_p^{i} \rangle_{L^2(\Omega_{i})}, \qquad \forall (l,p)=\{1,\cdots,N^2_{\phi}\}^2.
\end{eqnarray*}
\item For each  $i \in \{1,\cdots,L^2\}$, we restrict $\phi_0$ to $\Omega_{i}$: $\phi_{0|\Omega_{i}}$.
\item For any $i \in \{1,\cdots,L^2\}$, we project $\phi_{0|\Omega_{i}}$ on each local basis functions, and construct the coefficients $\widetilde{{\bf c}}^{i} = \{\widetilde{c}^i_l\}_l$ defined by
\begin{eqnarray*}
\widetilde{c}_l^{i} = \langle \phi_{0|\Omega_{i}},v_{l}^{i} \rangle_{L^2(\Omega_{i})}, \qquad \forall l=\{1,\cdots,N^2_{\phi}\}.
\end{eqnarray*}
\item Then ${\bf c}^{i}$ is solution to ${\bf A}_i{\bf c}^{i} = \widetilde{\bf c}^{i}$, using GMRES \cite{saad}.
\item We can then reconstruct the local solution as follows: $\phi_{0;i}(x_1,x_2)=\sum_{l=1}^{N^2_{\phi}}c^{i}_{l}v_{l}^{i}(x_1,x_2)$.
\item  We finally denote by $\phi^{(0)}_0(x_1,x_2)$ the reconstructed initial data, which is then given, for $i \in \{1,\cdots,L^2\}$,  by
\begin{eqnarray*}
\phi^{(0)}(x_1,x_2) =
\left\{
\begin{array}{ll}
\phi^{(0)}_{i}(x_1,x_2), & \hbox { if } (x_1,x_2)\in \Omega_{i} \, \hbox{ only} \\
\cfrac{1}{k}\sum_{l=1}^{k}\phi^{(0)}_{k}(x_1,x_2)&, \hbox { if } (x_1,x_2) \in \cap_{l=1}^k\Omega_{i_k} \, \hbox{ with } i_1,\cdots,i_k \in \{1,\cdots,L^2\}  
\end{array}
\right.
\end{eqnarray*}
\begin{figure}[!ht]
\begin{center}
\hspace*{1mm}\includegraphics[height=6cm, keepaspectratio]{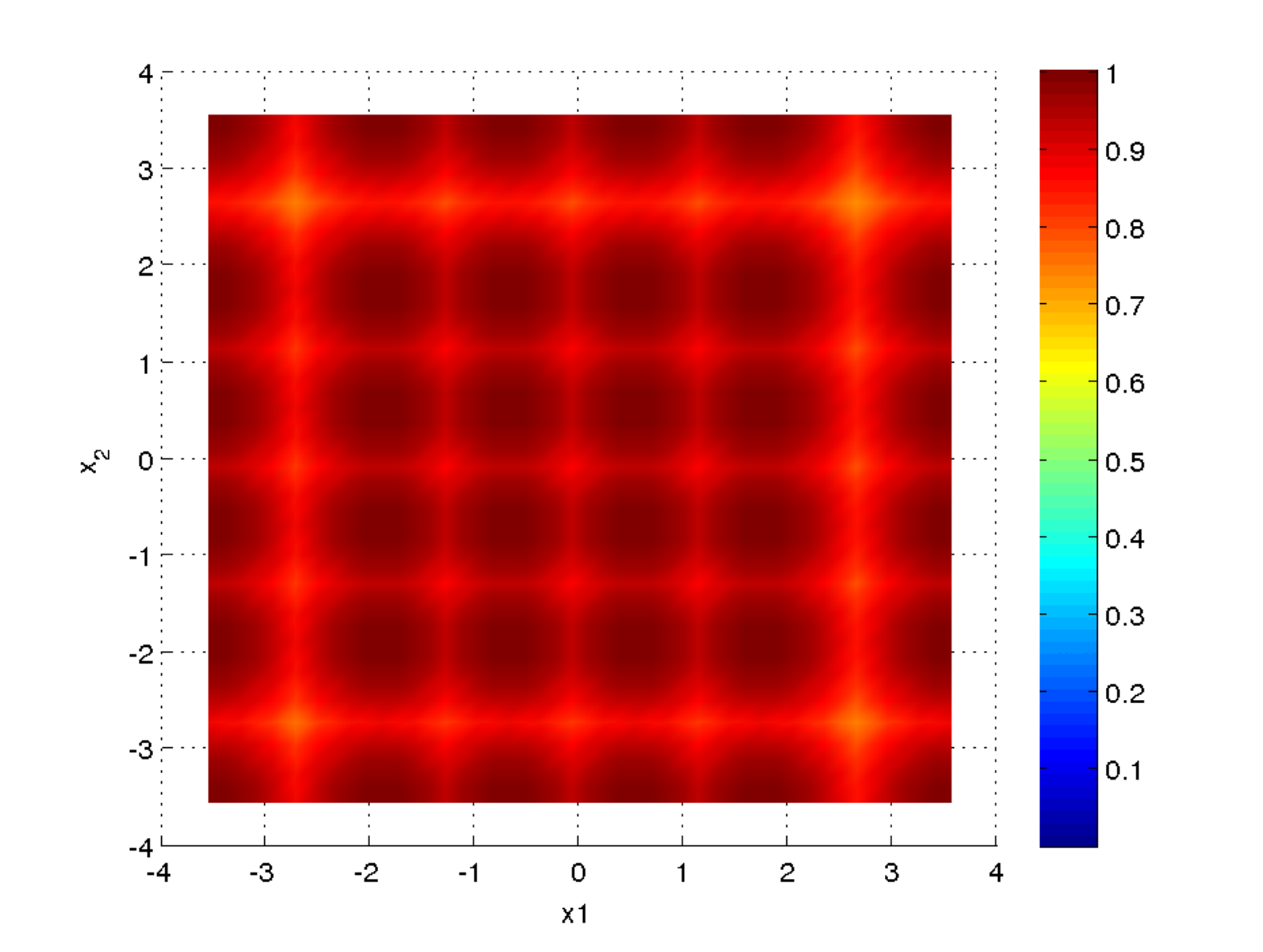}
\hspace*{1mm}\includegraphics[height=6cm, keepaspectratio]{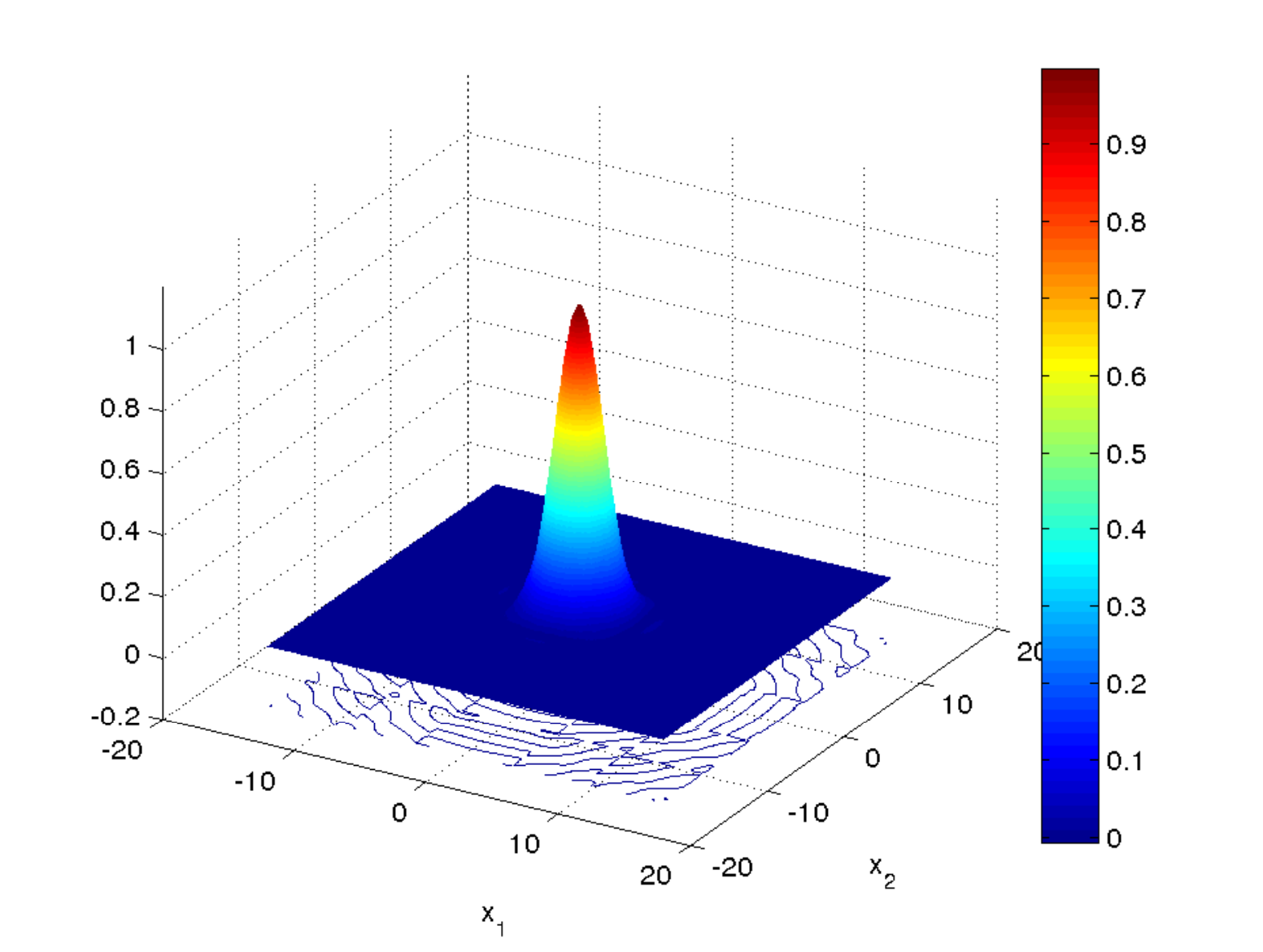}
\caption{(Left) $36$ Gaussian basis functions in a subdomain. (Right) Reconstruction of a given function from local basis functions.}
\label{CauchyRec}
\end{center}
\end{figure}
\end{itemize}
We show in Fig. \ref{CauchyRec} (Right) the reconstructed wavefunction $\phi^{(0)}$.
\subsection{Test 1.b: Heat equation}
\begin{figure}[!ht]
\begin{center}
\hspace*{1mm}\includegraphics[height=4cm, keepaspectratio]{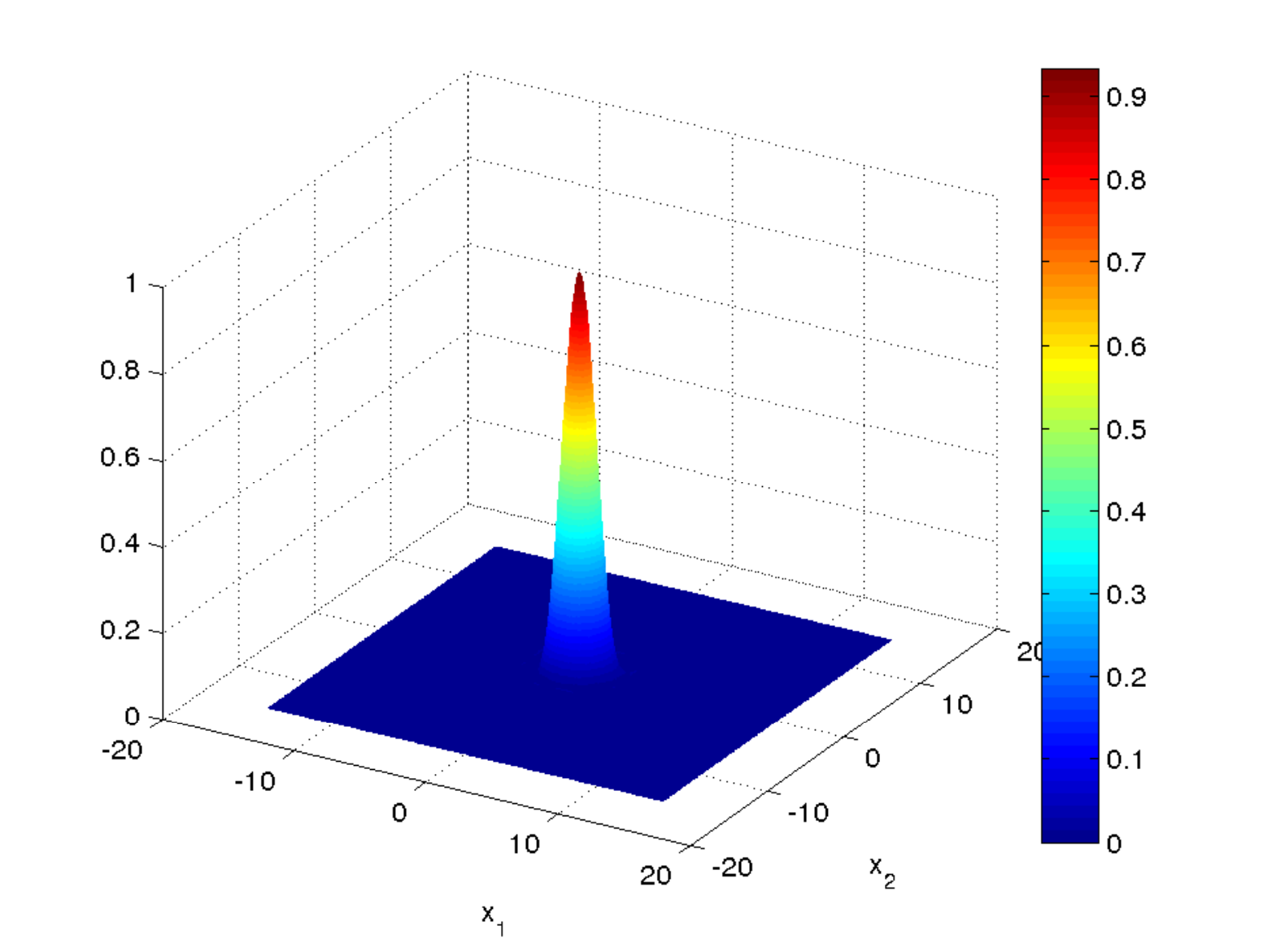}
\hspace*{1mm}\includegraphics[height=4cm, keepaspectratio]{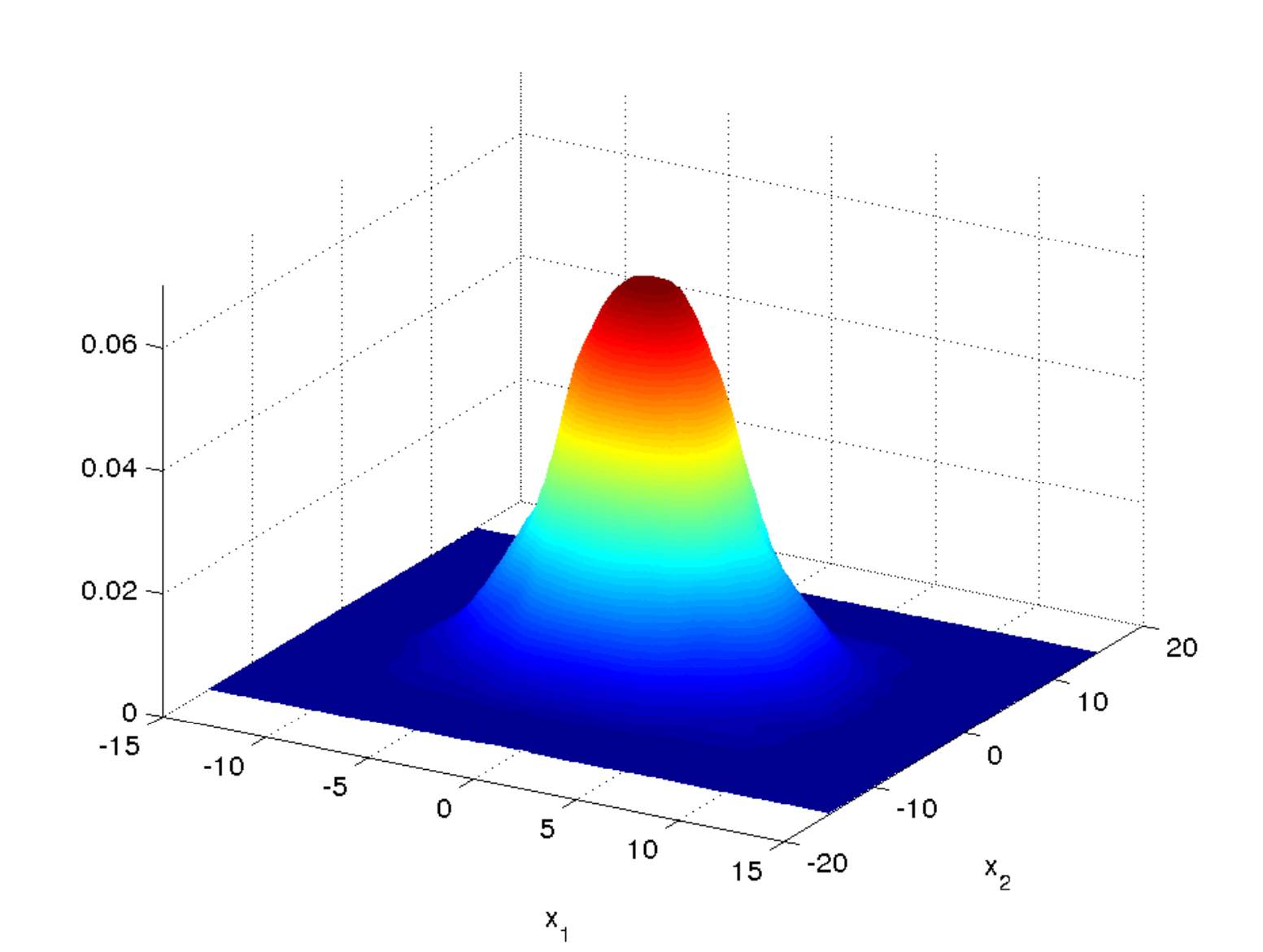}
\hspace*{1mm}\includegraphics[height=4cm, keepaspectratio]{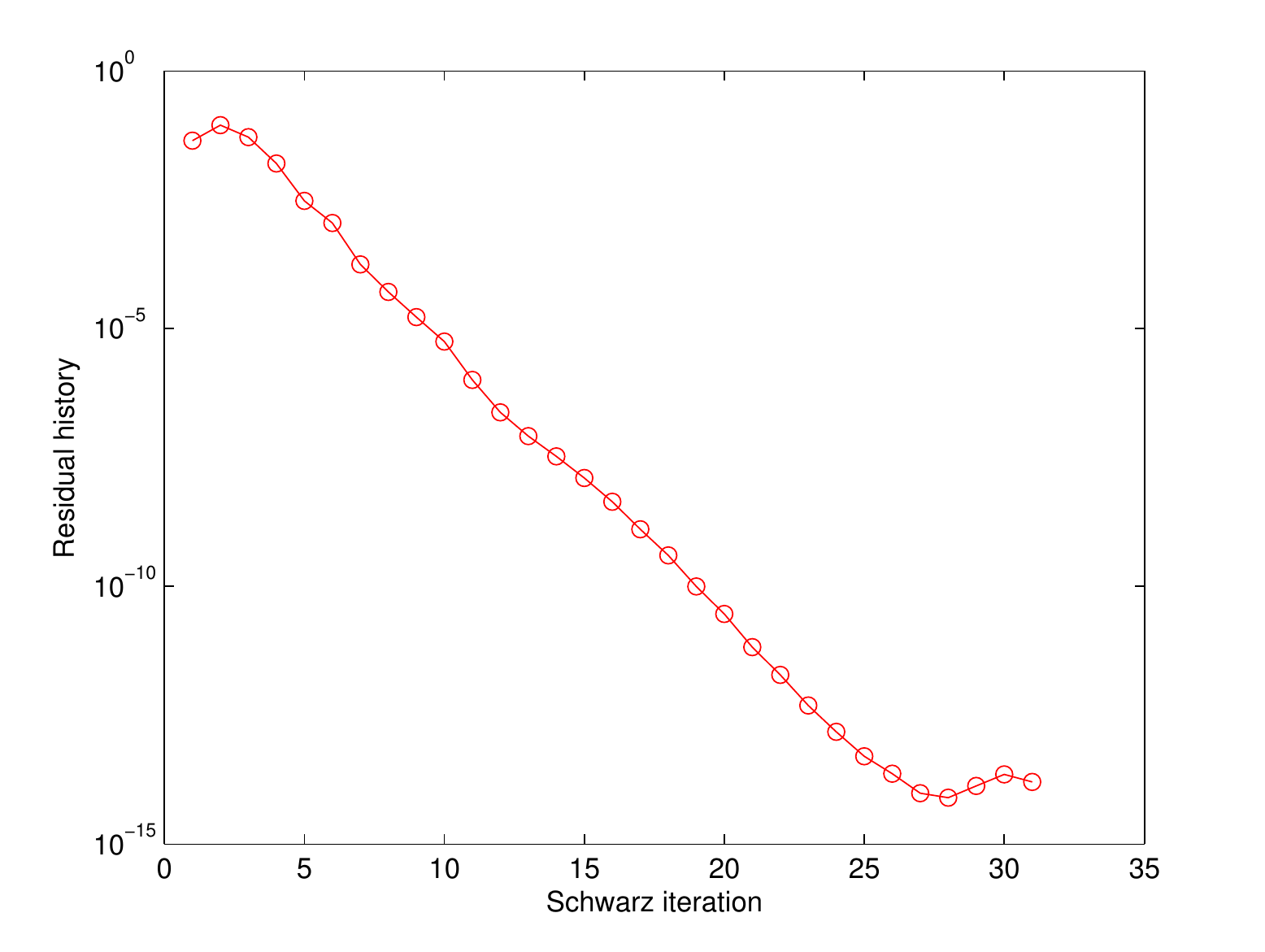}
\caption{$25$ subdomains. (Left) Initial data. (Middle) Reconstructed solution at time $T=16$. (Right) Residual history.}
\label{Heat1}
\end{center}
\end{figure}

A second preliminary test is dedicated to the computation by SWR-DDM with Robin-TC to the heat equation
\begin{eqnarray*}
\phi_t(x_1,x_2,t) -\triangle \phi(x_1,x_2,t)=0
\end{eqnarray*}
on $\Omega\times (0,T)$,  with the following initial data (see also Fig. \ref{Heat1} (Left))
\begin{eqnarray*}
\phi_0(x_1,x_2) = \exp\big(-0.2(x_1^2+x_2^2)\big).
\end{eqnarray*} 
The computational domain defined at the beginning of Section \ref{NumGauss}, and $\phi^{(0)}$ is constructed following the algorithm proposed in Subsection \ref{testA}. This is a relevant test, as the imaginary time method which will be implemented below is basically based on the solution of a normalized heat equation. The set-up is as follows. The global domain $(-15,15)$ is decomposed in $L^2=25$ subdomains, and the total number of grid points is $N^{(x_1)}\times N^{(x_2)}=201^2$. On each subdomain, a total of $N_{\phi}^2=36$ Gaussian local basis functions is constructed. The Robin-SWR algorithm is then implemented with $\mu=1$ defined in \eqref{RTC}, and we provide in Fig. \ref{Heat1} (Middle) the converged reconstructed solution at time $T=16$ where $\Delta t \approx 0.213$.
Here, we define the residual history as follows
\begin{eqnarray}\label{residueRT}
\textrm{Res}(k):=\Big(\int_0^{T}\sum_{i=1}^{L^2}\int_{\partial \Omega_i}|\phi_i^{(k)}-\phi_i^{(k-1)}|^2dx_1dx_2dt\Big)^{1/2}
\end{eqnarray}
and is reported in logscale as a function of the Schwarz iteration, $\big\{\big(k,\log(\textrm{Res}(k))\big), \, k \in \N\big\}$, in Fig. \ref{Heat1} (Right). Notice that for numerical convenience, the convergence criterion we use is a bit different from the one defined Section \ref{SWR}.
\subsection{Test 2.a: Ground state construction I}\label{testB}
\begin{figure}[!ht]
\begin{center}
\hspace*{1mm}\includegraphics[height=6.0cm, keepaspectratio]{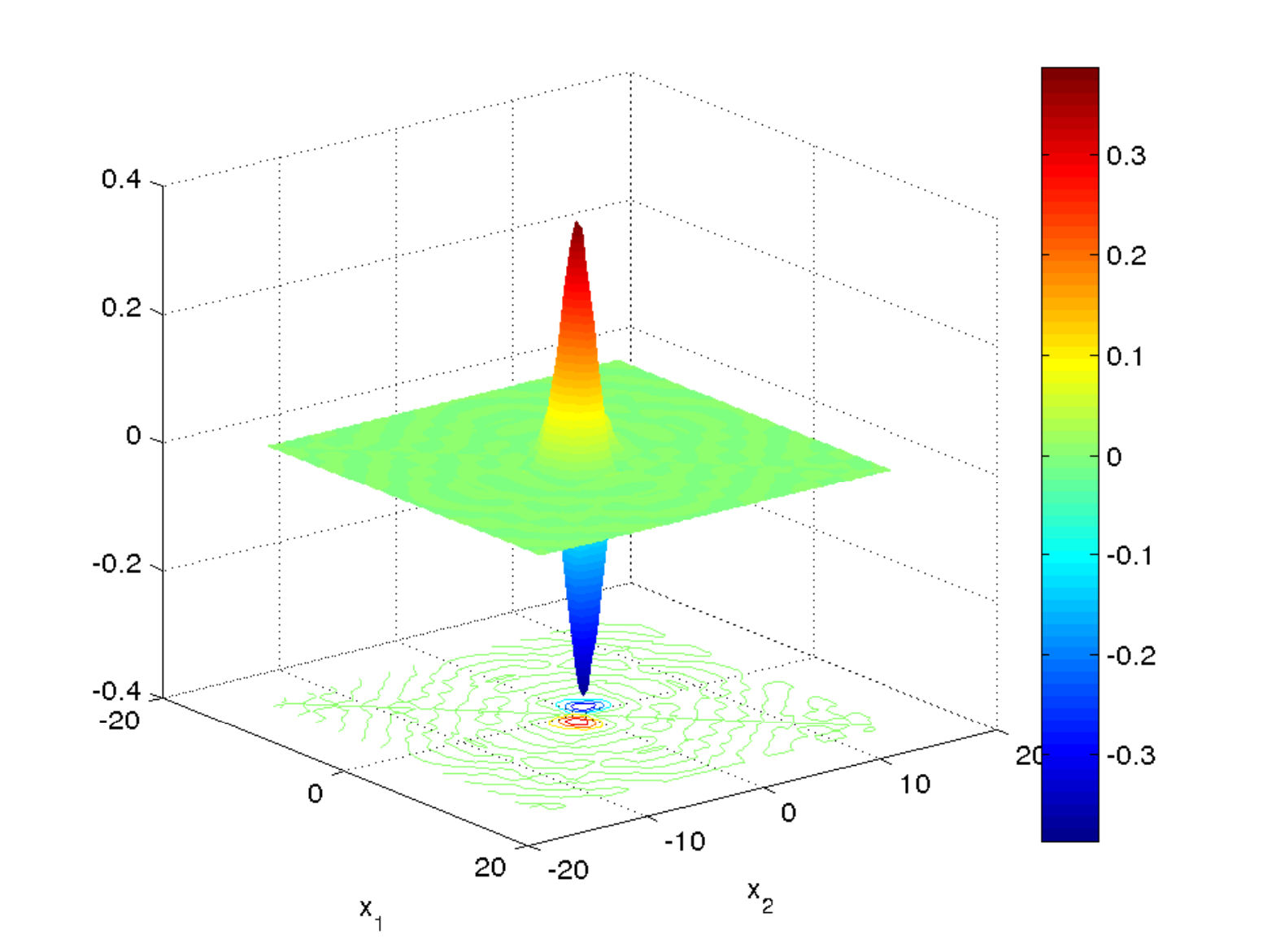}
\hspace*{1mm}\includegraphics[height=6.0cm, keepaspectratio]{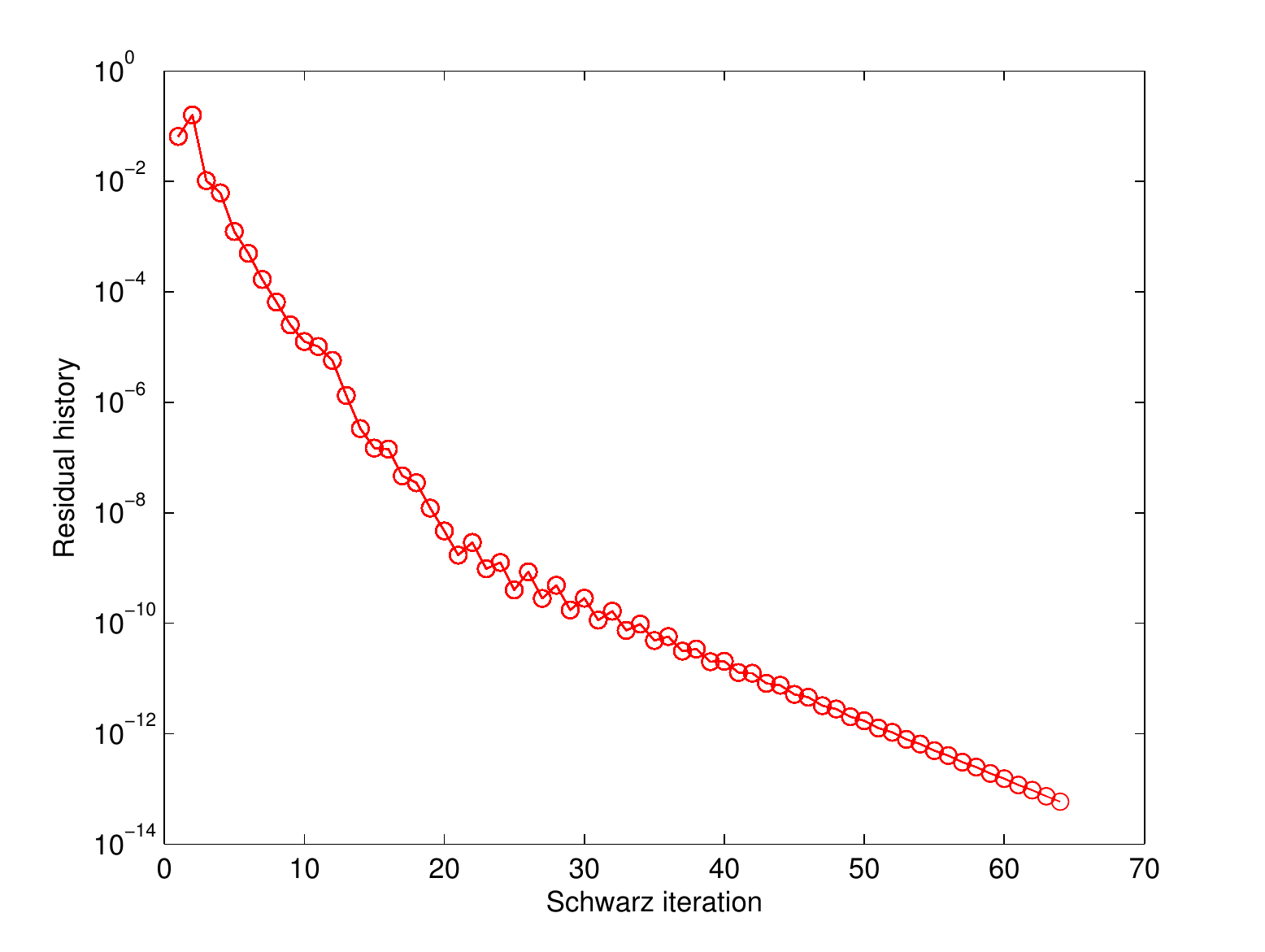}
\caption{$H_2$-molecule ground state: $25$ subdomains. (Left) Antisymmetric wavefunction. (Right) Residual error.}
\label{ground_Gauss}
\end{center}
\end{figure}

\begin{figure}[!ht]
\begin{center}
\hspace*{1mm}\includegraphics[height=6.0cm, keepaspectratio]{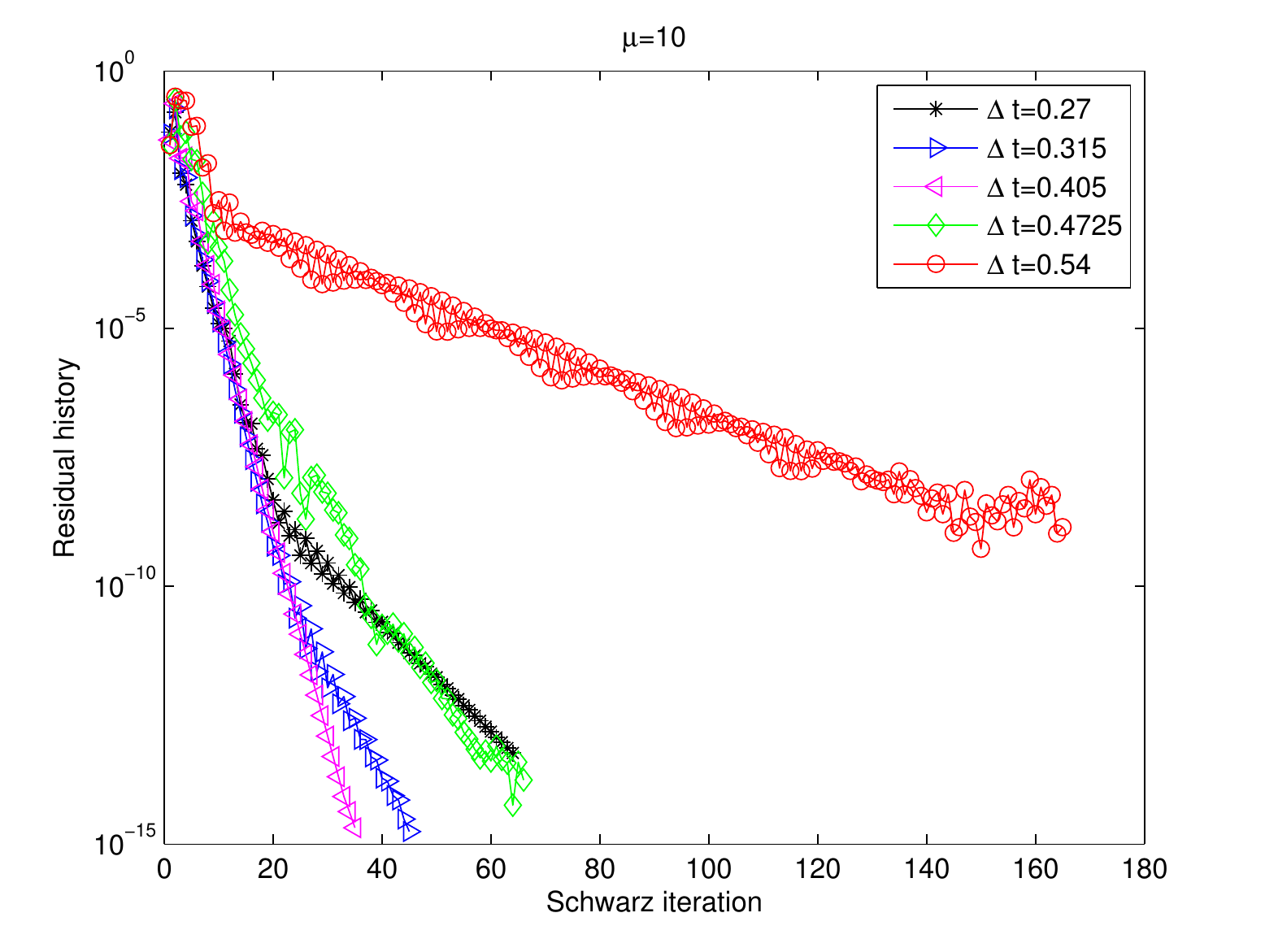}
\hspace*{1mm}\includegraphics[height=6.0cm, keepaspectratio]{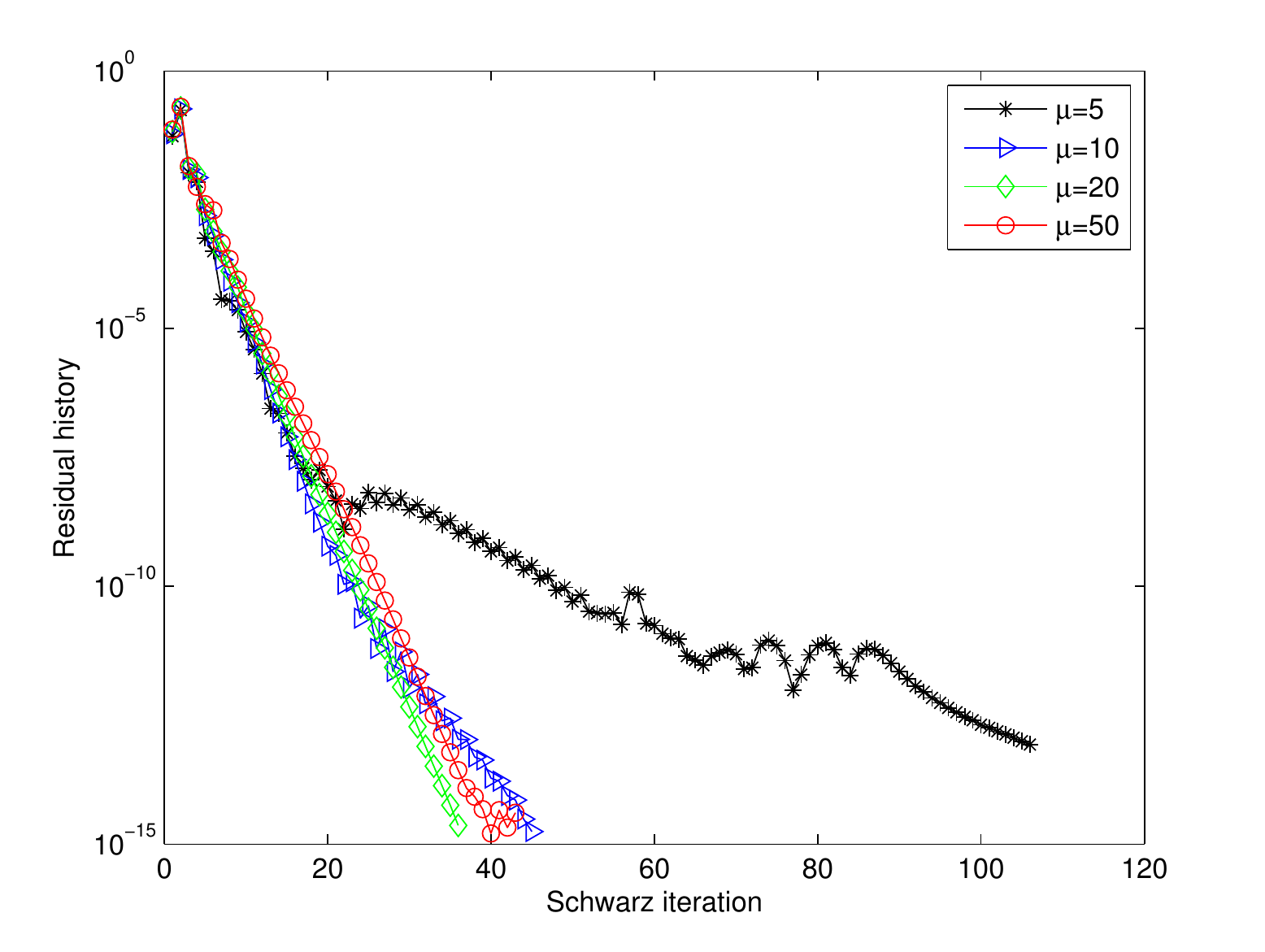}
\caption{Convergence rate as a function of the time step, for $\mu=10$. Convergence rate as a function of $\mu$, for $\Delta t=3.15\times 10^{-2}$.}
\label{DT}
\end{center}
\end{figure}

In this next experiment, we apply the imaginary time method for constructing the ground state of a 2-electron problem, more specifically a $H_2$-molecule with fixed nuclei. Again, the overall domain $(-15,15)$ is decomposed in $L^2=25$ subdomains. On each subdomain a total of $N_{\phi}^2=36$ Gaussian local basis functions, with $\delta=0.5$ in \eqref{gauss1d}, are used to construct the local solutions. The Robin-SWR algorithm described in Section \ref{SWR1} is applied with a Robin constant $\mu=10$ in \eqref{RTC}. Notice that a deeper analysis would help to select the optimized value of $\mu$ (that is ensuring the fastest convergence), see \cite{halpern2} in 1-d and for $2$ subdomains. At each Schwarz iteration $k$, we then solve $25$ independent imaginary time problems, from $0$ to $T^{\textrm{cvk},(k)}$ corresponding to the converged (imaginary) time of the imaginary time method. The wave transmission from one time iteration to another occurs through Robin transmission conditions. The test which is proposed here is as follows. The position of the nuclei is respectively $x_A=-1.25$ and $x_B=1.25$, while their charge is fixed to $Z_A=Z_B=1$. We use a regularized potential to avoid the singularity, thanks to the parameter $\eta=0.2$ in 
\begin{eqnarray}\label{pseudo}
V(x) = -1/\sqrt{(x-x_A)^2+\eta^2} -1/\sqrt{(x-x_B)^2+\eta^2},
\end{eqnarray}
Notice that the nuclei are located in the central subdomain. The total number of grid points is $N^{(x_1)}\times N^{(x_2)}=101^2$, and the overlap zone between 2 subdomains is $\approx 10\%$. We pick an initial guess as the following Gaussian function
\begin{eqnarray*}
\phi_0(x_1,x_2) = \exp\big(-(x_1^2+x_2^2)\big).
\end{eqnarray*}
At iteration $k$, we denote by $\widetilde{\phi}_g=\phi^{\textrm{cvg},(k)}$ the reconstructed solution at the converged time $T^{\textrm{cvg},(k)}$. We report in logscale Fig. \ref{ground_Gauss} (Right), the residual history as a function of the Schwarz iterations, $\big\{\big(k,\log(\textrm{Res}(k))\big), \, k \in \N\big\}$, where we numerically evaluate
\begin{eqnarray}\label{residueIT}
\textrm{Res}(k):=\Big(\int_0^{T^{\textrm{cvg},(k)}}\sum_{i=1}^{L^2}\int_{\partial \Omega_i}|\phi_i^{\textrm{cvg},(k)}-\phi_i^{\textrm{cvg},(k-1)}|^2dx_1dx_2dt\Big)^{1/2}.
\end{eqnarray}
The chosen time step is $\Delta t=4.5\times 10^{-1}$. Notice, that the constructed ground state $\widetilde{\phi}_g$, is not a priori antisymmetric. An a posteriori antisymmetrization of the reconstructed wavefunction $\widetilde{\phi}_g$, is possible thanks the operator $\mathcal{A}$ defined by:
\begin{eqnarray*}
\phi_g(x_1,x_2) = \mathcal{A}\widetilde{\phi}_{g}(x_1,x_2)=\left\{
\begin{array}{l}
\widetilde{\phi}_g(x_1,x_2),  \, x_2 \leq x_1,\\
-\widetilde{\phi}_g(x_1,x_2), \, x_2 > x_1
\end{array}
\right.
\end{eqnarray*}
In Fig. \ref{ground_Gauss} (Middle) we report the antisymmetrized computed energy state. Notice however, that this antisymmetric wavefunction is, a priori, associated to an eigenenergy higher than the groundstate energy.\\
We are now interested in the residual history as a function of the time step, for $\mu=10$. We basically observe that the smaller the time step, the faster the convergence of the SWR algorithm. This is coherent with \cite{halpern3,lorin-TBS}, where it is proven that the convergence rate for CSWR (based on Dirichlet boundary conditions) is exponential in $-1/\sqrt{\Delta t}$ for a one-dimensional two-domain problem. For the Robin-SWR algorithm however, the convergence rate is dependent on the choice of $\mu$, as we can observe in Fig. \ref{DT} (Left). 
We also represent in Fig. \ref{DT} (Right), the residual history as a function of $\mu$. Notice that in \cite{halpern3} was established for Robin-SWR, a criterion to optimize the value of $\mu$ for the one-dimensional linear advection diffusion reaction equation  with constant coefficients. In the presented simulations, we can not really conclude about an optimized choice of $\mu$, as a fine mathematical analysis would be required to extend ideas from \cite{halpern3}. The proposed set-up, multidimensional with non-constant potentials, would make the analysis quite complex and quite beyond the current state of the art in this research field.
\subsection{Test 2.b: Ground state construction II}\label{testB}
\begin{figure}[!ht]
\begin{center}
\hspace*{1mm}\includegraphics[height=4cm, keepaspectratio]{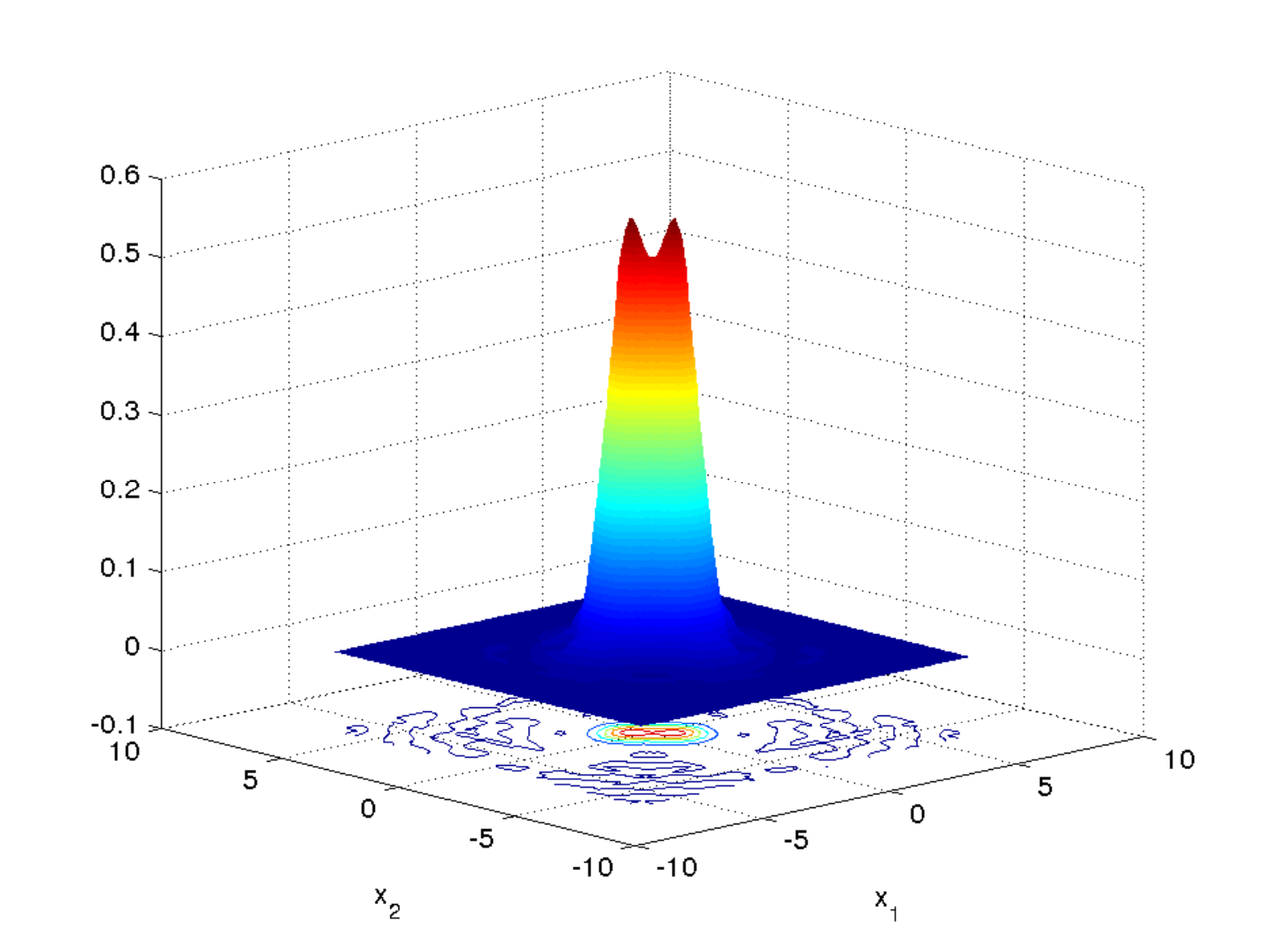}
\hspace*{1mm}\includegraphics[height=4cm, keepaspectratio]{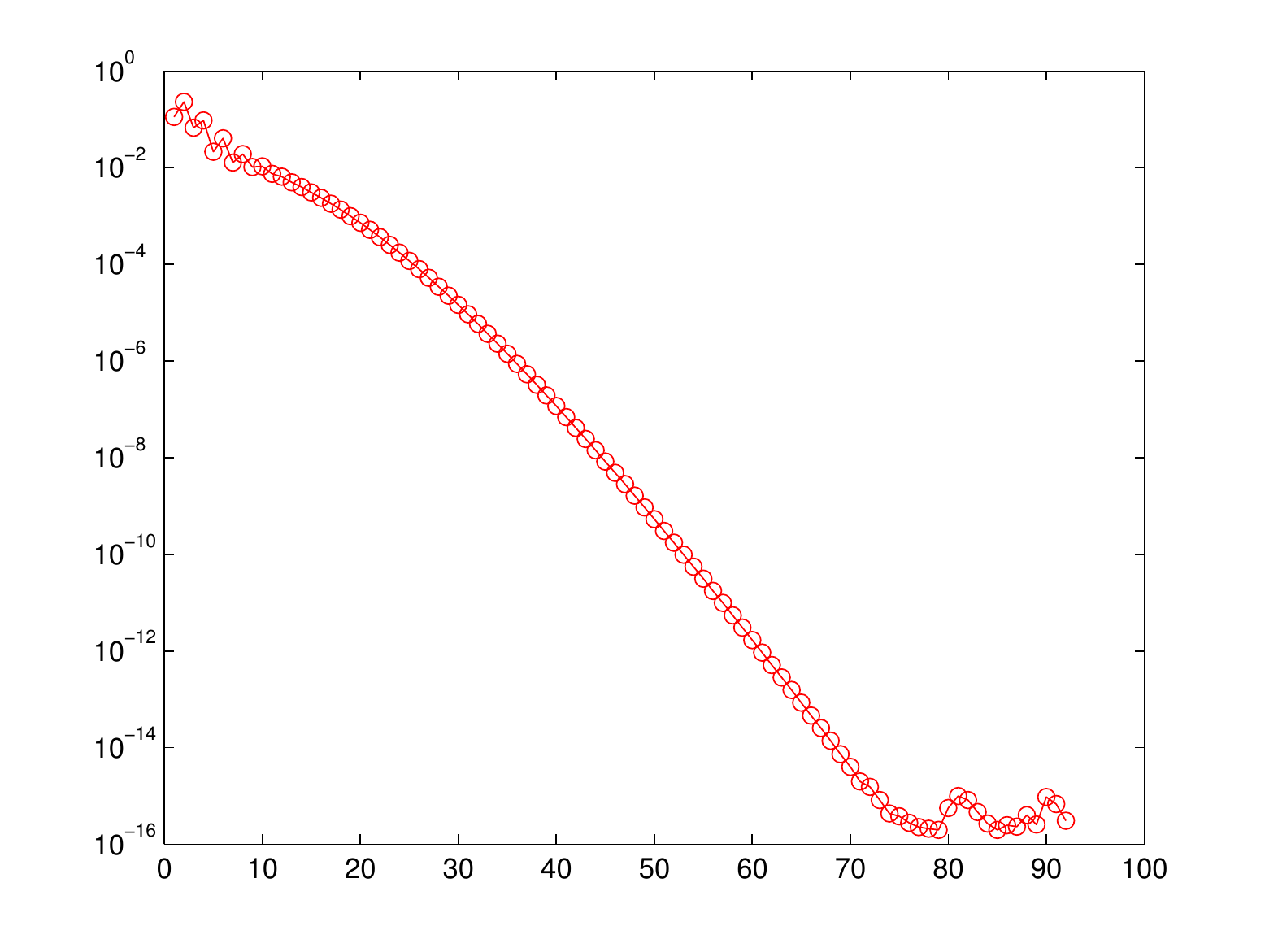}
\caption{Converged ground state: $25$ subdomains. (Left) Modulus of the converged ground state. (Right) Residual error.}
\label{ground_Gauss2}
\end{center}
\end{figure}

In this next test, we again apply the imaginary time method for constructing the ground state of the 2-electron problem. In this test however, the overall computational domain is smaller and the ground state support will cover several subdomains, slowing down the DDM convergence. The overall domain $(-6.5,6.5)$ is decomposed in $L^2=25$ subdomains, with Gaussian basis functions as defined in the previous section. On each subdomain a total of $N_{\phi}^2=36$ Gaussian local basis functions, with $\delta=2$ in \eqref{gauss1d}, are used to construct the local solutions. The chosen Robin constant is $\mu=10$  in \eqref{RTC}.  The nucleus singularities are still located in the central subdomain, respectively in $x_A=-0.5$ and $x_B=0.5$ and the charge while their charge is $Z_A=Z_B=1$. We again use a regularized potential to avoid the singularity. The total number of grid points is $N^{(x_1)}\times N^{(x_2)}=151^2$, the overlap zone between 2 subdomains is $\approx 10\%$. We pick an initial guess as the following Gaussian function $\phi_0=\widetilde{\phi}_0/\|\widetilde{\phi}_0\|_0$, where
\begin{eqnarray*}
\widetilde{\phi}_0(x_1,x_2) = \exp\big(-0.8(x_1^2+x_2^2)\big).
\end{eqnarray*}
We report in Fig. \ref{ground_Gauss2} (Left), the converged ground state obtained by the Robin-SWR method. In addition, we report in logscale Fig. \ref{ground_Gauss2} (Right), the residual history as a function of the Schwarz iterations, $\big\{\big(k,\log(\textrm{Res}(k))\big)\big\}$, defined in \eqref{residueIT}. The chosen time step is $\Delta t=10^{-2}$.
\subsection{Test 3: Real time experiment}
\begin{figure}[!ht]
\begin{center}
\hspace*{1mm}\includegraphics[height=6cm, keepaspectratio]{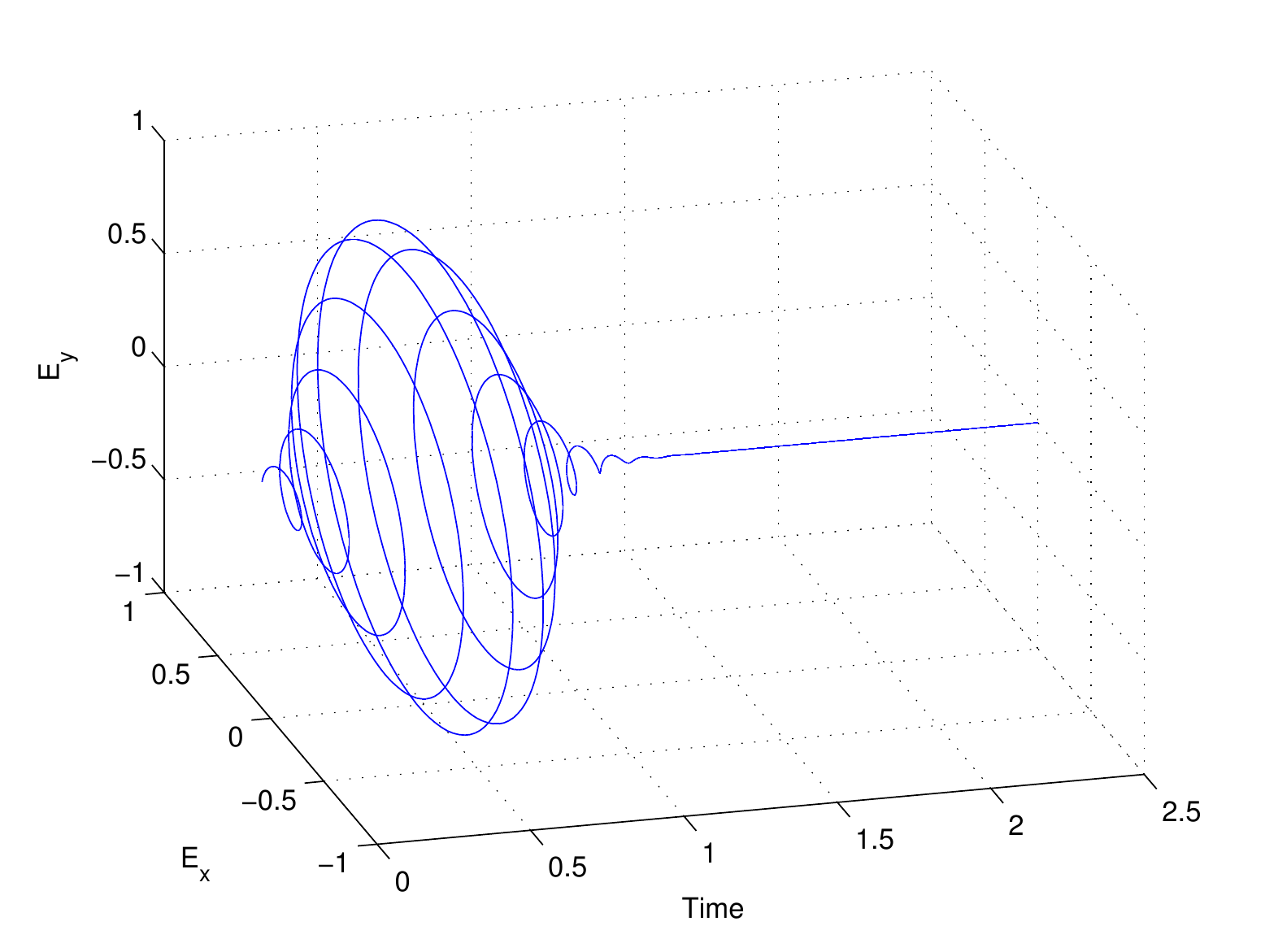}
\hspace*{1mm}\includegraphics[height=6cm, keepaspectratio]{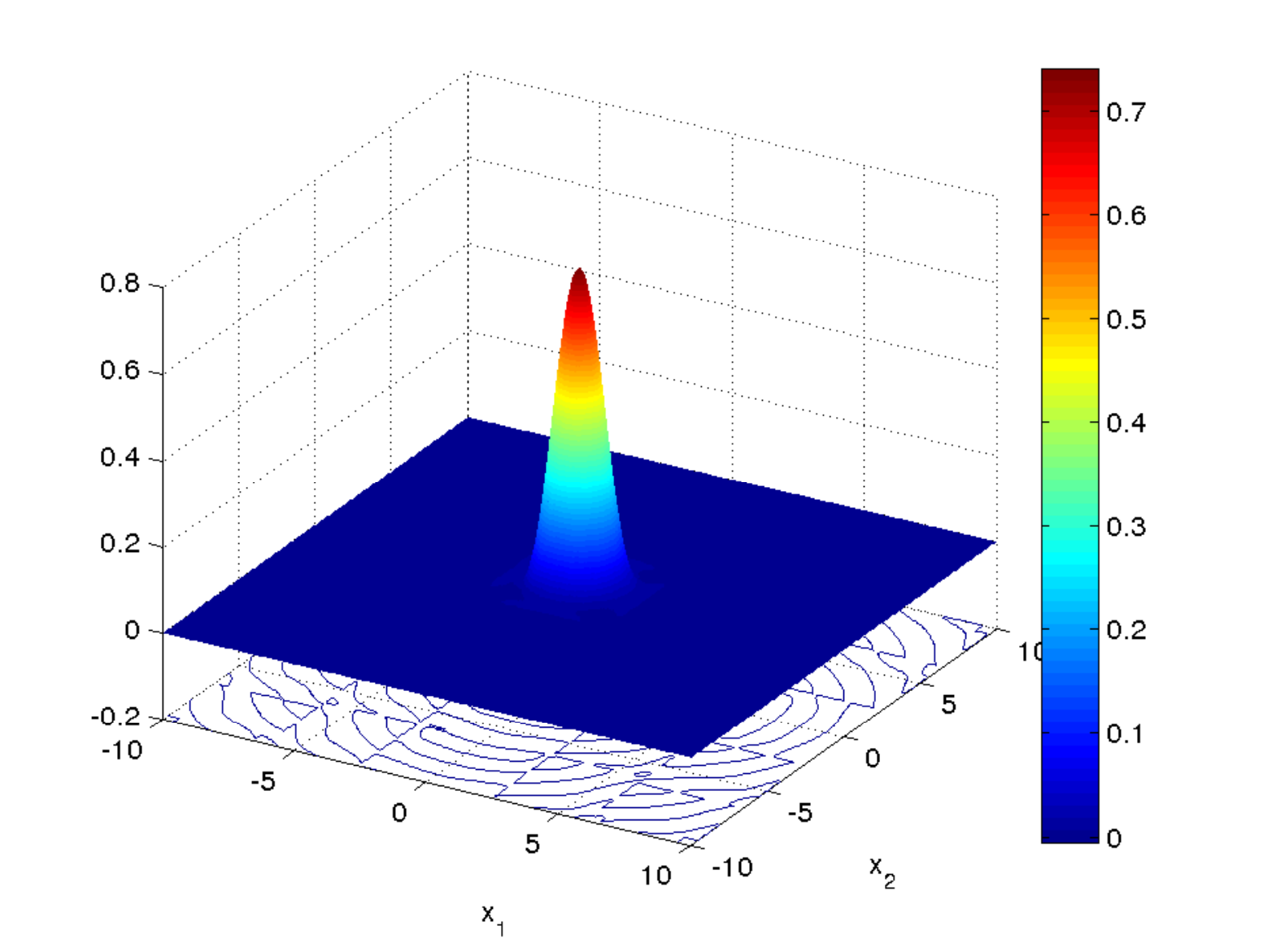}
\caption{(Left) Circular laser electric field from time $0$ to $T=2.5$. (Right) Initial state for TDSE.}
\label{laser_circ}
\end{center}
\end{figure}

\begin{figure}[!ht]
\begin{center}
\hspace*{1mm}\includegraphics[height=4cm, keepaspectratio]{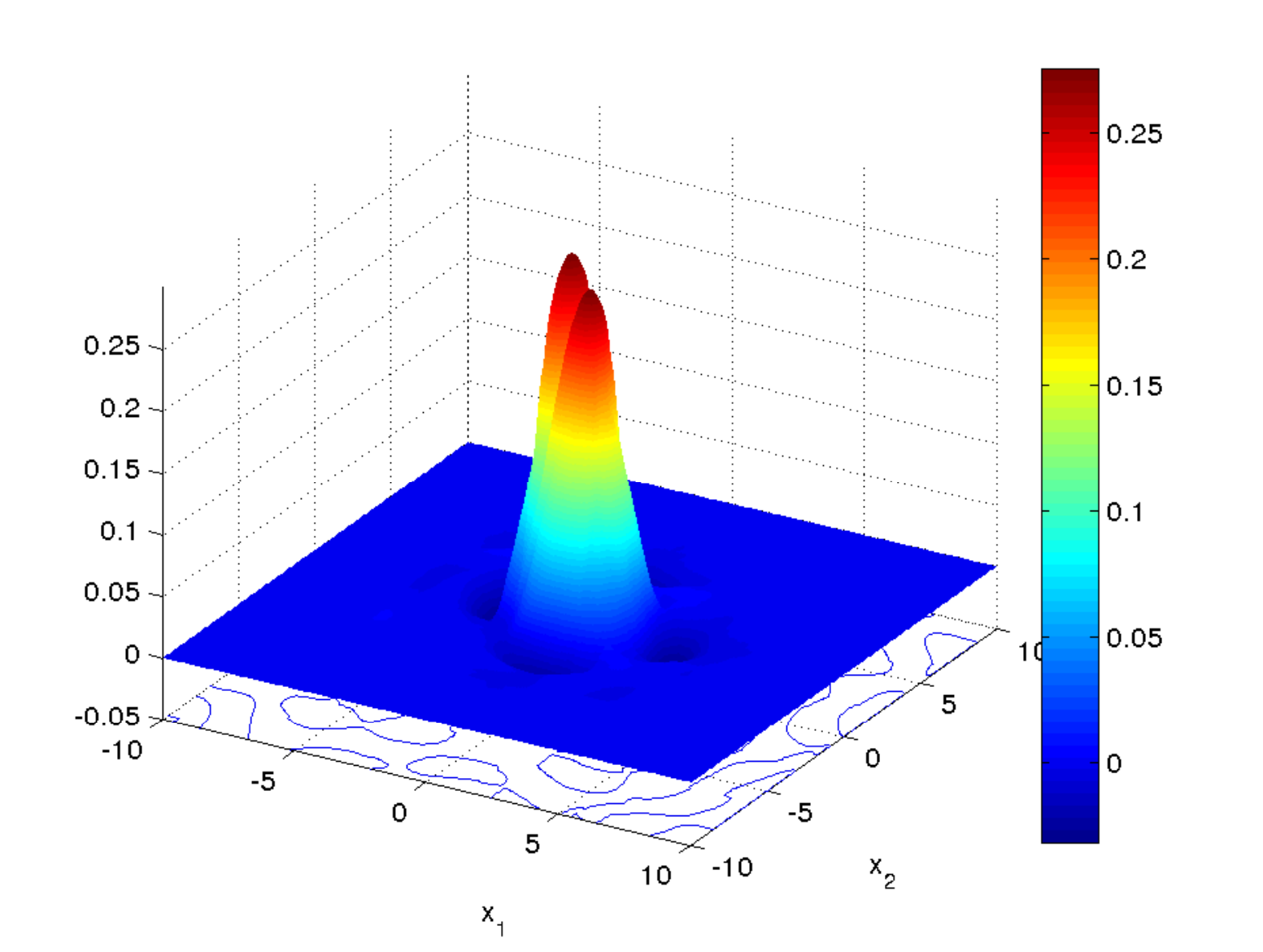}
\hspace*{1mm}\includegraphics[height=4cm, keepaspectratio]{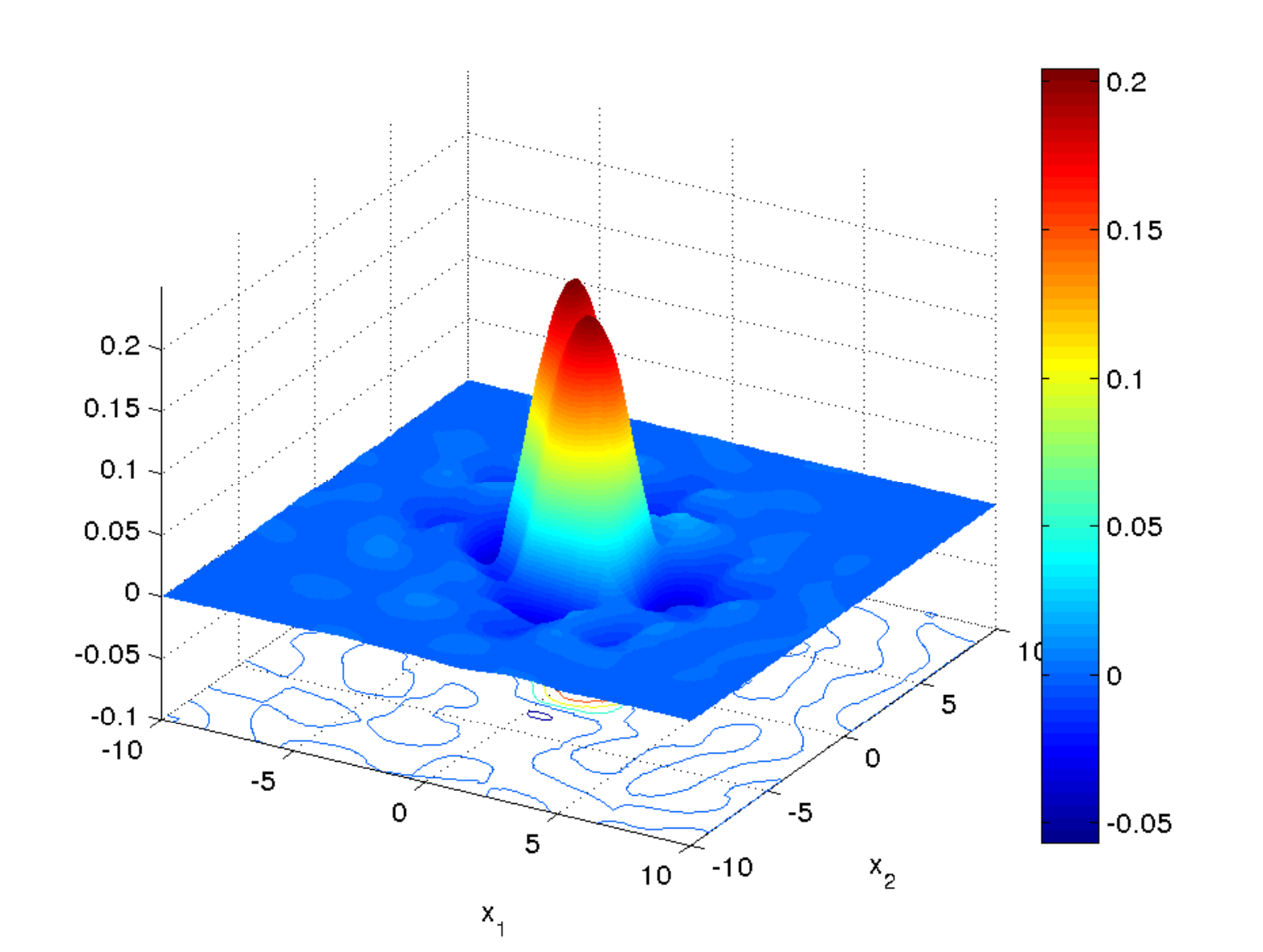}
\hspace*{1mm}\includegraphics[height=4cm, keepaspectratio]{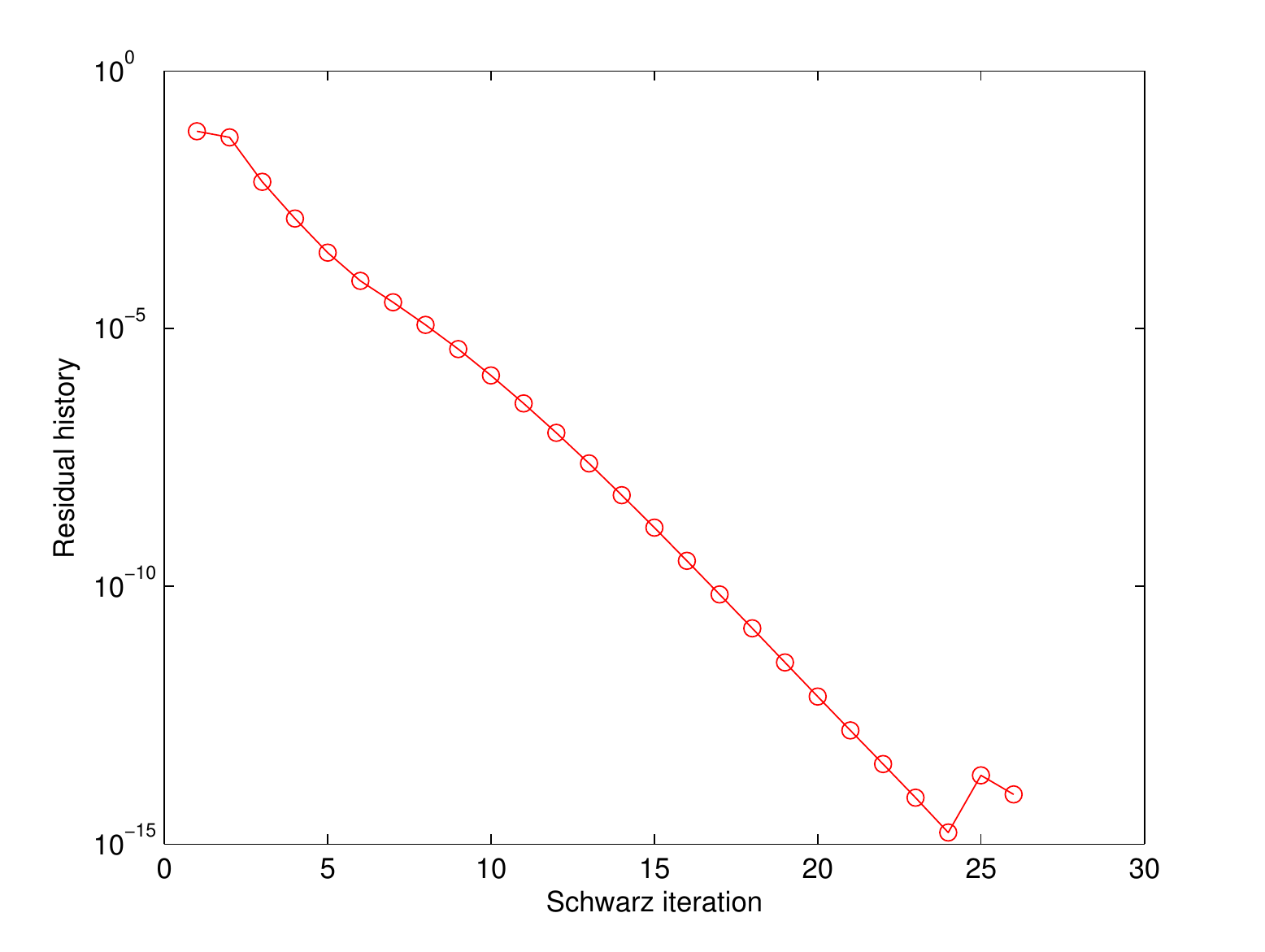}
\caption{$H_2$-molecule subject to electric field at time $T=2.5$: $25$ subdomains. (Left). Imaginary part of solution at final time after the first Schwarz iteration. (Middle) Solution at final time at Schwarz iteration $k=k^{(\textrm{cvg})}$. (Right) Residual error.}
\label{CVTDSE2}
\end{center}
\end{figure}

The following test is devoted to the evolution of a 2-electron wavefunction, subject to an external circular electric ${\bf E}(t)=(E_x(t),E_y(t))$ defined by:
\begin{eqnarray*}
\left\{
\begin{array}{lcl}
E_x(t) & = & E_0\cos(\omega_0t)\exp\big(-\nu_0(T/2-t)^2\big), \\
E_y(t) & = & E_0\sin(\omega_0t)\exp\big(-\nu_0(T/2-t)^2\big)
\end{array}
\right.
\end{eqnarray*}
where $E_0=1$ $\omega_0=8$, $\nu_0=10$, $T=2.5$, see Fig. \ref{laser_circ} (Left). We are interested in the convergence of the Robin-SWR algorithm. 
 As before we assume that the global domain $\Omega=(-10,10)^2$ is decomposed in $L^2=25$ square subdomains, with an overlap corresponding to $10\%$. We fix to $N^2_{\phi}=36$, the number of the Gaussian local basis functions, and in \eqref{RTC}, we impose $\mu=-10{\tt i}$. The total number of grid points is $N^{(x_1)}\times N^{(x_2)}=201^2$. The time step is given by $\Delta t=5\times 10^{-2}$ and final time $T=2.5$. The initial data is a Gaussian wave Fig. \ref{laser_circ} (Right). The potential \eqref{pseudo} is regularized using a parameter $\eta=0.5$.
\begin{eqnarray*}
\phi_0(x_1,x_2) = \exp\big(-(x_1^2+x_2^2)\big).
\end{eqnarray*} 
We report the solution at time $T=2.5$, at the end of the first Schwarz iteration the imaginary part of the wavefunction $\mathcal{I}\psi^{(1)}$ ($k=1$) in Fig. \ref{CVTDSE2} (Left). We represent in Fig. \ref{CVTDSE2} (Middle), the imaginary part of the converged solution $\mathcal{I}\psi^{(k^{(\textrm{cvg})})}$, $k=k^{(\textrm{cvg})}$ at time $T=2.5$. The residual history $\big\{\big(k,\log(\textrm{Res}(k))\big), \, k \in \N\big\}$ is represented in Fig. \ref{CVTDSE2} (Right).
%
\section{Numerical experiments: Local Slater's Determinants}\label{NumSlater}
This section is devoted to some numerical experiments in imaginary time  with local Slater's determinant basis functions, with $d=1$ and $N=2$. The geometry and domain decomposition is identical to Section \ref{NumGauss}.  For realistic applications, and as discussed above, an appropriate choice of $K_i$ as a function of the position of the nuclei will be useful in order to accurately reduce the overall computational complexity of the method. A forthcoming paper will be dedicated to some exhaustive experiments in real time.\\
The set-up is the same as above, except that the local basis functions are here assumed to be local Slater's determinants, constructed from $1$-electron orbitals, see Section \ref{1D-2E}. 
\subsection{Test 1: Local Slater Determinants construction}\label{test1}
\begin{figure}[!ht]
\begin{center}
\hspace*{1mm}\includegraphics[height=4cm, keepaspectratio]{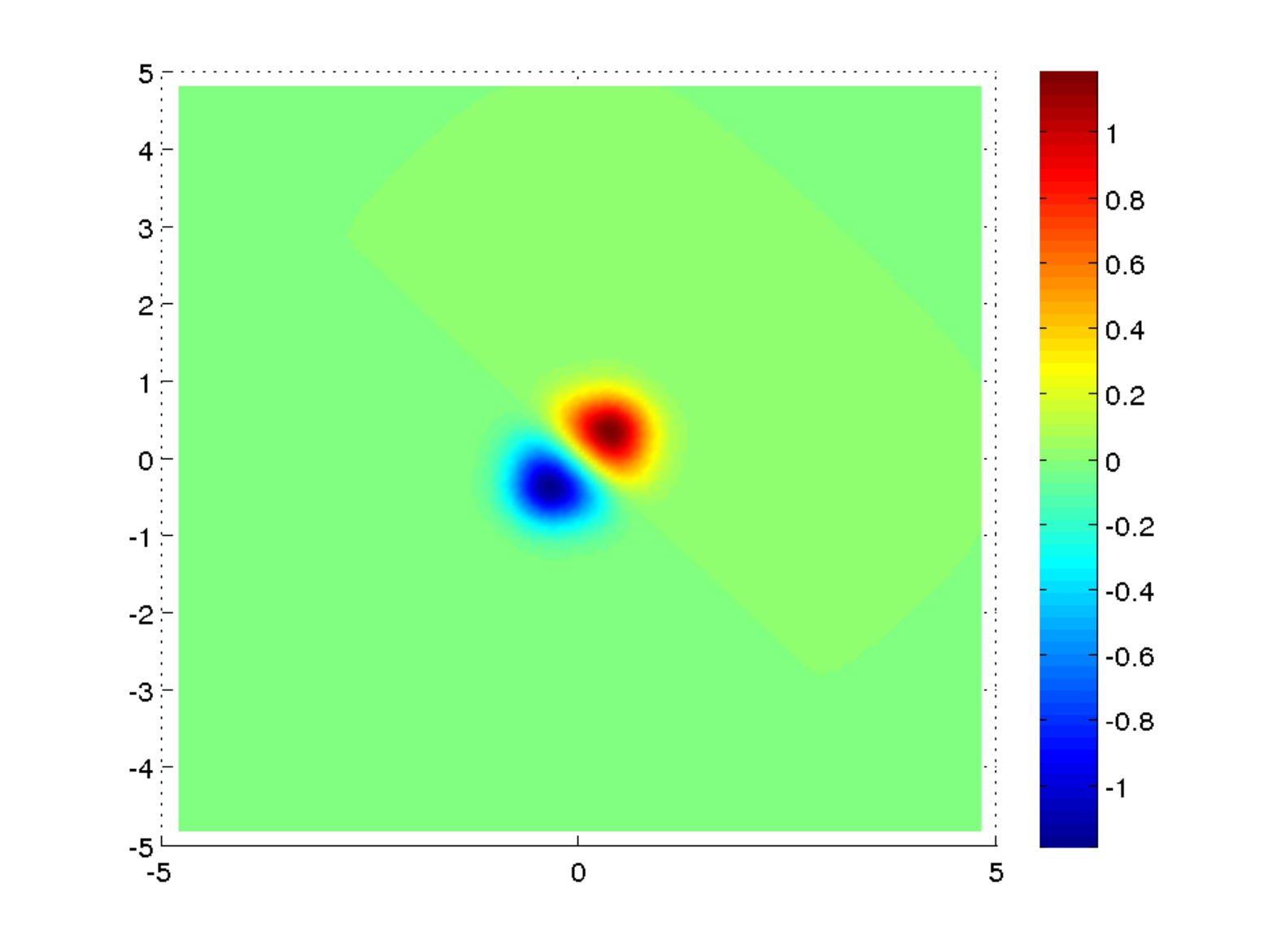}
\hspace*{1mm}\includegraphics[height=4cm, keepaspectratio]{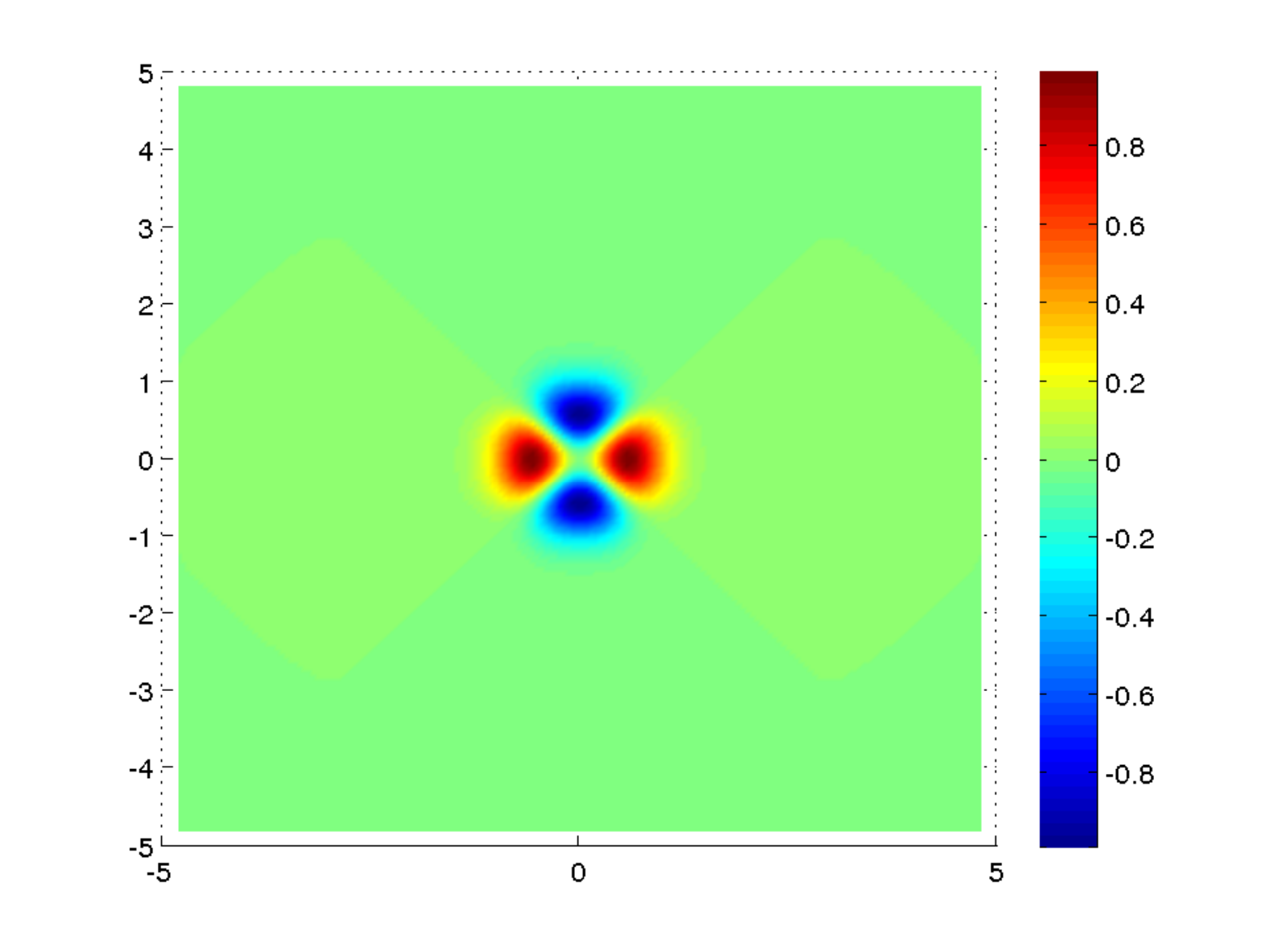}
\hspace*{1mm}\includegraphics[height=4cm, keepaspectratio]{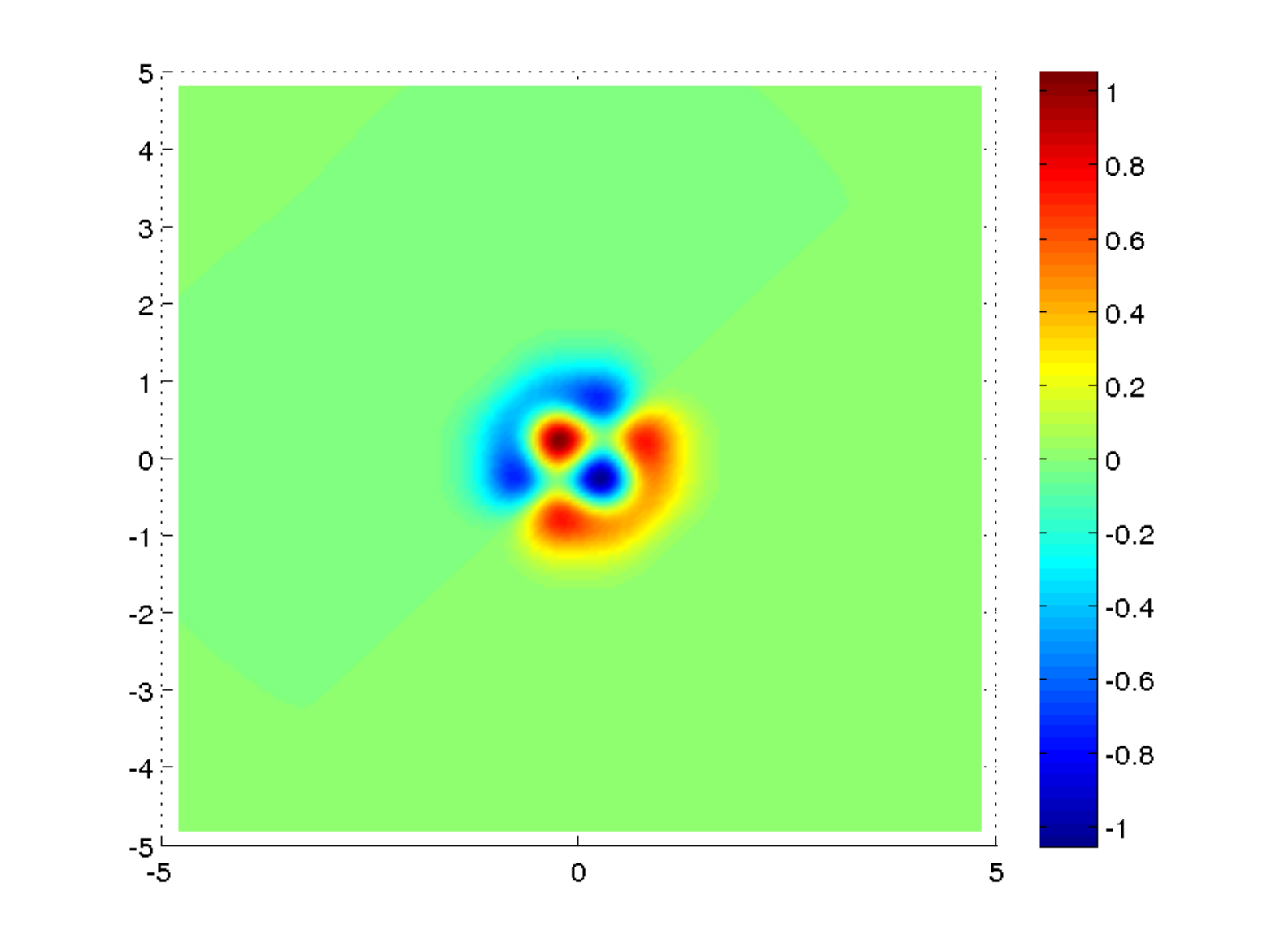}
\hspace*{1mm}\includegraphics[height=4cm, keepaspectratio]{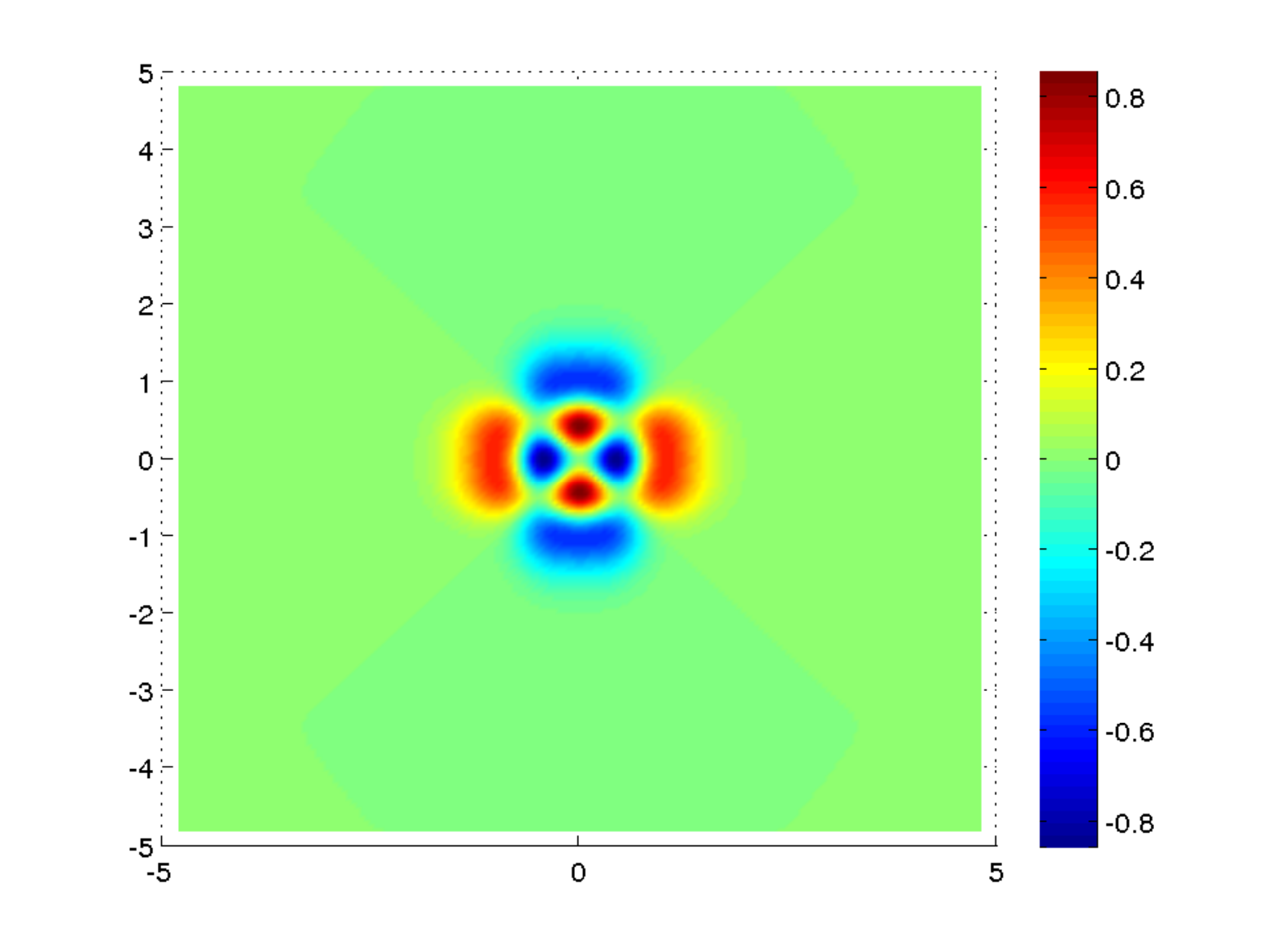}
\hspace*{1mm}\includegraphics[height=4cm, keepaspectratio]{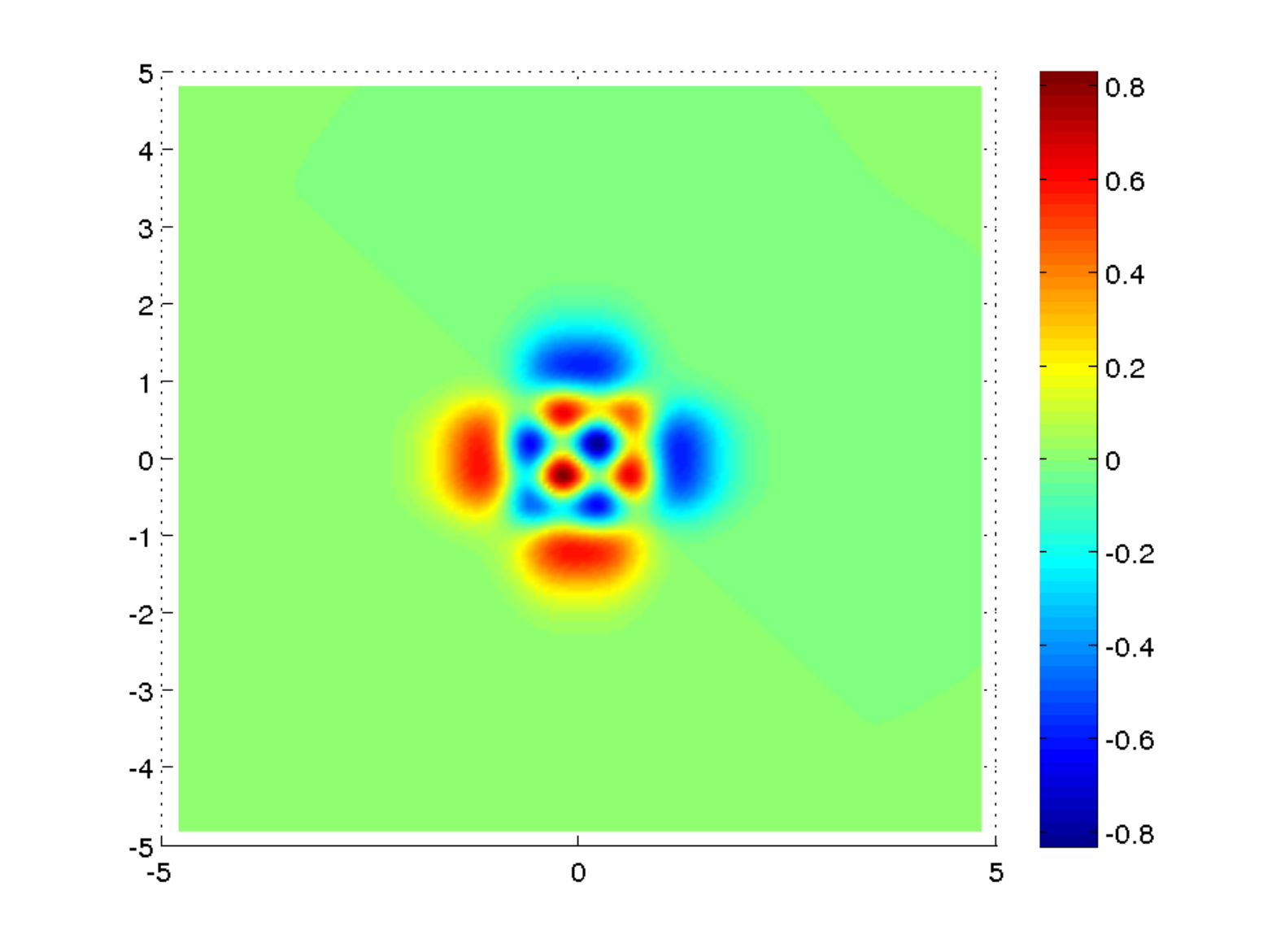}
\hspace*{1mm}\includegraphics[height=4cm, keepaspectratio]{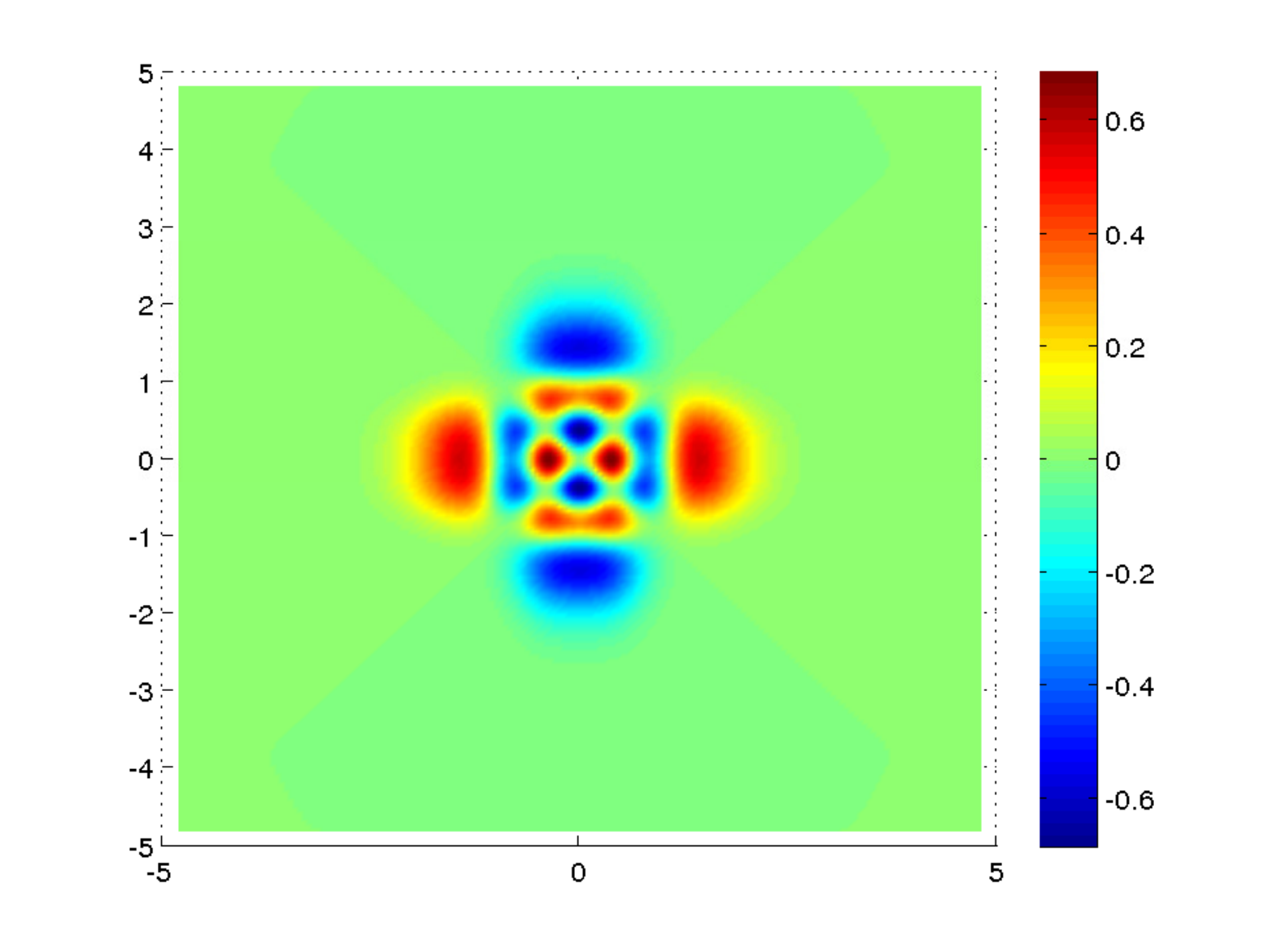}
\caption{First $6$ LSD's in $\Omega_{13}=[a_3-\epsilon^{(x_1)}/2,b_3+\epsilon^{(x_1)}/2]\times[c_3-\epsilon^{(x_2)}/2,d_3+\epsilon^{(x_2)}/2]$.}
\label{SDBF_3x3}
\end{center}
\end{figure}

This first test is dedicated to the construction of the LSD's for a 2-nucleus problem, with charge $Z_A=Z_B=1$. We choose $K=N_{\phi}(N_{\phi}+1)/2=45$, $L=5$ (for a total of $25$ subdomains $\Omega_1, \cdots,\Omega_{25}$) and $a=c=-8$, $b=d=8$ and $\epsilon^{(x_{1,2})}=3.2\times 10^{-1}$. The 2 nuclei are located in the central subdomain $\Omega_{13}=[a_3-\epsilon^{(x_1)}/2,b_3+\epsilon^{(x_1)}/2]\times[c_3-\epsilon^{(x_2)}/2,d_3+\epsilon^{(x_2)}/2]$ with $a_3-\epsilon^{(x_1)}/2=c_3-\epsilon^{(x_2)}/2=1.6$. We choose $\epsilon=0.1$ and $\sigma_1=315/256$ in \eqref{Beps}, and $N^{(x_1)}=N^{(x_2)}=301$. As an illustration, we represent the first 6 local Slater determinants (LSD's) in $\Omega_{13}=[a_3-\epsilon^{(x_1)}/2,b_3+\epsilon^{(x_1)}/2]\times[c_3-\epsilon^{(x_2)}/2,d_3+\epsilon^{(x_2)}/2]$ in Fig. \ref{SDBF_3x3}, as well as the first 6 LSD's in $\Omega_{14}=[a_4-\epsilon^{(x_1)}/2,b_4+\epsilon^{(x_1)}/2]\times[c_3-\epsilon^{(x_2)}/2,d_3+\epsilon^{(x_2)}/2]$. Recall in practice, that there is no need to explicitly construct these LSD's, see \cite{CAM15-10}.
\begin{figure}[!ht]
\begin{center}
\hspace*{1mm}\includegraphics[height=4cm, keepaspectratio]{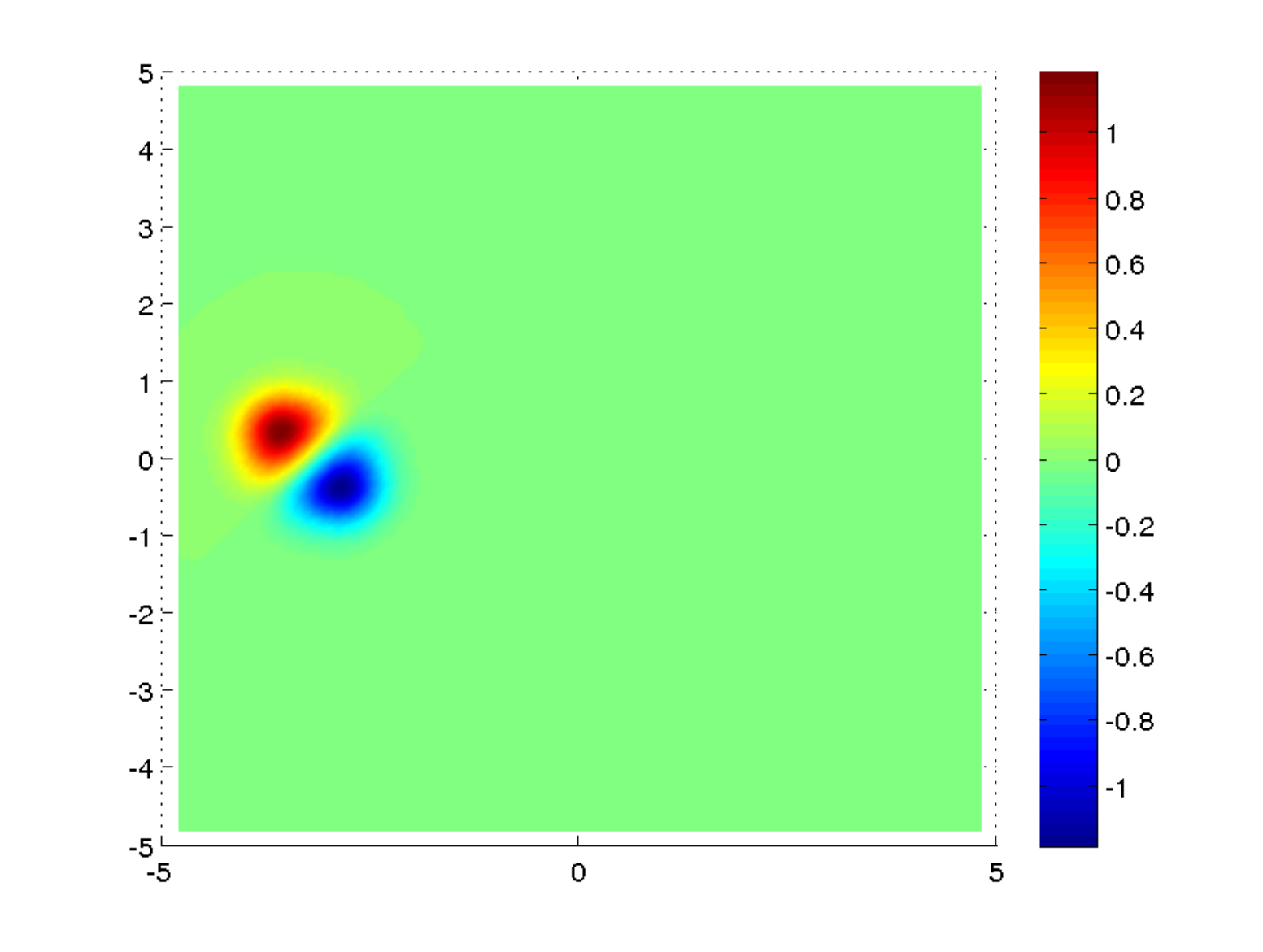}
\hspace*{1mm}\includegraphics[height=4cm, keepaspectratio]{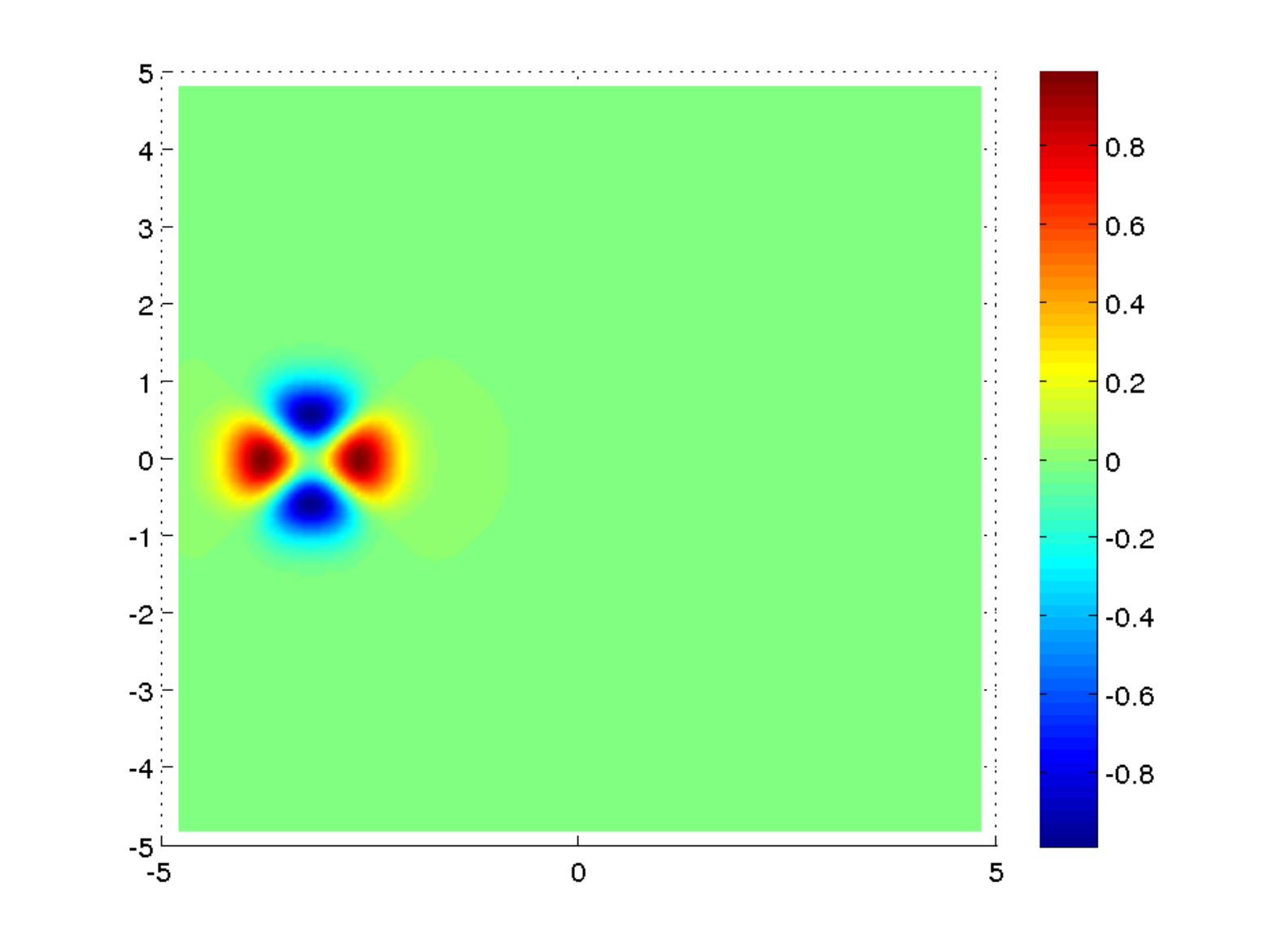}
\hspace*{1mm}\includegraphics[height=4cm, keepaspectratio]{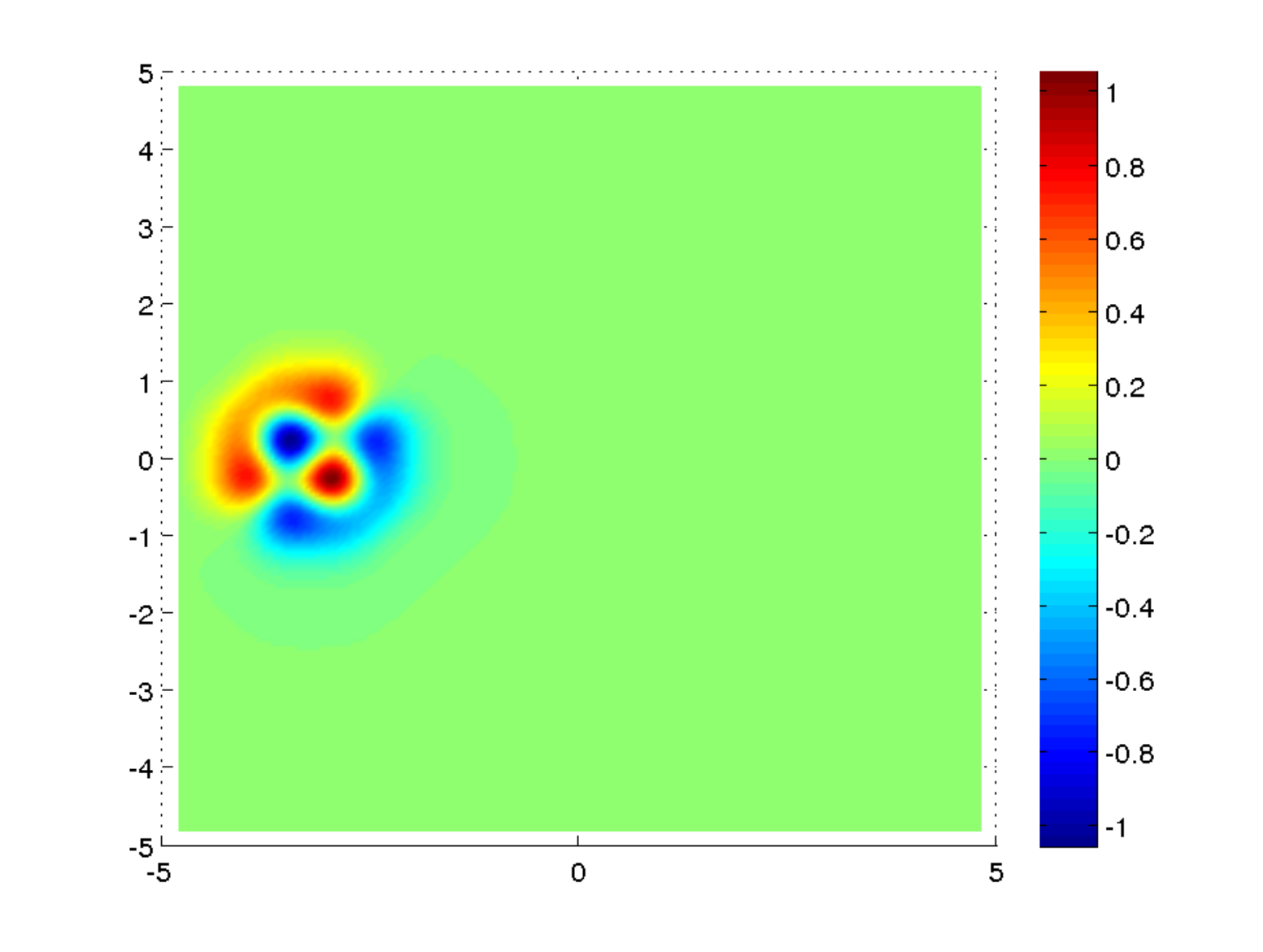}
\hspace*{1mm}\includegraphics[height=4cm, keepaspectratio]{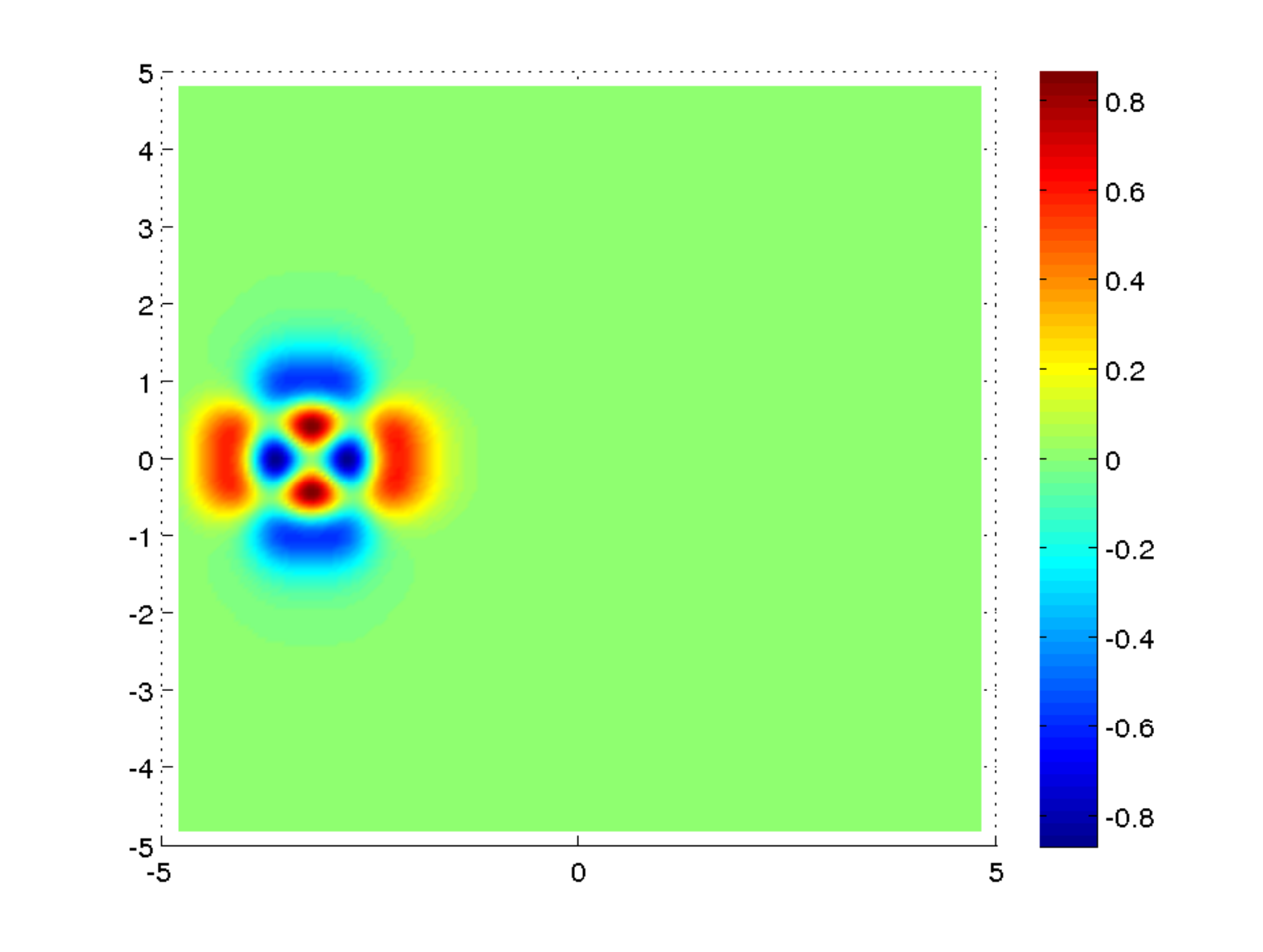}
\hspace*{1mm}\includegraphics[height=4cm, keepaspectratio]{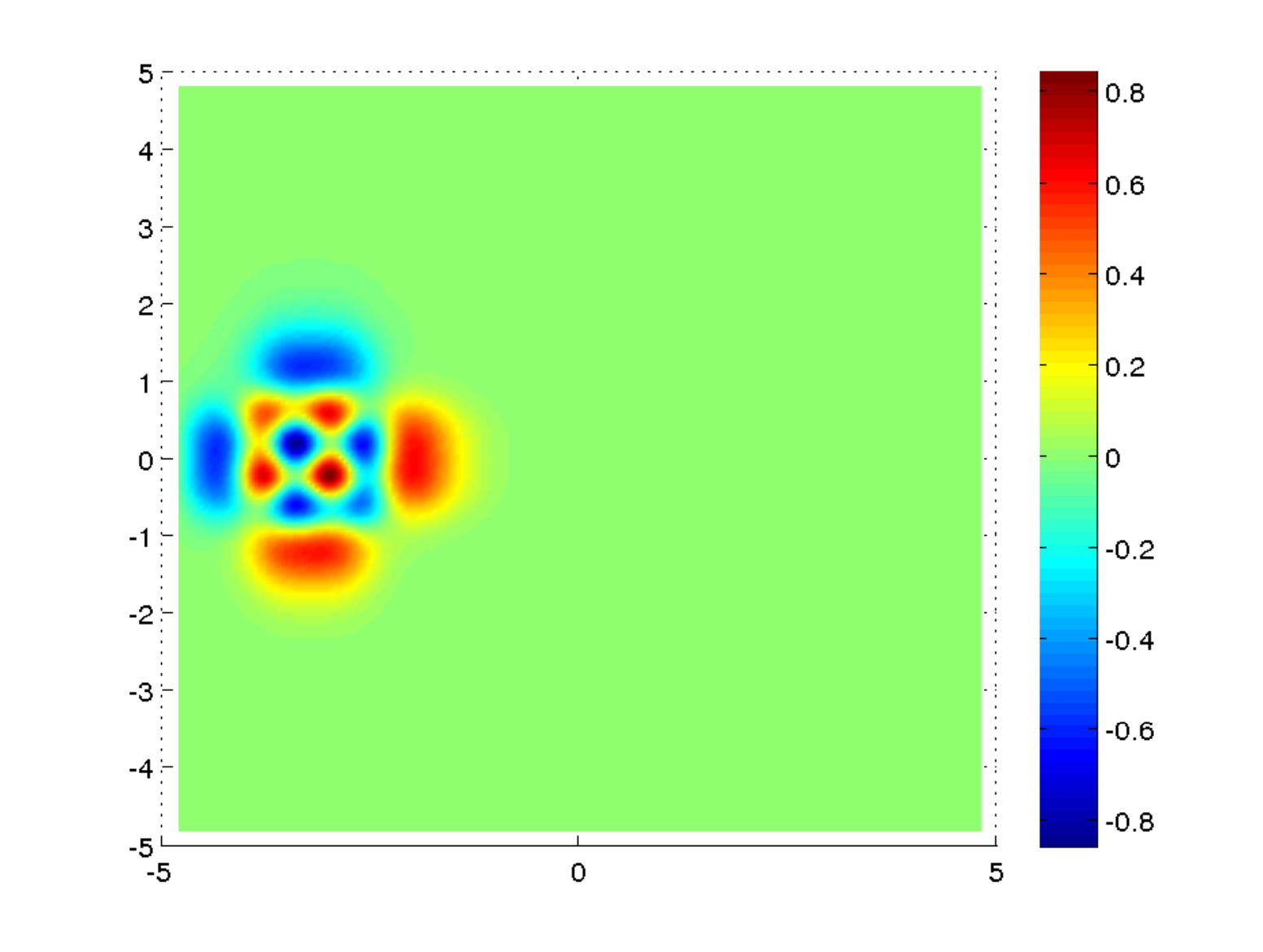}
\hspace*{1mm}\includegraphics[height=4cm, keepaspectratio]{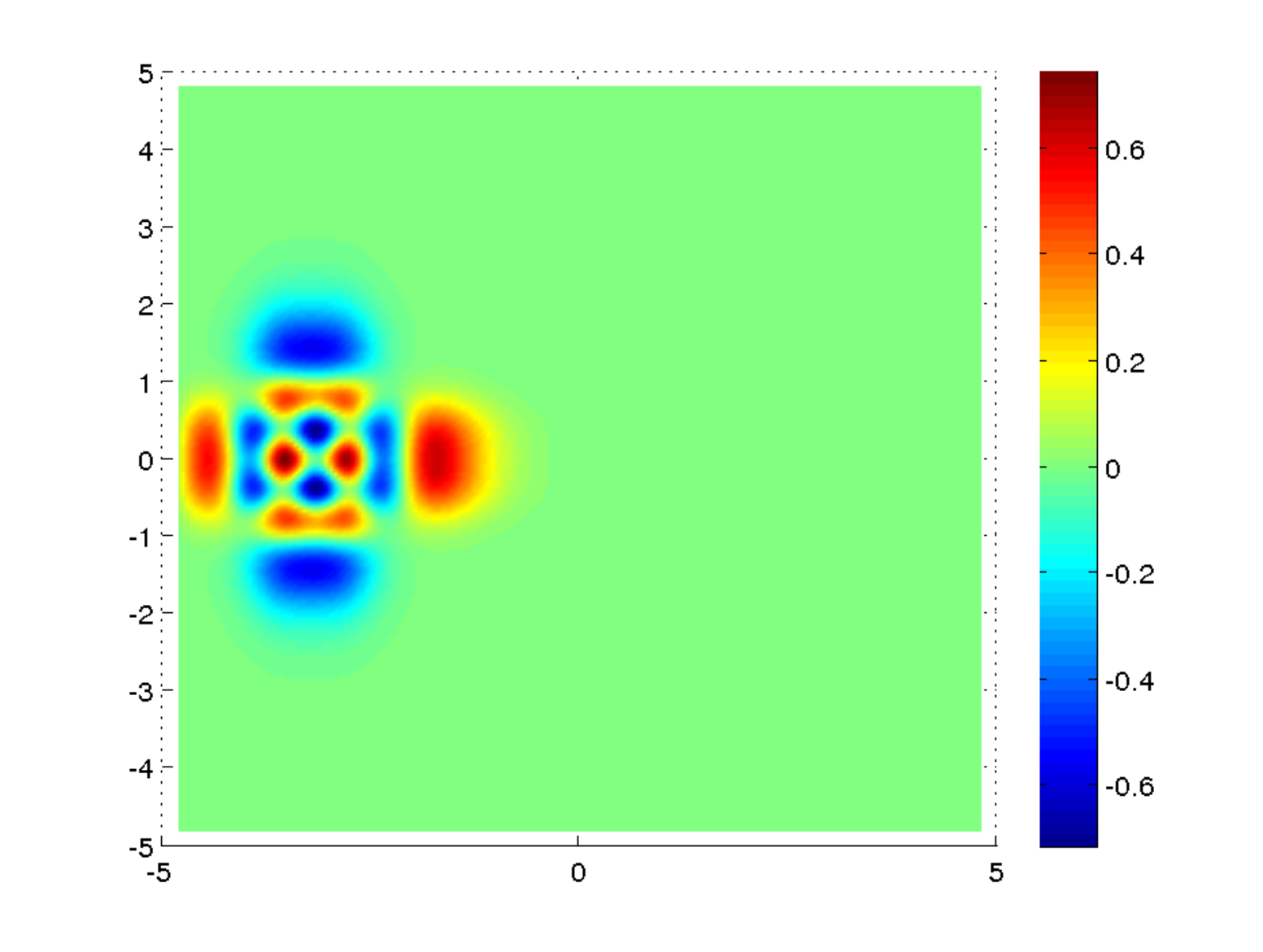}
\caption{First $6$ LSD's in $\Omega_{14}=[a_4-\epsilon^{(x_1)}/2,b_4+\epsilon^{(x_1)}/2]\times[c_3-\epsilon^{(x_2)}/2,d_3+\epsilon^{(x_2)}/2]$}
\label{SDBF_4x3}
\end{center}
\end{figure}
We also report in Table \ref{table1} the $L^{2}$-inner products in $[a,b]$, of the first 4 LSD's in $\Omega_{13}$ denoted here $v_{1},\cdots v_4$, with themselves and with the first 2 LSD's in $\Omega_{14}$ denoted $w_{1},w_2$, in order to illustrate their approximate orthogonality. Notice that the product of the $v_i$ with $w_j$ are naturally dependent on the overlap size as well as the mollifier parameters and the parameters used to build the localized 1-electron orbitals.
\begin{table}
\caption{LSD approximate orthogonality: $\langle v_i,v_j\rangle$ for all $i=1,\cdots,4$ and $j=i,\cdots 4$,  and $\langle v_i,w_j\rangle$ for $i=1,\cdots,4$ and $j=1,\cdots 2$}
\centering
\begin{tabular}{ccccccc}
Inner prod. & $v_1$ & $v_2$ & $v_3$ & $v_4$ & $w_1$ & $w_2$ \\
\hline
$v_1$ & $1$ & $2 \times 10^{-17}$ & $2\times 10^{-16}$ & $-2\times 10^{-17}$ & $-2\times 10^{-5}$ & $-5\times 10^{-5}$\\
$v_2$ & $2 \times 10^{-16}$ & $1$ & $-4\times 10^{-16}$ & $7\times 10^{-17}$ & $5\times 10^{-5}$ & $2\times 10^{-4}$\\
$v_3$ & $2\times 10^{-16}$  & $-4\times 10^{-16}$ & $1$ & $-5\times 10^{-16}$  & $-2\times 10^{-4}$ & $-7\times 10^{-4}$ \\
$v_4$ & $-2\times 10^{-17}$  &  $7\times 10^{-17}$ & $-5\times 10^{-16}$ & $1$ & $-1\times 10^{-3}$ & $3\times 10^{-3}$ \\
\hline
\end{tabular}
\label{table1}
\end{table}
We deduce from this table that, in general, the matrices ${\bf A}_i$ defined in Subsection \ref{expcons} are not exactly the identity matrices. 
\subsection{Test 2.a: Imaginary time experiment I}
\begin{figure}[!ht]
\begin{center}
\hspace*{1mm}\includegraphics[height=6cm, keepaspectratio]{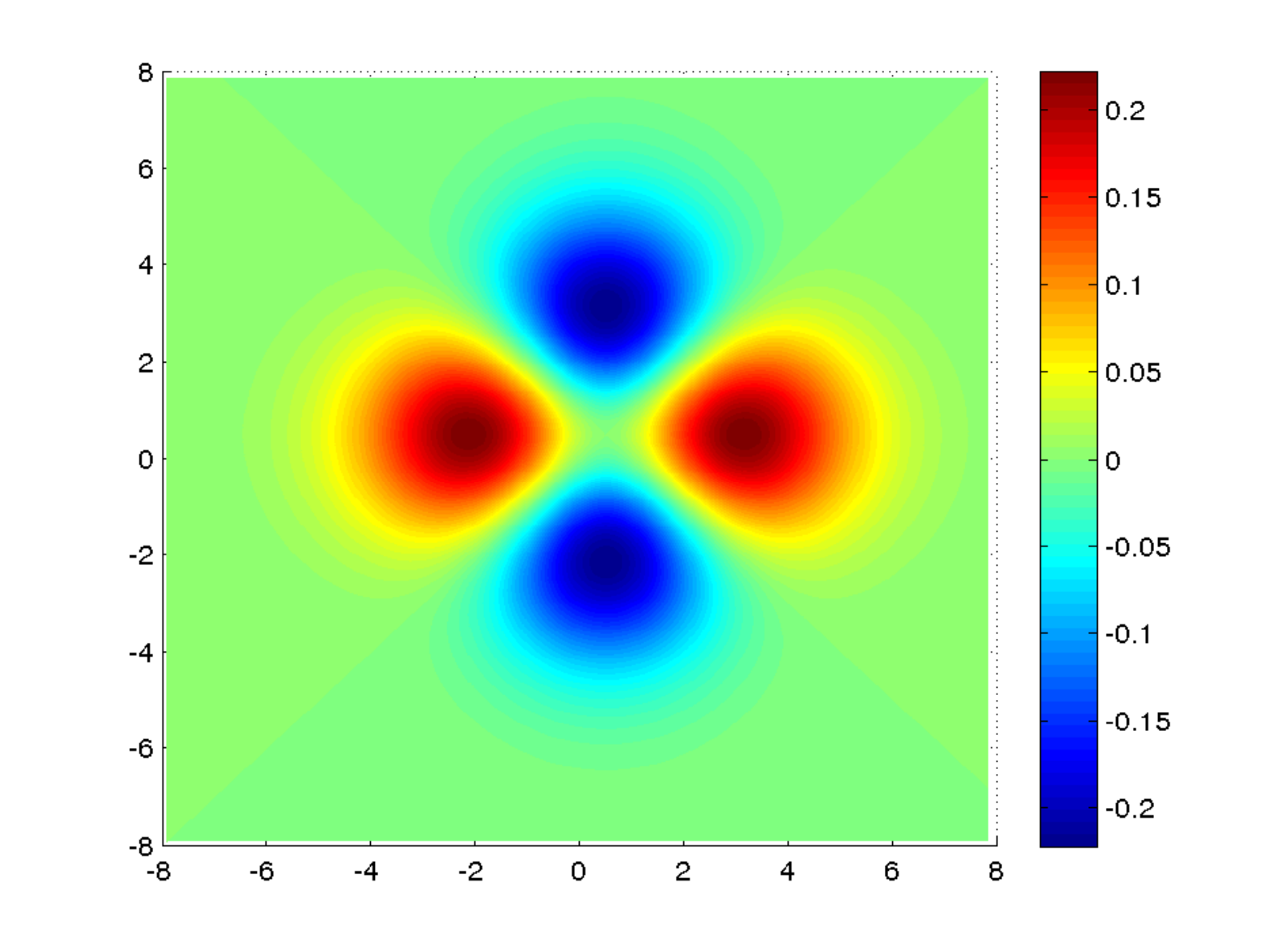}
\hspace*{1mm}\includegraphics[height=6cm, keepaspectratio]{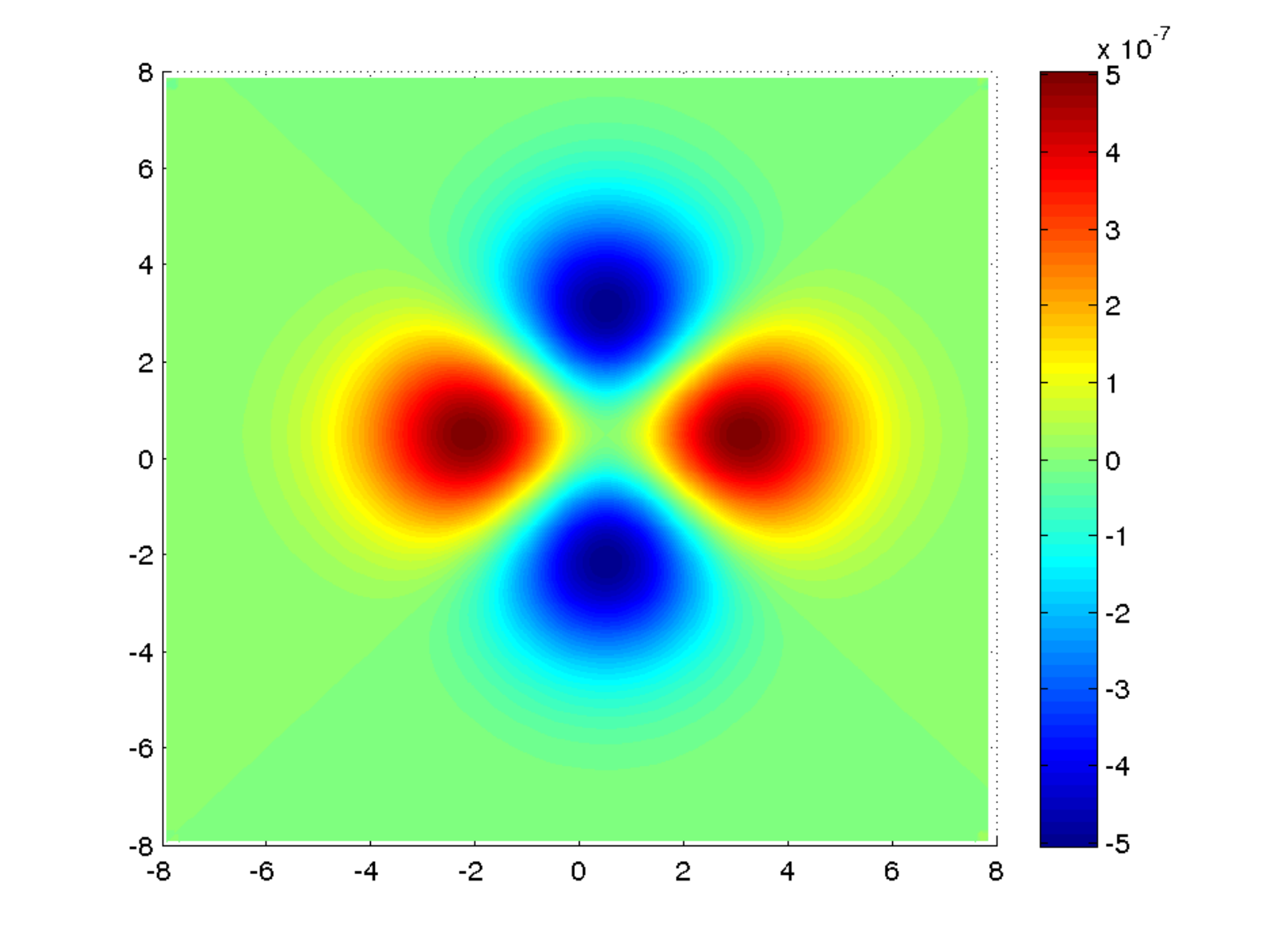}
\caption{(Left) Initial guess. (Right) Initial data reconstruction error: $\phi^{(k)}(x_1,x_2,0) - \phi_0(x_1,x_2)$.}
\label{phi0}
\end{center}
\end{figure}

\begin{figure}[!ht]
\begin{center}
\hspace*{1mm}\includegraphics[height=4cm, keepaspectratio]{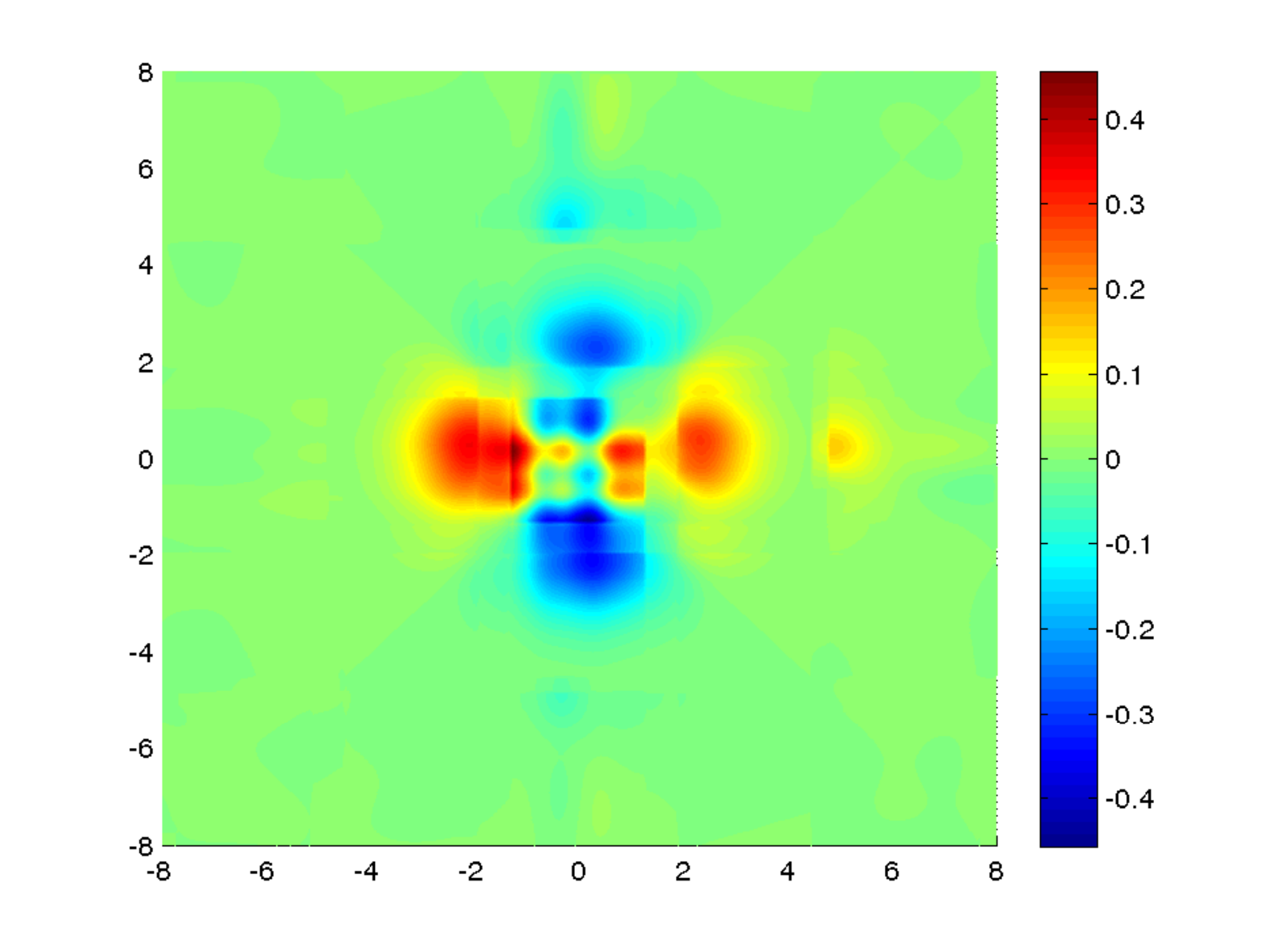}
\hspace*{1mm}\includegraphics[height=4cm, keepaspectratio]{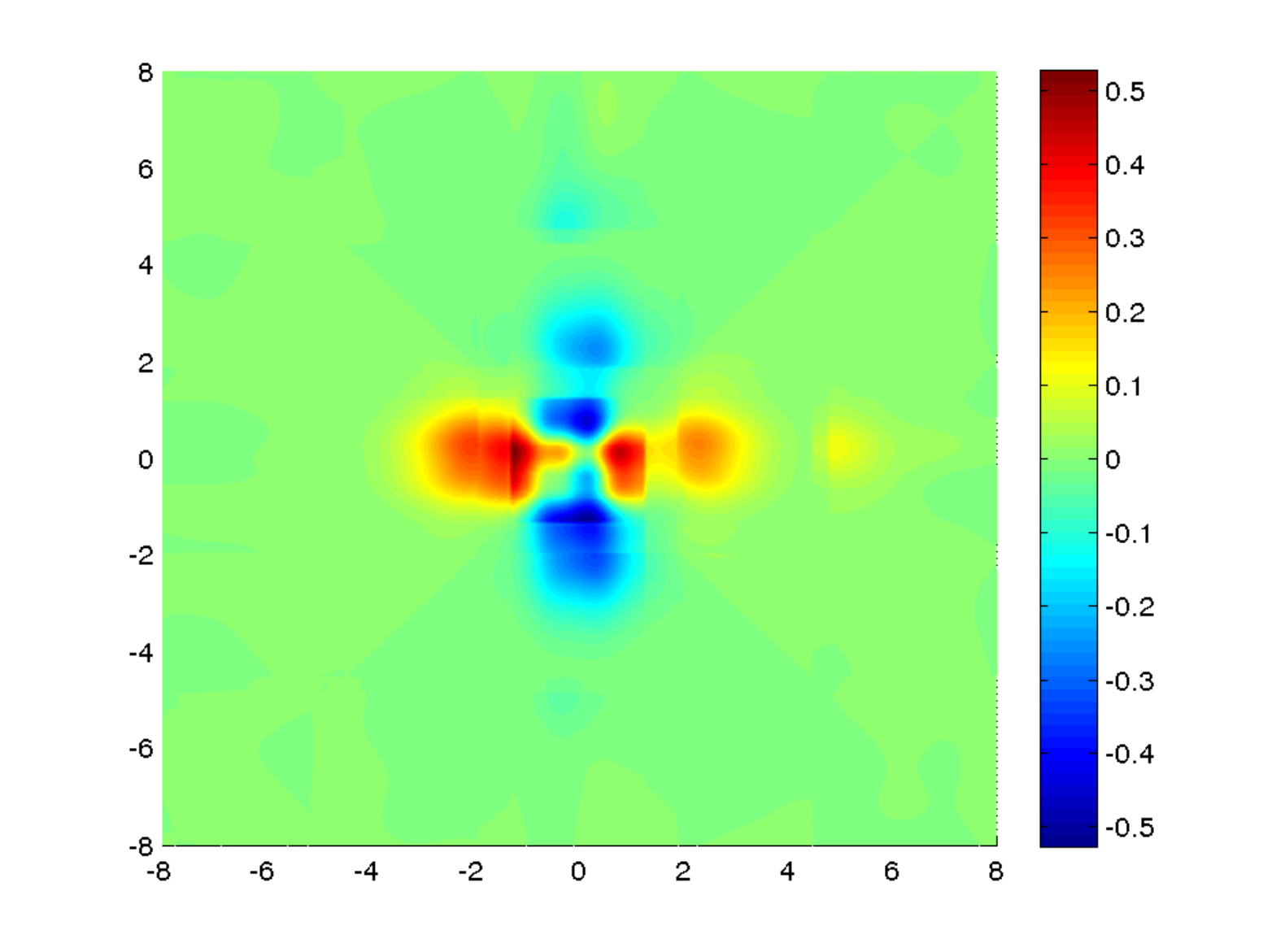}
\hspace*{1mm}\includegraphics[height=4cm, keepaspectratio]{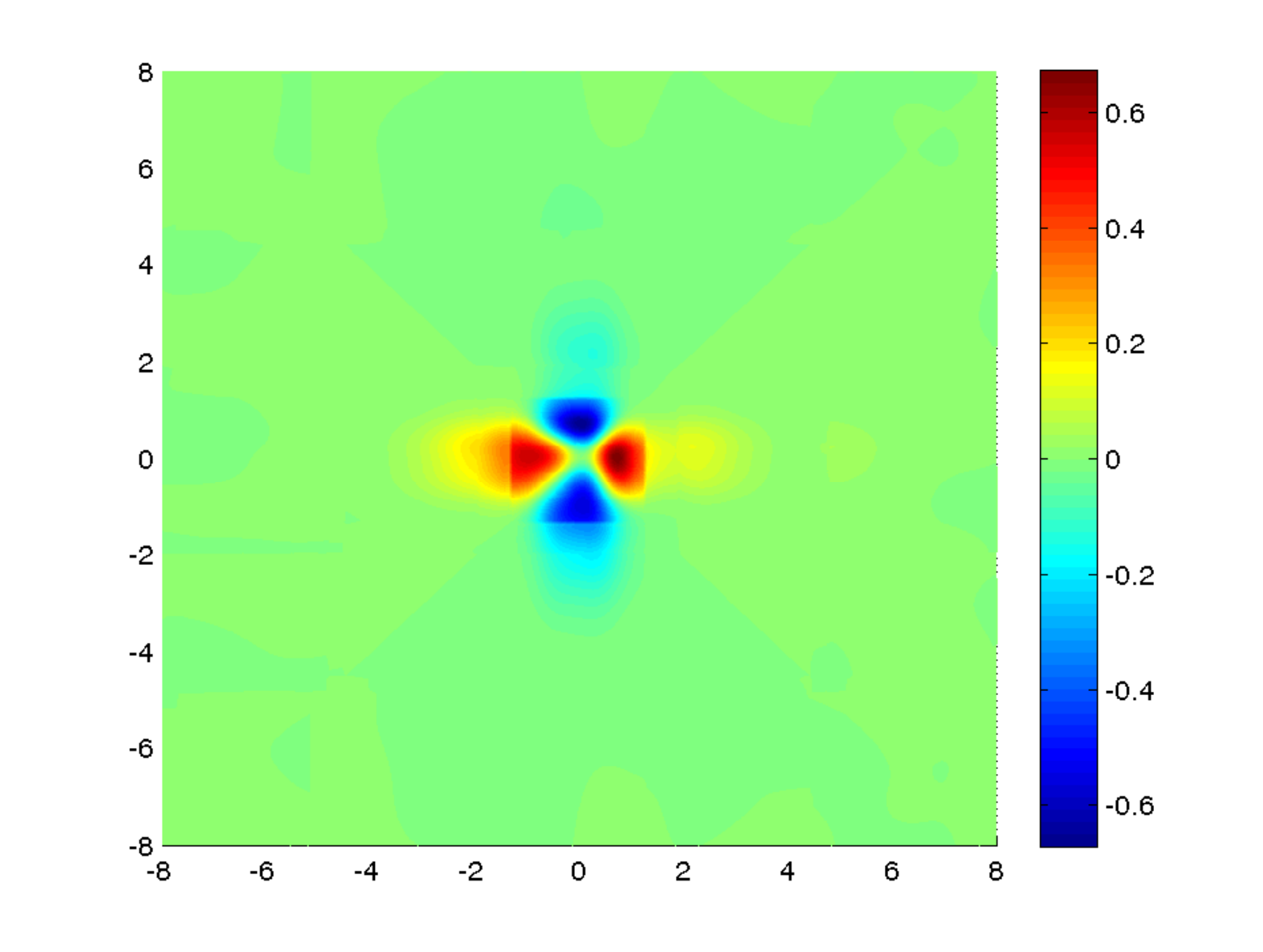}
\hspace*{1mm}\includegraphics[height=4cm, keepaspectratio]{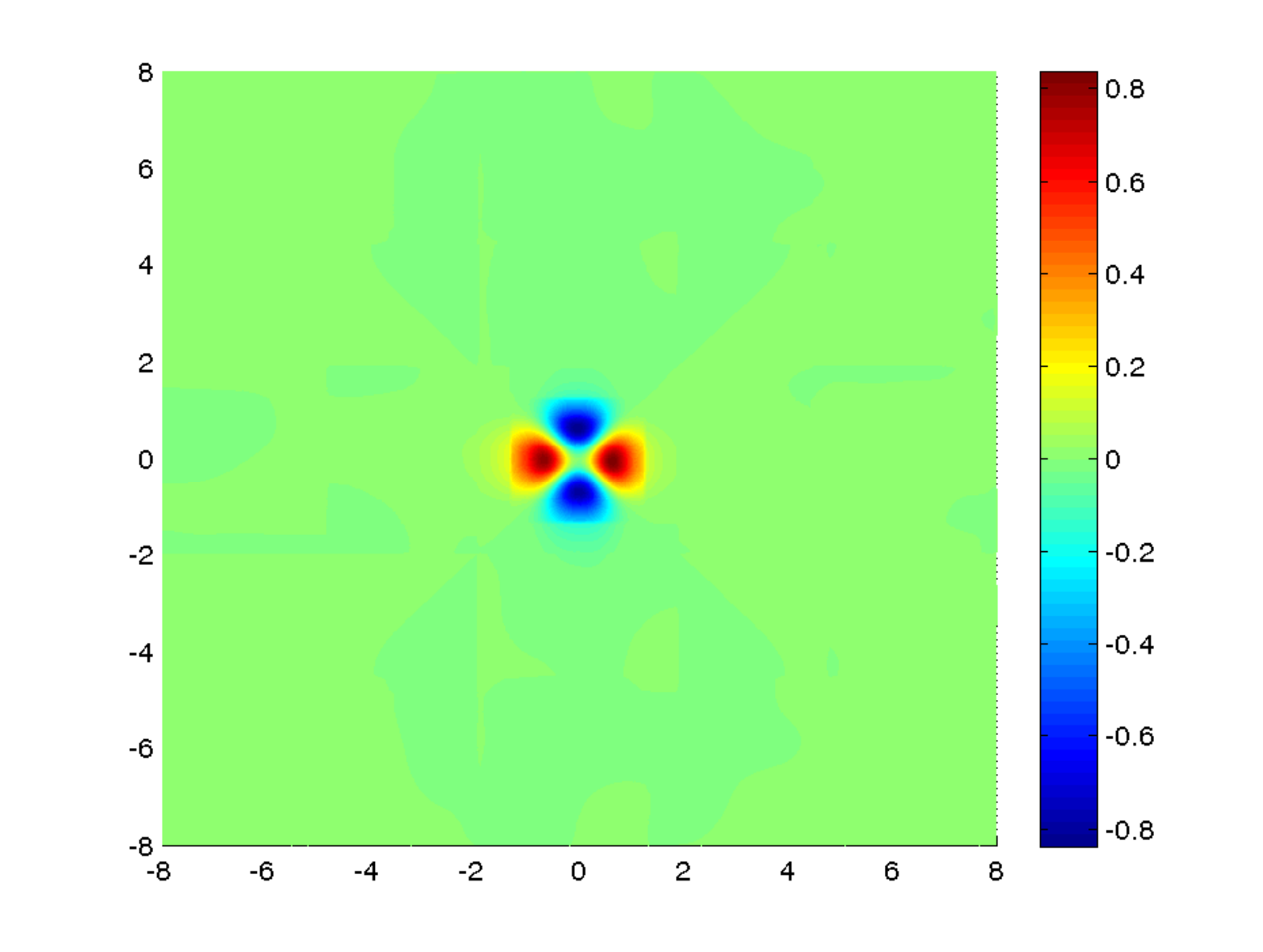}
\hspace*{1mm}\includegraphics[height=4cm, keepaspectratio]{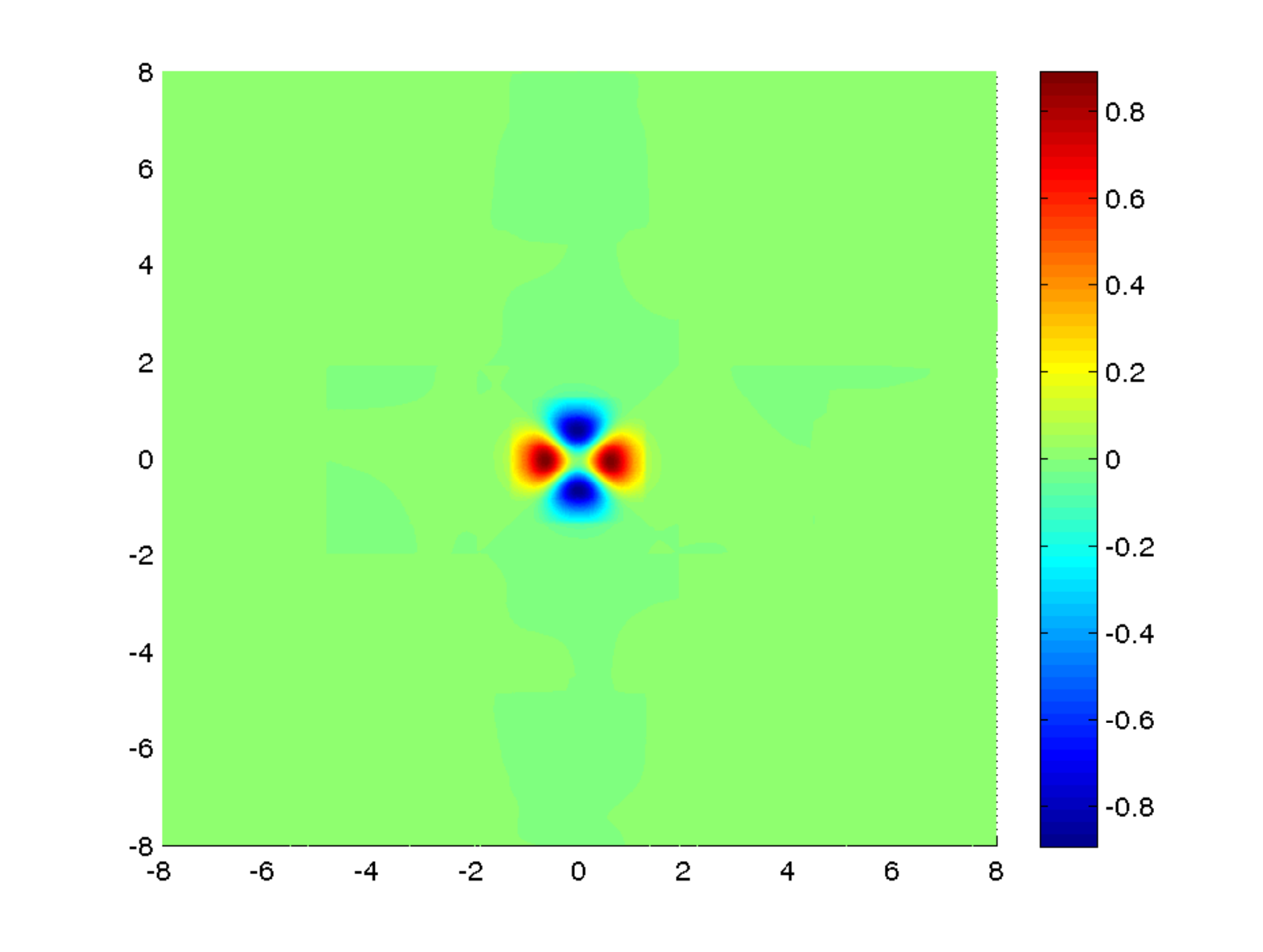}
\hspace*{1mm}\includegraphics[height=4cm, keepaspectratio]{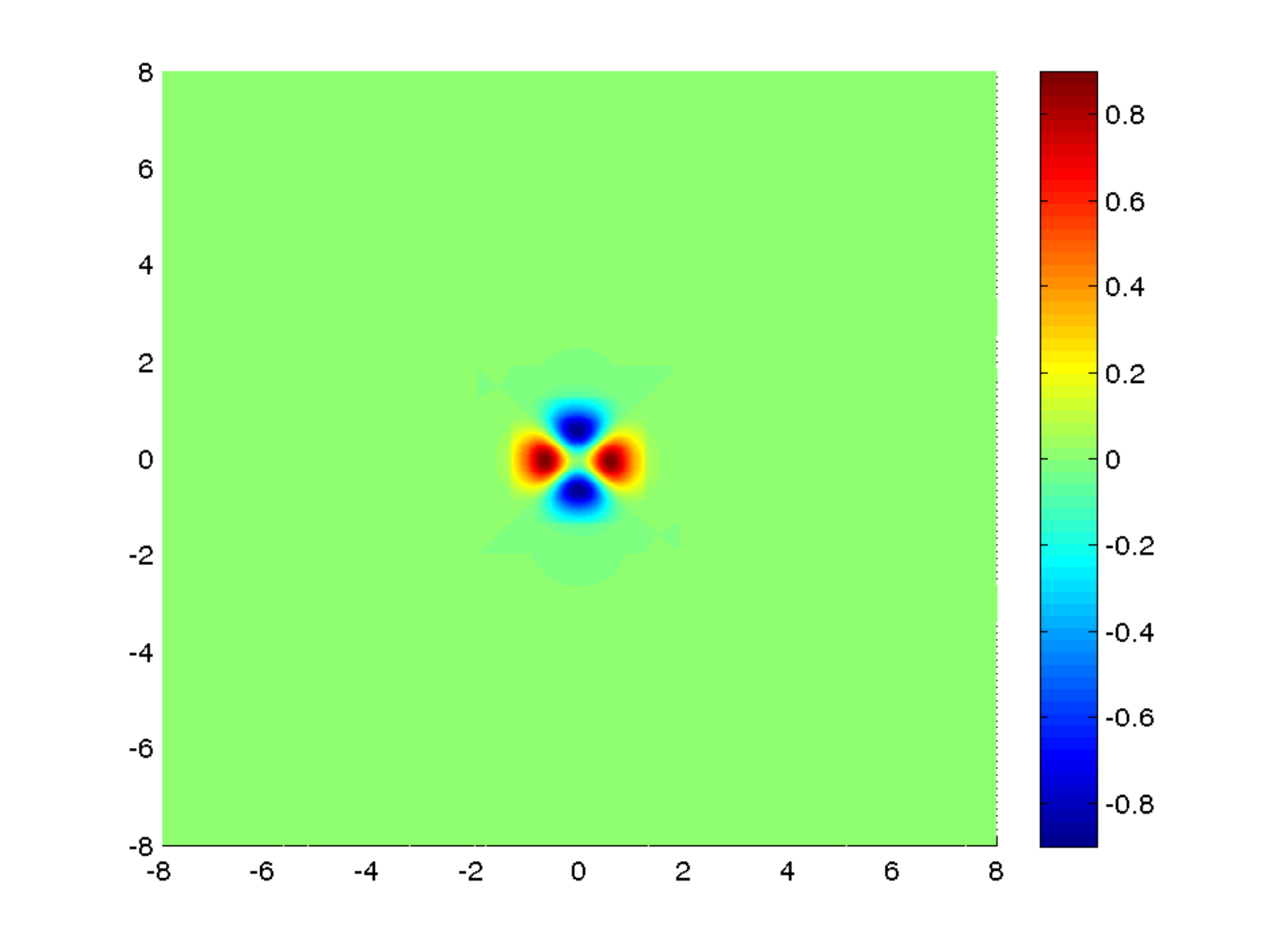}
\hspace*{1mm}\includegraphics[height=4cm, keepaspectratio]{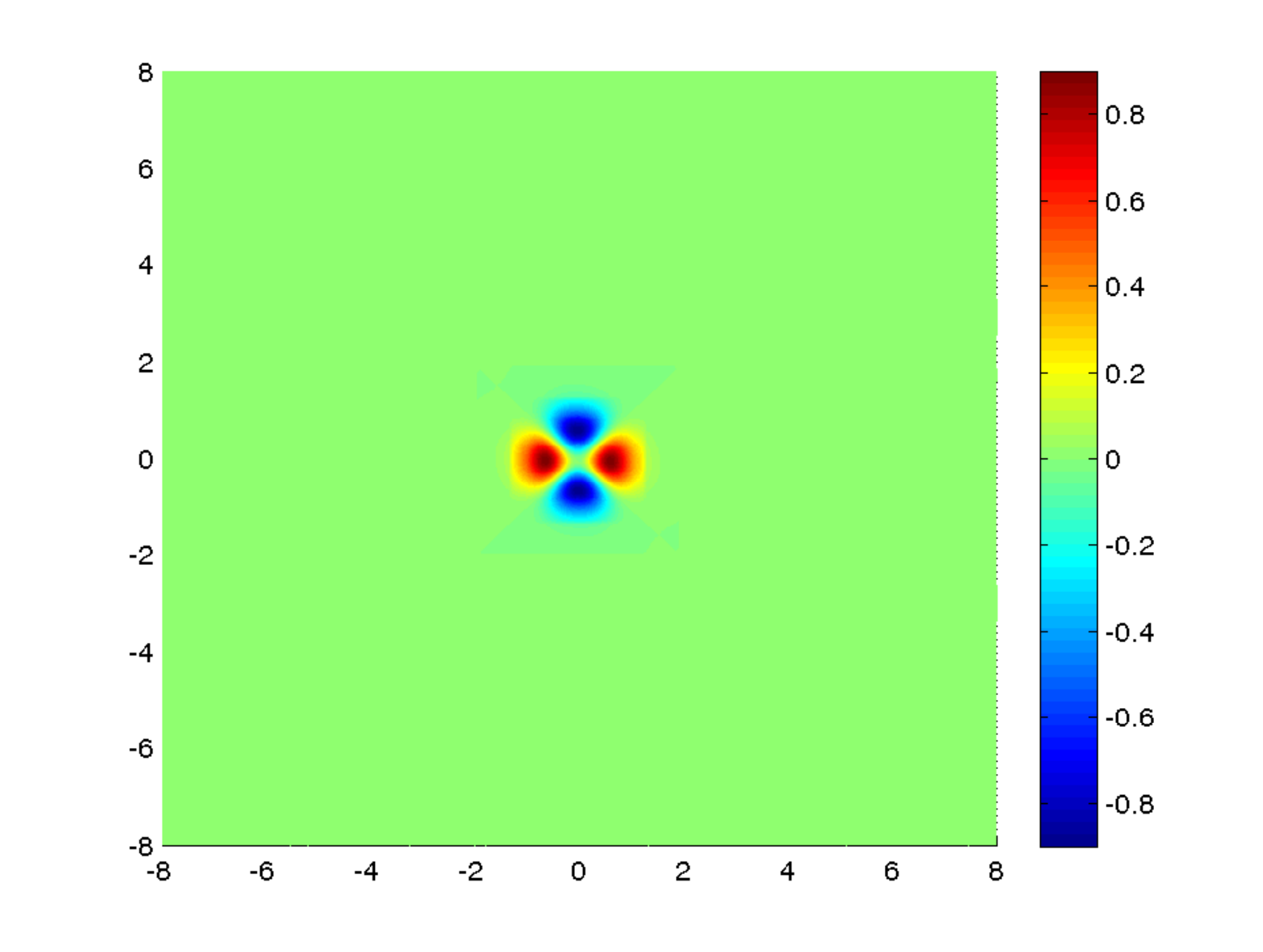}
\hspace*{1mm}\includegraphics[height=4cm, keepaspectratio]{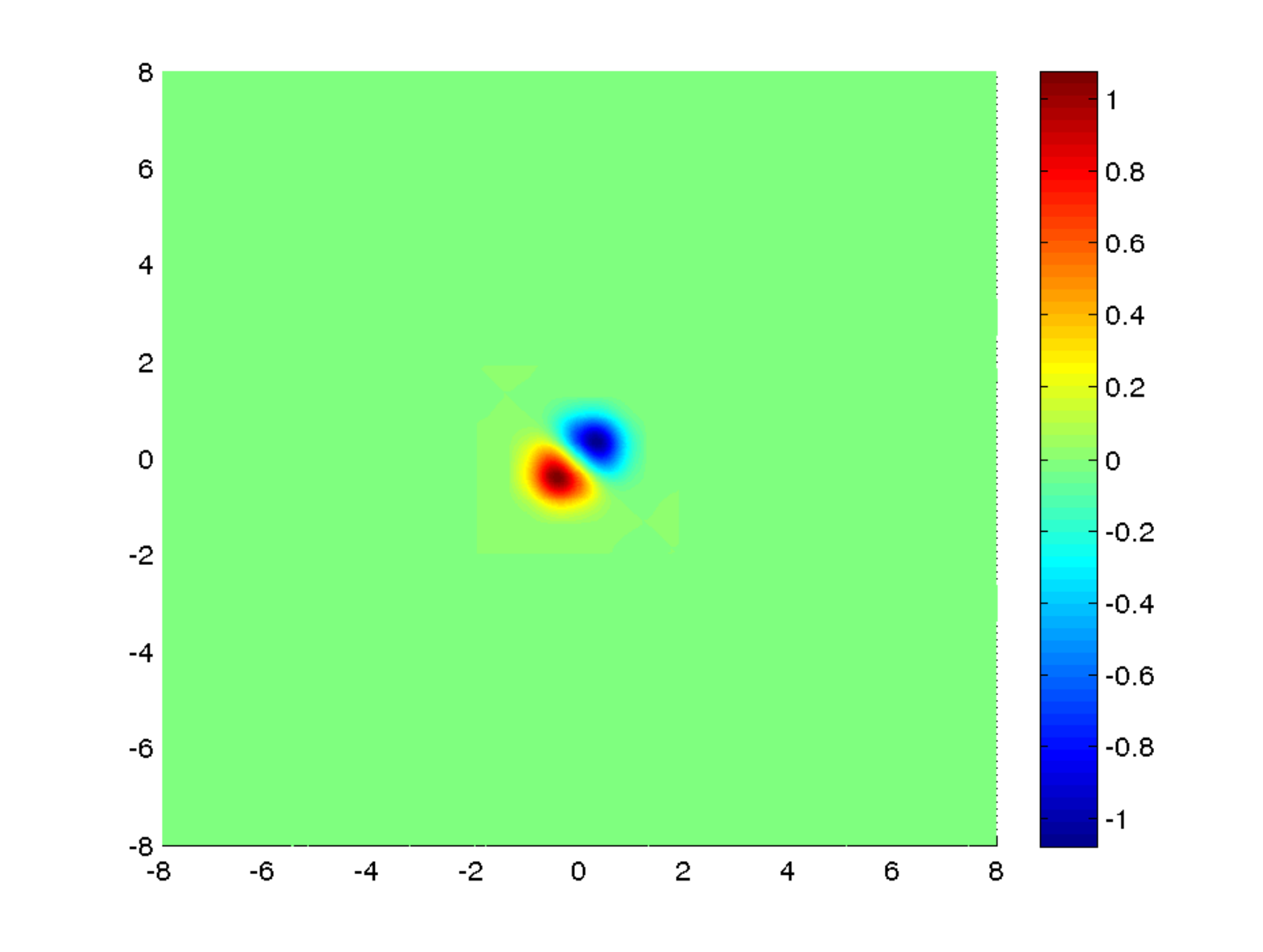}
\caption{Reconstructed wavefunction $\phi^{(0)}(x_1,x_2,t)$ with $t=n\Delta t$ with  $n=10,20,40,80,160,320,640,1280$.}
\label{phi0t}
\end{center}
\end{figure}

\begin{figure}[!ht]
\begin{center}
\hspace*{1mm}\includegraphics[height=6cm, keepaspectratio]{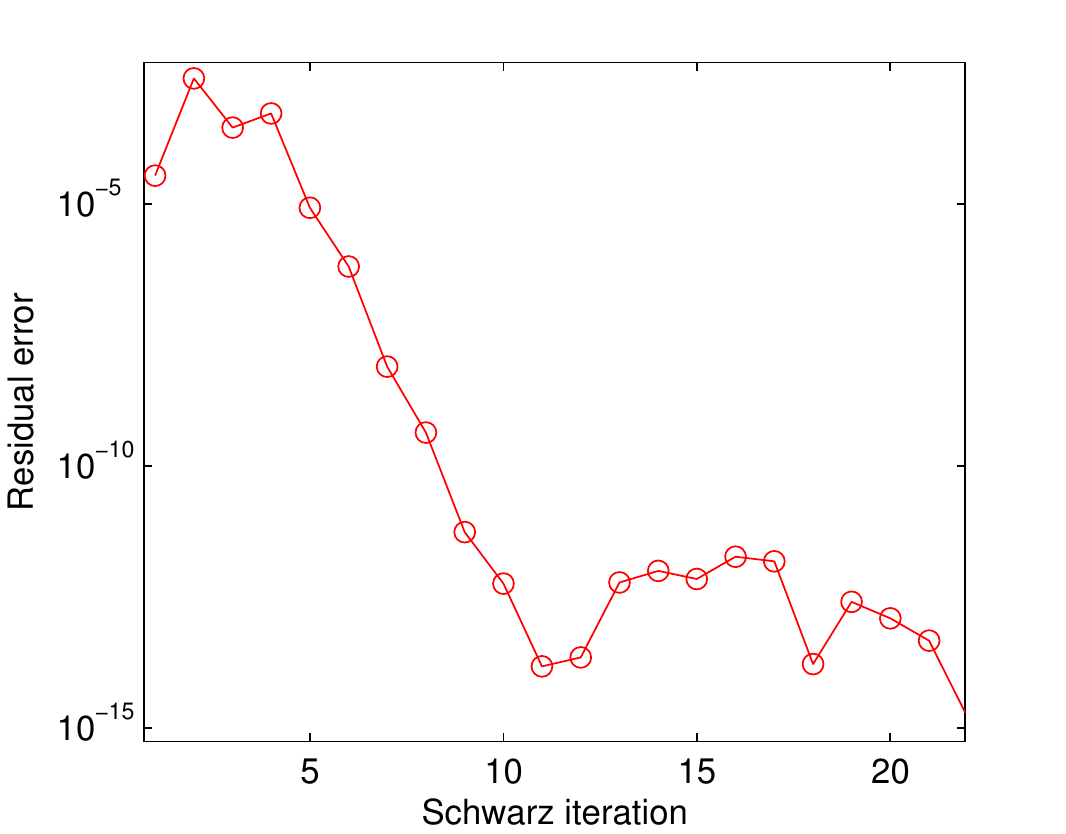}
\hspace*{1mm}\includegraphics[height=6cm, keepaspectratio]{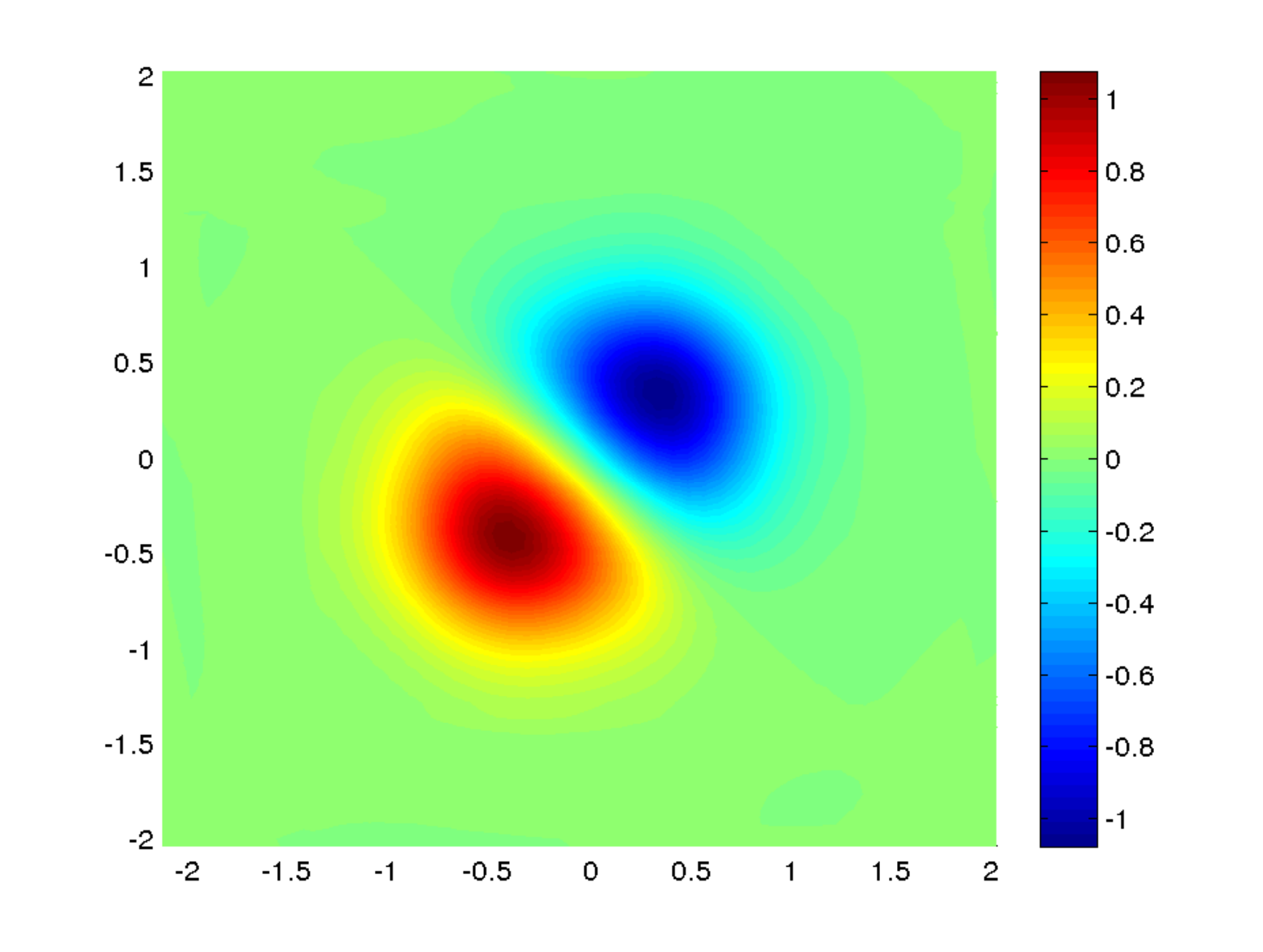}
\caption{(Left) Residual error in logscale as a function of Schwarz iteration. (Right) Converged ground state $\phi^{n^{\textrm{(cvg)}},k^\textrm{(cvg)}}$.}
\label{CV_IT}
\end{center}
\end{figure}

We implement the imaginary time method within a SWR domain decomposition framework in order to determine the ground state of $H_2$ with $x_A=-0.5$, $x_B=0.5$ (and $Z_A=Z_B=1$). The numerical data are as follows: $L=5$ (for a total of $25$ subdomains), $K=N_{\phi}(N_{\phi}+1)/2=28$ LSD's per subdomain (for a total of $25\times 28 = 700$ LSD's),  $a=c=-8$, $b=d=8$ and $\epsilon^{(x_{1,2})}=3.2\times 10^{-1}$.  The 2 nuclei are located in the central subdomain $\Omega_{13}=[a_3-\epsilon^{(x_1)}/2,b_3+\epsilon^{(x_2)}/2]\times[c_3-\epsilon^{(x_2)}/2,d_3+\epsilon^{(x_2)}/2]$ with $a_3-\epsilon^{(x_1)}/2=c_3-\epsilon^{(x_2)}/2=1.6$. We choose $\epsilon=0.1$ and $\sigma_1=315/256$  in \eqref{Beps}, and $N^{(x_1)}=N^{(x_2)}=201$. The LSD's are computed as in Subsection \ref{test1}. The matrices $\widetilde{{\bf H}}_i$ and $\widetilde{{\bf A}}_i$ are both sparse, and belong to $M_{28}(\R)$ for $i=1,\cdots,L^2$. \\
We choose Robin transmission conditions with $\mu=1$ in \eqref{RTC}, in the SWR algorithm. The initial guess is an antisymmetric function $\phi_0(x_1,x_2)=\widetilde{\phi}_0(x_1,x_2)/\|\widetilde{\phi}_0\|_{0}$ where
\begin{eqnarray*}
\left.
\begin{array}{lcl}
\widetilde{\phi}_0(x_1,x_2) & =  & \exp\big(-(x_1-1/2)^2/10-(x_2-1/2)^2/5\big)\\
& & - \exp\big(-(x_2-1/2)^2/10-(x_1-1/2)^2/5\big)
\end{array}
\right.
\end{eqnarray*}
and is represented in Fig. \ref{phi0} (Left).
Before iterating in Schwarz and imaginary time, we first need to construct the local projections of $\phi_0$ onto $\Omega_{i}=[a_i-\epsilon^{(x_1)}/2,b_i+\epsilon^{(x_1)}/2]\times[c_i-\epsilon^{(x_2)}/2,d_i+\epsilon^{(x_2)}/2]$, for each $i=1,\cdots,L$, with LSD's $\big\{v_j^{i}\big\}_{1\leq j \leq K}$ ($K=28$, $L=5$), that is we compute
\begin{eqnarray*}
\phi_i^{(0)}(x_1,x_2)  = \sum_{j=1}^K\langle \phi_0,v_j^{i}\rangle v^i_{j}(x_1,x_2).
\end{eqnarray*}
We reconstruct the initial data $\phi^{(0)}(\cdot,0)$ according to the algorithm presented in Subsection \ref{testA}, and we report in Fig. \ref{phi0} (Right) the reconstruction error: $\phi^{(0)}(\cdot,0) - \phi_0$. We then report in Fig. \ref{phi0t}, the reconstructed solution computed at Schwarz iteration $k=0$ and imaginary time $t_n=n\Delta t$, with $n=10,20,40,80,160,320,640,1280$ and with $\Delta t=1.4\times 10^{-3}$. As we can see, global convergence is almost reached at CNFG convergence at the first Schwarz iteration. This is due to the fact that the nuclei are located at the center of the central subdomain. In the subsequent Schwarz iterations, the residual error is still decreasing.
We then represent in Fig. \ref{CV_IT} (Left), as a function of the Schwarz iteration $k$, the residual error Res($k$) defined in \eqref{residueIT}. We also represent the converged solution $\phi^{\textrm{(cvg)}}$ in Fig. \ref{CV_IT} (Right), which was then reconstructed from $25$-local Schr\"odinger equations, showing the rapid convergence of the SWR algorithm.
\subsection{Test 2.b: Imaginary time experiment II}
\begin{figure}[!ht]
\begin{center}
\hspace*{1mm}\includegraphics[height=6cm, keepaspectratio]{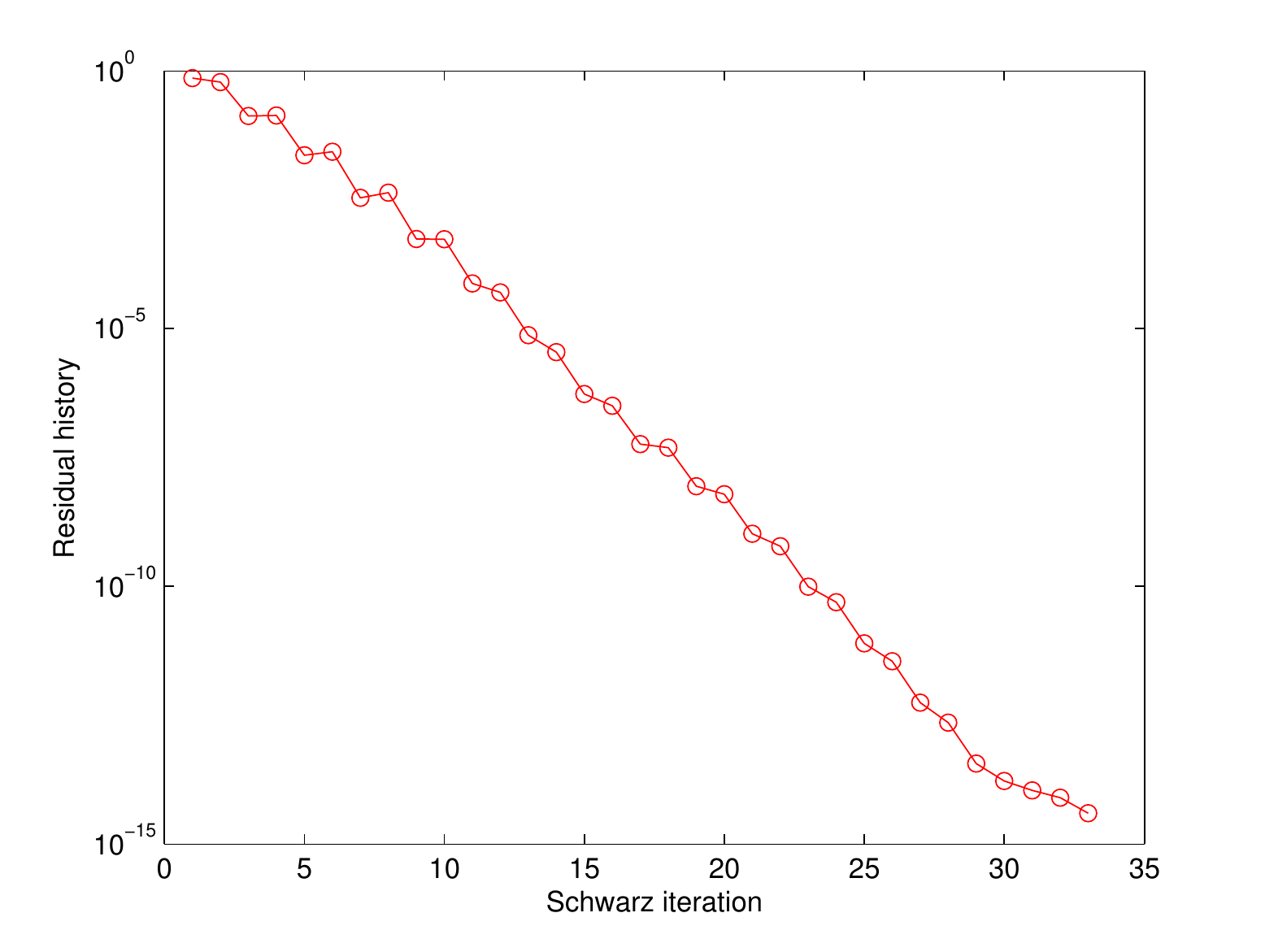}
\hspace*{1mm}\includegraphics[height=6cm, keepaspectratio]{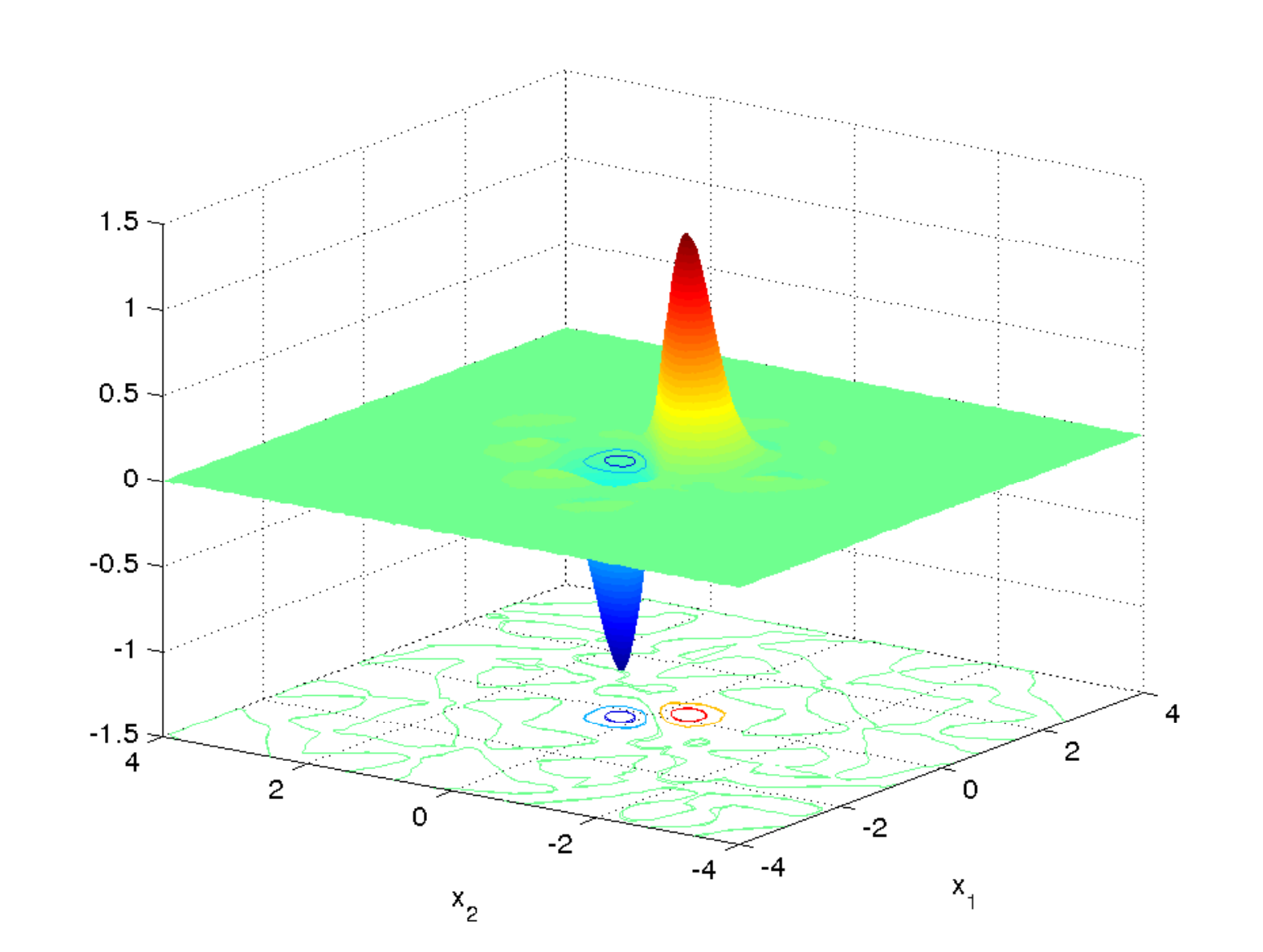}
\caption{(Left)  Residual error in logscale as a function of Schwarz iteration. (Right) Converged ground state $\phi^{n^{\textrm{(cvg)}},k^\textrm{(cvg)}}$}
\label{CV2_IT}
\end{center}
\end{figure}

We now present a stiffer problem. In the above test, the support of the ground state was mainly included in the central subdomain. In this new test, the support of the ground state of $H_2$ covers several subdomains. We take $x_A=-0.1$, $x_B=0.1$ and $\eta=0.6$ in \eqref{pseudo}. The basis is augmented by Gaussian basis functions as described in Subsection \ref{subsec:CB}, with $\delta=3$ in \eqref{gauss1d}, at the boundary of the subdomains, in order to ensure a better transmission. The global domain is $(-4,4)$ and is decomposed in $25$ subdomains. We select only $N_{\phi}=15$ local Slater's determinants per subdomain. The grid possesses a total of $N^{(x_1)}\times N^{(x_2)}=101^2$ points. Robin transmission conditions are imposed between subdomain with $\mu=10$, see \eqref{Heat1}.  The initial guess is an antisymmetric function $\phi_0(x_1,x_2)=\widetilde{\phi}_0(x_1,x_2)/\|\widetilde{\phi}_0\|_{0}$ where
\begin{eqnarray*}
\left.
\begin{array}{lcl}
\widetilde{\phi}_0(x_1,x_2) & =  & \exp\big(-4(x_1^2+x_2^2)\big).
\end{array}
\right.
\end{eqnarray*}
We first reconstruct the initial data $\phi^{(0)}(\cdot,0)$ as described in Subsection \ref{testA}. We report in Fig. \ref{CV2_IT} (Left) as a function of the Schwarz iteration $k$, the residual error Res($k$) in logscale. We also represent the converged solution $\phi^{\textrm{(cvg)}}$ in Fig. \ref{CV2_IT} (Right) which was then reconstructed from $25$-local Schr\"odinger equations, showing the rapid convergence of the SWR algorithm.
\section{Conclusion}\label{conclusion}
This paper was devoted to the derivation of a domain decomposition method for solving the $N$-body Schr\"odinger equation. More specifically a Schwarz waveform relaxation algorithm with Robin transmission conditions was proposed, along with a pseudospectral method, for computing in parallel ``many'' local Schr\"odinger equations, and from which is constructed a global wavefunction.  In order to improve the accuracy while keeping efficiency, local Slater's determinant functions were alternatively selected as basis functions, allowing in principle, for reduced local bases and as a consequence lower dimensional local approximate Hamiltonians. Some preliminary tests show a promising approach, which will be further developed on more elaborated cases.
\\
\\
\noindent{\bf Acknowledgments}. The author would like to thank Prof. C.R. Anderson (UCLA) for helpful discussions about mollifiers and grid-based methods for solving the N-body Schr\"odinger equation.

\newpage
\appendix
\section{Antisymmetric wavefunction reconstruction}\label{APXA}
\noindent In order to construct an antisymmetric wavefunction, we propose a specific decomposition of $\R^{dN}$ with {\it antisymmetric} local basis functions. In order to simplify the presentation, we will assume that i) the subdomains $\{\Omega_i\}_{i \in \{1,\cdots,L^{dN}\}}\in \R^{dN}$ are hypercubes of identical size and $L\in 2\N+1$, and ii) there is no overlap between the subdomains. We first define:
\begin{definition}
 We denote by $\sigma(i;p,q) \in \{1,\cdots,L^{dN}\}$ the subdomain index such that for $({\bf x}_1,\cdots,{\bf x}_p,\cdots,{\bf x}_q,\cdots,{\bf x}_N) \in \Omega_{i}$: $({\bf x}_1,\cdots,{\bf x}_q,\cdots,{\bf x}_p,\cdots,{\bf x}_N) \in \Omega_{\sigma(i;p,q)}$. The index $i$ refers to the subdomain index, and $(p,q)$ to the permutation coordinate indices.
\end{definition}
Notice that we naturally have $\sigma\big(\sigma(i;p,q);q,p\big)=i$ and that $\sigma\big(\sigma(i;p,q);q,p\big)$ is unique as there is no subdomain overlap. 
Antisymmetry of the wavefunction would then occur if for all $i \in\{1,\cdots,L^{dN}\}$, the local basis $\big\{v^{i}_j\big\}_{1\leq j\leq K_i}$ coincides with the local bases $\big\{v^{\sigma(i;p,q)}_j\big\}_{1\leq j \leq K_{\sigma(i;p,q)}}$, for all $(p,q) \in\{1,\cdots,N\}^2$.  We define
\begin{eqnarray*}
\Sigma(i)=\big\{\sigma(i;p,q)\in \{1,\cdots,L^{dN}\}, \, \forall (p,q)\in \{1,\cdots,N\}^2\big\}
\end{eqnarray*}
If there is no overlap between subdomains, $\Sigma(i)$ is actually a singleton. In order to guarantee the antisymmetry of the overall wavefunction, at any time and any Schwarz iteration, we proceed as follows. \\
\noindent{\it Antisymmetric wavefunction}. For any $i\in \{1,\cdots,L^{dN}\}$ and for any $l \in \Sigma(i)$:
\begin{itemize}
\item we solve a (local) TDSE on $\Omega_i$ in the form
\begin{eqnarray}\label{statWF}
\psi^{(k)}_i(\cdot,t) = \sum_{j=1}^{K_i}c_j^{i,(k)}(t)v^i_{j}
\end{eqnarray}
\item and we deduce:
\begin{eqnarray}\label{statWF2}
\psi^{(k)}_{\Sigma(i)}(\cdot,t) = -\sum_{j=1}^{K_{\Sigma(i)}}c_j^{i,(k)}(t)v^{\Sigma(i)}_{j}
\end{eqnarray}
\end{itemize}
The global wavefunction is then reconstructed according to the algorithm derived in Section \ref{NAWF} without overlap. We deduce by construction, the following proposition.
\begin{prop}
The reconstructed wavefunction $\psi^{(k)}$ with \eqref{statWF} and \eqref{statWF2}, is antisymmetric.
\end{prop} 
\noindent{\bf Proof.} The proof is trivial. For any $({\bf x}_1,\cdots,{\bf x}_p,\cdots,{\bf x}_q,\cdots,{\bf x}_N) \in \Omega_i$, \\
we have $({\bf x}_1,\cdots,{\bf x}_q,\cdots,{\bf x}_p,\cdots,{\bf x}_N)\in \Omega_{\Sigma(i)}$, then 
\begin{eqnarray*}
\left.
\begin{array}{lcl}
\psi^{(k)}({\bf x}_1,\cdots,{\bf x}_p, \cdots,{\bf x}_q,\cdots, {\bf x}_N,t) & =& \psi_i^{(k)}({\bf x}_1,\cdots,{\bf x}_p, \cdots,{\bf x}_q,\cdots, {\bf x}_N,t)\\
\\
&  = & -\psi_{\Sigma(i)}^{(k)}({\bf x}_1,\cdots,{\bf x}_q, \cdots,{\bf x}_p,\cdots, {\bf x}_N,t)\\
\\
&  =&  -\psi^{(k)}({\bf x}_1,\cdots,{\bf x}_q, \cdots,{\bf x}_p,\cdots, {\bf x}_N,t).
\end{array}
\right.
\end{eqnarray*}
 $\Box$
\\
\\
{\bf Discrete local Hamiltonian construction with local Slater's determinants}. The construction of discrete local Hamiltonians $\widetilde{{\bf H}}_{i}$ is similar to the procedure described in Subsection \ref{subsec:SLO}, which is an application of \cite{CAM15-10}, in a DDM framework. Due to the antisymmetry constraint detailed above, we however need additional specifications. The strategy presented in Subsection \ref{subsec:SLO}, and extended to $d$ dimensions and $N$ particles allows for the construction for any subdomain $\Omega_i$, $1\leq i \leq L^{dN}$, of the SLO's $\big\{\varphi_j^{i}({\bf x})\big\}_{j}$. These SLO's are smooth, have compact support and possess orthogonality properties, which  are described at the end of Subsection \ref{subsec:SLO}. We here summarize the explicit construction of a local Hamiltonian $\widetilde{{\bf H}}_i$, $1 \leq i \leq L^{dN}$, say for $\Omega_i$. Following the notations used above, we need to compute:
\begin{eqnarray*}
\left.
\begin{array}{lcl}
I_{mpqr}^{ijkl} &= &\int_{\R^{dN}} \varphi^i_m({\bf x})\varphi^j_p({\bf x})\cfrac{\varphi^k_q({\bf s})\varphi^l_r({\bf s})}{|{\bf s}-{\bf x}|}d{\bf x}d{\bf s} \\
&= &\int_{\hbox{Supp}\varphi_m^i\cap \hbox{Supp}\varphi_p^j\cap \hbox{Supp}\varphi_q^k\cap \hbox{Supp}\varphi_r^l}\varphi^m_j({\bf x})\varphi^j_p({\bf x})\Phi_{qr}^{kl}({\bf x})d{\bf x}.
\end{array}
\right.
\end{eqnarray*} 
where 
\begin{eqnarray*}
\Phi_{qr}^{kl}({\bf x}) = \int_{\hbox{Supp}\varphi^k_q\cap \hbox{Supp}\varphi^l_r}\cfrac{\varphi_q^k({\bf s})\varphi_r^l({\bf s})}{|{\bf s}-{\bf x}|}d{\bf s}
\end{eqnarray*} 
can be achieved using mollifiers:
\begin{eqnarray*}
\Phi_{qr}^{kl}({\bf x}) \approx \int_{\hbox{Supp}\varphi^k_q\cap \hbox{Supp}\varphi^l_r}\varphi_q^k({\bf s})\varphi_r^l({\bf s})B_{\epsilon}({\bf s}-{\bf x})d{\bf s}
\end{eqnarray*} 
or alternatively, if $d=3$
\begin{eqnarray*}
-\triangle \Phi_{qr}^{kl} = 4\pi\varphi_q^k\varphi_r^l
\end{eqnarray*}
Then, for $m=l,l\pm 1$
\begin{eqnarray*}
I_{lm}^{ij} = \int_{\hbox{Supp}\varphi^i_l\cap \hbox{Supp}\varphi_m^j}\cfrac{1}{2}\varphi_l^i({\bf x})\triangle\varphi_m^j({\bf x})-\sum_{A=1}^P\cfrac{Z_A}{|{\bf x}-{\bf x}_A|}\varphi_l^i({\bf x})\varphi_m^j({\bf x})d{\bf x}.
\end{eqnarray*} 
Notice that for constructing  ${\bf Q}_i^{x,y,z}$, we need to compute
\begin{eqnarray*}
J^{ij}_{lm} = \int_{\hbox{Supp}\varphi^i_l\cap \hbox{Supp}\varphi_m^j}{\bf x}\varphi_l^i({\bf x})\varphi_m^{j}({\bf x})d{\bf x}
\end{eqnarray*}
which does not present any additional difficulty compared to one-domain problems, see again \cite{CAM15-09}.  We then efficiently construct the local Hamiltonians using the strategy presented in \cite{CAM15-10}. Once the matrix local Hamiltonians are constructed, we can determine the LSD's and solve the time-independent and dependent Schr\"odinger equations.

\section{Computational complexity analysis of the SWR-DDM solver for $N$-body equation}\label{APXC}
\noindent In this appendix, we analyze the computational complexity of the SWR method for both the time independent and time dependent cases. An obvious consequence of the use of a SWR$/$FCI method with orbital basis functions, is that the number of degrees of freedom (dof) is expected to be much smaller compared to finite difference$/$volume methods (FDM$/$FVM) or low degree finite element methods (FEM). For instance, in a subdomain $\Omega_i$, a cell center FVM ($Q_0$) consists of choosing the basis functions as $\big\{{\bf 1}_{V^i_j}({\bf x}_1,\cdots,{\bf x}_N)/|V^i_j|\big\}_{1\leq j\leq N_i}$ with finite volumes $V^i_j$, such that $\cup_{j=1}^{N_i}V_j^i=\tau_h(\Omega_i)$. That is:
\begin{eqnarray*}
\psi_i^{(k)}({\bf x}_1,\cdots,{\bf x}_N,t) = \sum_{j=1}^{N_i}\cfrac{1}{|V^i_j|}_j{\bf 1}_{V_j^i}c^i_j(t), \, ({\bf x}_1,\cdots,{\bf x}_N) \in \Omega_i.
\end{eqnarray*}
Although these $Q_0$-basis functions are very simple, and that the corresponding matrices are sparse, in order to get a precise description of the wavefunction a ``very'' large number $N_i$ on $\Omega_i$, of finite volumes is necessary.  Slater's orbitals $\big\{v_j^i\big\}_{1\leq j\leq K_i}$, defined above, would typically contain a very large number of finite volumes $V_j^i$ ($N_i \gg K_i$). The consequence is that, although FVM-linear systems are much sparser than Galerkin-FCI systems, they are also of much higher dimension.\\
Below, we study the overall computational complexity (CC) and the scalability of the of the SWR-DDM in $d$-dimension and for $N$ particles. The analysis will be provided for both the stationary and unstationary $N$-body Schr\"odinger equations. \\
\\
\noindent{\bf Computational complexity and scalability for the time-independent $N$-body equation.} From now on, we assume that the computational domain is decomposed in $L^{dN}$-subdomains, $\{\Omega_{p}\}_{1\leq p \leq L^{dN}}$, and that $K_{p}$ (resp. $\mathcal{K}_p$) with $1 \leq p \leq L^{dN}$, local basis functions are selected (resp. the number of degrees of freedom) per subdomain. We denote by $K_{\textrm{Tot}}:=\sum_{p=1}^{dN} K_{p}$ (resp. $\mathcal{K}_{\textrm{Tot}}:=\sum_{p=1}^{dN} \mathcal{K}_{p}$) the total number of local basis functions (resp. degrees of freedom). In the following, we will assume for simplicity that there is a fixed number of local basis functions per subdomain, that is $\mathcal{K}_{p} \approx \mathcal{K}_{\textrm{Tot}}/L^{dN}$, for any $1 \leq p \leq L^{dN}$. Instead of dealing with a full discrete Hamiltonian in $M_{\mathcal{K}_{\textrm{Tot}}}(\R)$, we then rather deal with $L^{dN}$ local discrete Hamiltonians in $M_{\mathcal{K}_{\textrm{Tot}}/L^{dN}}(\R)$. The very first step then consists of constructing the local basis functions, then of the $L^{dN}$ local Hamiltonians, which is pleasingly parallel (perfect distribution of the integral computations). We focus on the complexity and scalability for computing the eigenenergies from these (local or global) discrete Hamiltonians. \\
\\
 We then recall the main ingredients necessary to study the computational complexity analysis of the SWR method presented in Section \ref{SWR} for solving the time-independent Schr\"odinger equation using the NGF-method. We have decomposed the spatial domain $\Omega \subset \R^{dN}$ in $L^{dN}$ subdomains and solve an imaginary-time-dependent Schr\"odinger equation (or real-time normalized heat equation) on each subdomain. At a given Schwarz iteration $k \in \N$, we denote by $T^{(k)}_{p}$ (resp. $n_{p}^{(k)}$) the imaginary convergence time (resp. number of time iterations to converge) for the NGF-method in the subdomain $\Omega_{p}$, where $1\leq p \leq L^{dN}$. In addition, each imaginary time iteration requires $\mathcal{O}(\mathcal{K}_{p}^{\beta^{\textrm{(S)}}_p})$ operations, where $1 < \beta^{\textrm{(S)}}_p <  3$ (due to sparse linear system solver). The index $(\textrm{S})$ refers to the stationary Schr\"odinger equation. We denote by $k^{\textrm{(cvg)}}$, the total number of Schwarz iterations to reach convergence, as described in Section \ref{SWR}. Notice that $k^{\textrm{(cvg)}}$ is  strongly dependent on the type of transmission conditions \cite{lorin-TBS2}. We get
\begin{prop} 
The computational complexity CC$^{\textrm{\textrm{(S)}}}_{\textrm{SWR}}$ of the overall SWR-DDM method describe in Subsection \ref{SWR2} for solving the Schr\"odinger equation in the stationary case is given by
\begin{eqnarray}\label{CC2} 
\textrm{CC}^{\textrm{(S)}}_{\textrm{SWR}} = \mathcal{O}\Big(\sum_{k=1}^{k^{\textrm{(cvg)}}}\sum_{p=1}^{L^{dN}}n_{p}^{(k)}\mathcal{K}_p^{\beta^{\textrm{(S)}}_p}\Big).
\end{eqnarray}
Assuming that $\beta^{\textrm{(S)}}_p$, (resp. $n_{p}^{(k)}$) is $p$-independent (that is subdomain independent),  and then denoted $\beta^{\textrm{(S)}}$ (resp. $N^{(k)}$), we have
\begin{eqnarray*}
\textrm{CC}^{\textrm{(S)}}_{\textrm{SWR}} = \mathcal{O}\Big(\cfrac{\mathcal{K}_{\textrm{Tot}}^{\beta^{\textrm{(S)}}}}{L^{dN(\beta^{\textrm{(S)}}-1)}}\sum_{k=1}^{k^{\textrm{(cvg)}}} N^{(k)}\Big).
\end{eqnarray*}
\end{prop}
Thus
\begin{itemize}
\item assuming that the algorithm is implemented on a $P$-core machine, each core will deal with $\approx L^{dN}/P$ subdomains. The message passing load is dependent on the type of transmission conditions, but typically for classical or Robin SWR the communication load will be very light. As a consequence an efficiency ($T/PT_P$) close to $1$ is expected.
\item we are dealing with $L^{dN}$ linear systems with approximatly $\mathcal{K}_{\textrm{Tot}}/L^{dN}$ degrees of freedom (instead of a unique large system of $\mathcal{K}_{\textrm{Tot}}$ degrees of freedom if a huge discrete Hamiltonian was considered). As $\beta^{\textrm{(S)}}$ is strictly greater than $1$, we benefit from a scaling effect. 
\end{itemize}
The SWR-DDM approach will be attractive in the starionary case, if typically
\begin{eqnarray*}
\sum_{k=1}^{k^{\textrm{(cvg)}}} N^{(k)} \ll L^{dN(\beta^{\textrm{(S)}}-1)}
\end{eqnarray*}
In order to satisfy this condition i) an implicit solver for the heat equation will allow for a faster convergence (the bigger the time step, the smaller $N^{(k)}$) of the NGF method, appropriate transmission conditions will allow for a minimization of $k^{\textrm{(cvg)}}$.\\
\\
\noindent{\bf Computational complexity and scalability for the time-dependent $N$-body equation.} Thanks to the SWR approach and as in the stationary case, the computation of the time-dependent equation, does not involve a full discrete Hamiltonian in $M_{\mathcal{K}_{\textrm{Tot}}}(\C)$, but rather $L^{dN}$ discrete local Hamiltonians in $M_{\mathcal{K}_{\textrm{Tot}}/L^{dN}}(\C)$. Notice that if we use the same Gaussian basis functions in each subdomain, the local potential-free Hamiltonians are identical in each subdomain, and has to be performed only once. The contribution from the interaction potential and laser field in the local Hamitonians, are however subdomain-dependent, but are diagonal operators. We assume that the TDSE is computed from time $0$ to $T>0$. An implicit scheme ($L^2$-norm preserving) is implemented, which necessitates the numerical solution at each time iteration of a sparse linear system. We denote by $n_T$ the total number of time iterations to reach $T$, which will be assumed to be the same for both methods. We deduce that
\begin{prop}
The computational complexity $\textrm{CC}^{\textrm{(NS)}}_{\textrm{SWR}}$ of the overall SWR-DDM method describe in Subsection \ref{SWR1} for solving the Schr\"odinger equation in the time dependent case is given by:
\begin{eqnarray}\label{CC4} 
\textrm{CC}^{\textrm{(NS)}}_{\textrm{SWR}} = \mathcal{O}\big(n_Tk^{\textrm{(cvg)}}\sum_{p=1}^{L^{dN}}\mathcal{K}_p^{\beta^{\textrm{(NS)}}_p}\big).
\end{eqnarray}
with $1< \beta_p^{\textrm{(NS)}} <  3$, where the index $\textrm{(NS)}$ refers to the non-stationary case. Assuming that the $\beta^{\textrm{(NS)}}_p$  is $p$-independent (that is subdomain independent),  and denoted $\beta^{\textrm{(NS)}}$, we get
\begin{eqnarray*}
\textrm{CC}^{\textrm{(NS)}}_{\textrm{SWR}} = \mathcal{O}\Big(n_Tk^{\textrm{(cvg)}}\cfrac{\mathcal{K}_{\textrm{Tot}}^{\beta^{\textrm{(NS)}}}}{L^{dN(\beta^{\textrm{(NS)}}-1)}}\Big).
\end{eqnarray*}
\end{prop}
The SWR-DDM will then be attractive in the non-stationary case, if 
\begin{eqnarray*}
k^{\textrm{(cvg)}} \ll L^{dN(\beta^{\textrm{(NS)}}-1)} \, .
\end{eqnarray*}

\bibliographystyle{plain}
\bibliography{refs}

\end{document}